\documentclass[11pt,reqno,twoside]{article}
\def\COMPLETE{}
\input{package_header} %file containing all used packages

\begin{document}
%------------------------------------------------
\title{A Stochastic Bregman Primal-Dual Splitting Algorithm for Composite Optimization}
\author{Antonio Silveti-Falls\thanks{Toulouse School of Economics, University of Toulouse, France. E-mail: Tony.S.Falls@gmail.com. This work was done while the author was at GREYC-ENSICAEN.} \and Cesare Molinari
\thanks{Istituto Italiano di Tecnologia, Italy. E-mail: cecio.molinari@gmail.com} \and  
Jalal Fadili\thanks{Normandie Universit\'e, ENSICAEN, UNICAEN, CNRS, GREYC, France. E-mail: Jalal.Fadili@ensicaen.fr.}
}
\date{}
\maketitle
\begin{flushleft}\end{flushleft}
%%%%%%%%%%%%%%%%%%%%%%%%%%%%%%%%%%%
\begin{abstract}
We study a stochastic first order primal-dual method for solving convex-concave saddle point problems over real reflexive Banach spaces using Bregman divergences and relative smoothness assumptions, in which we allow for stochastic error in the computation of gradient terms within the algorithm. We show ergodic convergence in expectation of the Lagrangian optimality gap with a rate of $O\para{1/k}$ and that every almost sure weak cluster point of the ergodic sequence is a saddle point in expectation under mild assumptions. Under slightly stricter assumptions, we show almost sure weak convergence of the pointwise iterates to a saddle point. Under a relative strong convexity assumption on the objective functions and a total convexity assumption on the entropies of the Bregman divergences, we establish almost sure strong convergence of the pointwise iterates to a saddle point. Our framework is general and does not need strong convexity of the entropies inducing the Bregman divergences in the algorithm. Numerical applications are considered including entropically regularized Wasserstein barycenter problems and regularized inverse problems on the simplex.
\end{abstract}

\begin{keywords}
Bregman divergence; primal-dual splitting; noneuclidean splitting; saddle point problems; first order algorithms; convergence rates; relative smoothness; total convexity; Banach space.
\end{keywords}

\begin{AMS}
49J52, 65K05, 65K10.
\end{AMS}

%\tableofcontents

%%%%%%%%%%%%%%%%%%%%%%%%%%%%%%%%%%%
%----------------------------------
% These commands allow Gummi to
% edit each piece individually
% and still maintain previews.
%----------------------------------
\ifdefined\COMPLETE
\else
\documentclass[12pt]{article}
\input{tex_package_header} %file containing all the used libraries
\begin{document}
\fi
%----------------------------------
% Remove at your own risk!
% There is also a footer that ends
% the document if the main wasn't
% loaded.
%----------------------------------

\section{Introduction}
\begin{subsection}{Problem Statement and Algorithm}
The goal is to solve the following primal-dual, or saddle point, problem over the real reflexive Banach spaces $\X_p$ and $\X_d$, where the subscript $p$ refers to primal and $d$ to dual:
\begin{equation}\tag{$\mathscr{P}.\mathscr{D}.$}\label{probl}
\begin{split}
& \min\limits_{x\in\X_p} \max\limits_{\mu\in\X_d} \ \LL{x,\mu}
\end{split}
\end{equation}
where 
\begin{equation}\label{lagrangian}
\begin{split}
\LL{x,\mu}\eqdef f(x)+g(x)+\iota_{\C_p}(x)+\scal{Tx}{\mu}-h^{*}(\mu)-l^{*}(\mu)-\iota_{\C_d}(\mu)
\end{split}
\end{equation}
is the Lagrangian functional and $\iota_{\C_p}$ and  $\iota_{\C_d}$ are the indicator functions of the convex constraint sets $\C_p$ and $\C_d$, respectively, and $T: \X_p \to \X_d^*$ is a linear mapping. We denote the primal and dual problems as
\begin{equation}\label{PProb}\tag{$\mathrsfs{P}$}
\min\limits_{x\in\C_p} \brac{f(x)+g(x)+\para{\left(h \underset{\C_d}{\infconv} l\right)\circ T} (x)}
\end{equation}
\begin{equation}\label{DProb}\tag{$\mathrsfs{D}$}
\min\limits_{\mu\in\C_d} \brac{h^{*}(\mu)+l^{*}(\mu)+\para{\left(f^{*} \underset{\C_p}{\infconv} g^{*}\right)\circ \left(-T^*\right)} (\mu)}
\end{equation}
where $\para{f^*\underset{\C_p}{\infconv}g^*}^* \eqdef f + g + \iota_{\C_p}$, using $*$ to denote the Fenchel conjugate, and similarly for $\underset{\C_d}{\infconv}$. In the case in which $\C_p$ and $\C_d$ are trivial constraints, i.e., the entire spaces $\X_p$ and $\X_d$, the corresponding primal and dual problems related to \eqref{probl} are
\begin{equation*}
\begin{split}
& \min\limits_{x\in\X_p} \brac{f(x)+g(x)+\para{\left(h \infconv l\right)\circ T} (x)}\\
& \min\limits_{\mu\in\X_d} \brac{h^{*}(\mu)+l^{*}(\mu)+\para{\left(f^{*} \infconv g^{*}\right)\circ \left(-T^*\right)} (\mu)}
\end{split}
\end{equation*}
where $\infconv$ recovers the classical \emph{infimal convolution} defined by $f \infconv g(v) = \inf_{w \in \X_p} \para{f(w)+g(v-w)}$. The set of solutions for \eqref{PProb} and \eqref{DProb} are written as
\nnewq{\label{eq:solset}
\solsetp&\eqdef \argmin\limits_{x\in\C_p}\brac{\max\limits_{\mu\in\C_d}\brac{f\para{x} + g\para{x} +\ip{Tx,\mu}{}-h^*\para{\mu}-l^*\para{\mu}}}\\\solsetd&\eqdef \argmax\limits_{\mu\in\C_d}\brac{\min\limits_{x\in\C_p}\brac{f\para{x}+g\para{x} +\ip{Tx,\mu}{}-h^*\para{\mu} - l^*\para{\mu}}}.
}
The set of saddle points for the Lagrangian defined in \eqref{lagrangian} is denoted
\newq{\label{sadp}
\sadp \eqdef \brac{\para{\xs,\mus}\in\X_p\times\X_d: \quad\forall \para{x,\mu}\in\X_p\times\X_d,\quad \LL{\xs,\mu}\leq \LL{\xs,\mus} \leq \LL{x,\mus}}
}
which obeys the inclusion $\sadp\subset\solsetp\times\solsetd$.

Given a real reflexive Banach space $\X$, we denote by $\Gamma_0\left(\X\right)$ the space of proper convex lower semicontinuous functions from $\X$ to $\R\cup \left\{+\infty\right\}$. For a subset $\C$ of a Banach space, $\inte \C$ denotes its interior. We suppose the following standing hypotheses on the problem, which we collectively denote by \ref{ass:hyp}:\leqnomode
\begin{align}
\label{ass:hyp}\tag*{\bf{($\mathbf{H}$)}}
\left\{
\begin{tabular}{p{.825\textwidth}}
\begin{enumerate}[label=\bf{(\subscript{\mathbf{H}}{{\arabic*}})}]
	\item The Banach spaces $\X_p$ and $\X_d$ are real and reflexive, while $\C_p\subset \X_p$ and $\C_d\subset \X_d$ are nonempty convex closed subsets.
	\item The functions $f$ and $g$ belong to $\Gamma_0\left(\X_p\right)$ while $l$ and $h$ belong to $\Gamma_0\left(\X_d\right)$, with $\C_p \subset \dom(f)$ (resp. $\C_d \subset \dom(h^*)$) and $f$ (resp. $h^*$) is differentiable on $\inte \C_p$ (resp. $\inte \C_d$).
	\item $\C_p \cap \dom(g) \neq \emptyset$ and $\C_d \cap \dom(l^*) \neq \emptyset$.
	\item The operator $T: \X_p\to\X_d^*$ is linear and continuous.
	\item The set of saddle points $\sadp$ for \eqref{probl} is nonempty. %\label{assum:saddle}
\end{enumerate}
\end{tabular}
\right.
\end{align}
\reqnomode
It is well-known that $\sadp$ is non-empty under suitable domain qualification conditions.

Before introducing the method, we recall the definition of Bregman divergence which will be key to our algorithm and to the theoretical analysis of convergence.
\begin{definition}[Bregman divergence]\label{def:bregman}
Given a function $\phi: \X \to \R\cup\brac{\pinfty}$, often referred to as the entropy, differentiable on $\inte \dom\para{\phi}$, its Bregman divergence is defined by
	\newq{
	D_{\phi}(x,y) \eqdef \begin{cases}\phi(x)-\phi(y)-\langle \nabla\phi(y),x-y\rangle & \mbox{if }x\in \dom\para{\phi}\qandq y\in\inte\dom\para{\phi},\\ +\infty & \mbox{else}.\end{cases}
	}
\end{definition}
Notice that if $\phi$ belongs to $\Gamma_0\left(\X\right)$ then $D_\phi$ is always nonnegative by the subdifferential inequality.\\

The Stochastic Bregman Primal-Dual Splitting algorithm, SBPD for short, is presented in \algref{alg:ibpds}. We introduce two entropy functions, $\phi_p$ and $\phi_d$, and we denote by $D_p$ and $D_d$ their Bregman divergences, respectively. We further consider the possibility of some stochastic error in the computation of the gradients\footnote{The addition of stochastic error in the computation of $D$-$\prox$ operators associated to $g$ or $l^*$, while interesting, is problematic for the algorithm in the sense that the monotone inclusions may no longer hold and the iterates themselves might not remain in the interior of the domain as desired.} $\nabla f$ and $\nabla h^*$ which we will denote for $\nabla f\para{x_k}$ as $\pnk$ and for $\nabla h^*\para{\mu_k}$ as $\dnk$. 

\begin{algorithm}[h]
	\caption{\label{alg:ibpds}
	\label{alg:bpds}Stochastic Bregman Primal-Dual Splitting (SBPD).}
	\For{$k=0,1,\ldots$}{
		\begin{equation*}
			\begin{split}
				x_{k+1}&=\argmin_{x\in\C_p} \left\{g(x)+\langle \nabla f (x_k) + \pnk,x\rangle +\langle Tx, \wtilde{\mu}_k\rangle +\frac{1}{\lambda_k}D_p\left(x,x_k\right)\right\}\\
				\mu_{k+1}&=\argmin_{\mu\in\C_d} \left\{l^{*}(\mu)+\langle \nabla h^{*} (\mu_k) + \dnk,\mu\rangle - \langle T\wtilde{x}_k, \mu\rangle +\frac{1}{\nu_k}D_d\left(\mu,\mu_k\right)\right\}
			\end{split}
		\end{equation*}
		where $\wtilde{\mu}_k=\mu_k$ and $\wtilde{x}_k=2x_{k+1}-x_k$.
	}
\end{algorithm}

In the deterministic setting for the primal update, i.e., $\pnk=0$ for each $k\in\N$, the first step of the algorithm can be re-written in the following way:
\begin{equation*}
	\begin{split}
		x_{k+1}&=\argmin_{x\in\C_p} \left\{g(x)+f (x_k)+\langle \nabla f (x_k),x-x_k\rangle +\langle Tx, \wtilde{\mu}_k\rangle +\frac{1}{\lambda_k}D_p\left(x,x_k\right)\right\}\\
		&=\para{\nabla \phi_p+\lambda_k\partial g }^{-1} \para{\nabla \phi_p-\lambda_k \nabla \left(f(\cdot)+ \langle T\cdot, \wtilde{\mu}_k\rangle\right)  } \ \left(x_k\right)\\
		&=\para{\nabla \phi_p+\lambda_k\partial g }^{-1} \ \left(\nabla \phi_p\left(x_k\right)-\lambda_k\nabla f \left(x_k\right)-\lambda_k T^* \wtilde{\mu}_k\right).
	\end{split} 
\end{equation*}
Analogously, if $\dnk =0$ for all $k\in\N$, 
\begin{equation*}
	\begin{split}
		\mu_{k+1}&=\para{\nabla \phi_d+\nu_k\partial l^{*}}^{-1} \ \left(\nabla \phi_d\left(\mu_k\right)-\nu_k\nabla h^* \left(\mu_k\right)+\nu_k T \wtilde{x}_k\right).
	\end{split}
\end{equation*}
A priori, the mappings $\para{\nabla \phi_p+\lambda_k\partial g}^{-1}$ and $\para{\nabla\phi_d+\nu_k\partial l^*}^{-1}$, sometimes referred to as $D$-$\prox$ mappings, may be empty, may not be single-valued, or may not map $\inte\dom \para{\phi_p}$ (resp. $\inte\dom\para{\phi_d}$) to $\inte\dom\para{\phi_p}$ (resp. $\inte\dom\para{\phi_d}$). In light of this, we will only consider $\phi_p$ and $\phi_d$ for which these mappings are well-defined and map from $\inte\dom\para{\phi_p}$ to $\inte\dom\para{\phi_p}$ and the analog for $\phi_d$ (see \ref{ass:legendre}). In Section~\ref{sec:assump}, we will elaborate on the class of Legendre functions on a real reflexive Banach space given in \cite[Definition 2.2]{BauschkeBorweinCombettes} which will help us to ensure that the $D$-$\prox$ mappings are well-defined.
\end{subsection}

\begin{subsection}{Contribution and Prior Work}
The idea of using primal-dual methods to solve convex-concave saddle point problems has been around since the 1960s, e.g., \cite{moreau64}, \cite{rockafellar64}, \cite{lebedev67},  or \cite{mclinden74}. For an introduction into the use of primal-dual methods in convex optimization, we refer the reader to \cite{pesquetplaying}. More recently, without being exhaustive, there were the notable works \cite{combettes2012}, \cite{chambolle2011first}, \cite{condat2012primal}, \cite{vu2011splitting}, and \cite{chambolle2016ergodic} which examined problems quite similar to the one posed here using first order primal-dual methods. 

In particular, \cite{chambolle2016ergodic} studied \eqref{probl} using $D$-$\prox$ mappings, i.e., proximal mappings where the euclidean energy has been replaced by a suitable Bregman divergence, under the assumption that $f$ and $h^*$ are Lipschitz-smooth $\Gamma_0$ functions and that the entropies $\phi_p$ and $\phi_d$ are strongly convex. They show ergodic convergence of the Lagrangian optimality gap with a rate of $O\para{1/k}$ under mild assumptions and also faster rates, e.g., $O\para{1/k^2}$ and linear convergence, under stricter assumptions involving strong convexity. We generalize their results by relaxing the Lipschitz-smooth assumption to a relative smoothness assumption, by analyzing the totally convex and relatively strongly convex case, by introducing stochastic error to the algorithm, and by showing almost sure weak convergence of the pointwise iterates themselves. Additionally, the recent work \cite{jiang2021bregman} studied a variant of the problem considered in \cite{chambolle2016ergodic} focused on semidefinite programming with $D$-$\prox$ mappings and an adaptive step size. As in \cite{chambolle2016ergodic}, they assume that the entropies inducing the Bregman divergences are strongly convex, in contrast to our work. The authors in \cite{Nguyen15} proposed a Bregman primal-dual method that iteratively constructs the best Bregman approximation to an arbitrary point from the Kuhn-Tucker set of a composite monotone inclusion in real reflexive Banach spaces, and for which they established strong convergence of the iterates. When specialized to structured minimization, their framework covers \eqref{probl} but without the smooth parts nor infimal-convolutions or the constraint sets $\C_p$ and $\C_d$. Moreover, their algorithm necessitates a complicated Bergman projection step and they do not consider stochastic versions.

Generalizations of \cite{chambolle2016ergodic} involving inexactness already exist in the form of \cite{rasch2020inexact} and \cite{chambolle2018stochastic}, however, \cite{rasch2020inexact} only considers determinstic inexactness and proximal operators computed in the euclidean sense, i.e., with entropy equal to the euclidean energy, and requires Lipschitz-smoothness. It's worth noting that the inexactness considered in their paper allows for the inexact computation of the proximal operators, in contrast to our work. While \algref{alg:ibpds} allows for inexactness, in the form of stochastic error, it is only allowable in the computation of gradient terms. The paper \cite{chambolle2018stochastic} allows for a very particular kind of stochastic error in which one randomly samples a set of indices at each iteration in an arbitrary but fixed way, i.e., according to some fixed distribution. However, the stochastic error we consider in the present paper is more general while encompassing the previous cases, although with less sharp results if the noise is not well behaved.

Another related work is that of \cite{hamedani2018primal} which generalizes the problem considered in \cite{chambolle2016ergodic} by allowing for a nonlinear coupling $\Phi\para{x,\mu}$ in \eqref{probl} instead of $\ip{Tx,\mu}{}$, although they maintain essentially the same Lipschitz-smoothness assumptions as in \cite{chambolle2016ergodic} translated to $\Phi\para{x,\mu}$. They are able to show a $O\para{1/k}$ convergence rate for the ergodic Lagrangian optimality gap under mild assumptions and an accelerated rate $O\para{1/k^2}$ when $g$ in \eqref{probl} is strongly convex with another assumption on the coupling $\Phi\para{x,\mu}$.

The notion of relative smoothness is key to the analysis of differentiable but not Lipschitz-smooth optimization problems. The earliest reference to this notion can be found in an economics paper \cite{birnbaum2011distributed} where it is used to address a problem in game theory involving fisher markets. Later it was parallelly developed for Bregman Forward-Backward splitting in \cite{BauschkeBolteTeboulle} and then in \cite{lu2018relatively} (see also \cite{Nguyen17,Bui2021}), and coined relative smoothness in \cite{lu2018relatively}. This idea allows one to apply arguments involving descent lemmas which are normally relegated to Lipschitz-smooth problems and it has been extended, for instance to define relative Lipschitz-continuity in \cite{lu2019relative}, in \cite{lu2020generalized} for the stochastic generalized conditional gradient, and to define a generalized curvature constant for the generalized conditional gradient algorithm in \cite{silvetisiopt}. The analogous idea of relative strong convexity, while noted before in \cite{chambolle2016ergodic}, was not explored in detail; here we analyze our algorithm under such assumptions in combination with total convexity of the entropies.

To our knowledge, our paper is the first to analyze \eqref{probl} under a relative smoothness condition with $D$-$\prox$ mappings. Additionally, we are the first to include stochastic error in the computation of the gradient terms for \eqref{probl} under these assumptions.
\end{subsection}

\begin{subsection}{Paper Organization}

The rest of the paper is divided into four sections. In \secref{sec:prelim}, we recall some basic definitions to make precise all the notions used in the paper along with some useful elementary results regarding sequences of random variables. 

In \secref{sec:est}, we make explicit all the assumptions \ref{ass:legendre}-\ref{ass:strconvex} we will use on the objective functions, entropies, step sizes, etc. We go on to establish the main estimation of \lemref{estimate} under \ref{ass:hyp} and \ref{ass:legendre}-\ref{ass:littled} that will be used in the convergence analysis of the ergodic, pointwise, and relatively strongly convex cases. The key idea is to utilize the descent lemma given by relative smoothness along with the usual inequalities for $\Gamma_0$ functions to estimate the optimality gap $\LL{x_{k},\mu} - \LL{x,\mu_{k}}$ in terms of the Bregman divergences induced by the entropies $\phi_p$ and $\phi_d$. The proof of the estimation here is similar in spirit to the proof of the main estimation in \cite{chambolle2011first}, with the main difference being that we are unable to use Young's inequality to deal with the coupling terms, which we handle using \ref{ass:littled}. There are also some lemmas involving \ref{ass:hyp} and \ref{ass:legendre}-\ref{ass:error} regarding the stochastic error, culminating in a summability result for the sequences $\seq{\TEX{\ip{\wnk,w-w_{k+1}}{}}}$ and $\seq{\EX{\ip{\wnk,w-w_{k+1}}{}}{\filts_k}}$ which appear in the convergence analysis.

In \secref{sec:conv}, we use the estimation developed in \secref{sec:est} along with \ref{ass:hyp} and \ref{ass:legendre}-\ref{ass:strconvex} regarding the entropies $\phi_p$ and $\phi_d$ and the regularity of their induced Bregman divergences to show convergence of the algorithm; first convergence of the expectation of the Lagrangian optimality gap for the ergodic iterates under \ref{ass:hyp} and \ref{ass:legendre}-\ref{ass:error} and then almost sure weak convergence of the pointwise iterates under \ref{ass:hyp} and \ref{ass:legendre}-\ref{ass:Tconsistent}. Finally, we examine the case where \ref{ass:strconvex} holds, i.e., there is relative strong convexity of the objective functions with respect to the entropies, and total convexity of the entropies themselves. For the ergodic analysis, denote by $\para{x_\infty,\mu_\infty}$ an almost sure weak sequential cluster point of the ergodic primal-dual sequence $\seq{\para{\bar{x}_k,\bar{\mu}_k}}$. Then we show that its expectation, namely $\TEX{\para{x_\infty,\mu_\infty}}$, is a saddle point. We prove also, for every $x$ and $\mu$, convergence of the expectation of the Lagrangian optimality gap $\TEX{\LL{\bar{x}_{k},\mu}-\LL{x,\bar{\mu}_{k}}}$ with a rate of $O\para{1/k}$. For the pointwise analysis, we begin by showing an almost sure asymptotic regularity result for the primal-dual sequence $\seq{w_k}$. With this, we are then able to adapt the well known Opial's lemma (see \cite{opial1967weak}) to the Bregman primal-dual setting to establish almost sure weak convergence of the primal-dual sequence $\seq{w_k}$ to a saddle point $\ws$. In the final part of this section, we establish almost sure strong convergence of the primal-dual sequence $\seq{w_k}$ to a saddle point $\ws$ under \ref{ass:strconvex} and total convexity of the entropies.

Lastly, in \secref{sec:app}, we explore potential applications of the algorithm and demonstrate numerically its effectiveness when applied to two different problems. The first is a simple linear inverse problem on the simplex with total variation regularization, which we examine in the deterministic and stochastic case. The second is an application in optimal transport involving the entropically regularized Wasserstein distance and inverse problems. There is also a discussion of other possible applications of the algorithm to entropic Wasserstein barycenter problems.
\end{subsection}
%-----------------------
% Footer to allow Gummi
% to maintain previews
% for individual pieces
%-----------------------
\ifdefined\COMPLETE
\else
\end{document}
\fi

%----------------------------------
% These commands allow Gummi to
% edit each piece individually
% and still maintain previews.
%----------------------------------
\ifdefined\COMPLETE
\else
\documentclass[12pt]{article}
\input{tex_package_header} %file containing all the used libraries
\begin{document}
\fi
%----------------------------------
% Remove at your own risk!
% There is also a footer that ends
% the document if the main wasn't
% loaded.
%----------------------------------
\section{Notations and preliminary facts}\label{sec:prelim}
\subsection{Basic notation}
Given a real reflexive Banach space $\X$, we denote by $\X^*$ its topological dual and by $\ip{u,x}{}$ the duality pairing for $x\in\X$ and $u\in\X^*$. The norm on $\X$ is denoted $\norm{\cdot}{\X}$. The symbols $\rightharpoonup$ and $\to$ denote respectively weak and strong convergence.
%We say that a sequence $\seq{x_k}$ with $x_k\in\X$ for each $k\in\N$ converges strongly to some $x\in\X$ iff $\norm{x_k-x}{\X}\to 0$. On the other hand, we say a sequence $\seq{x_k}$ with $x_k\in\X$ for each $k\in\N$ converges weakly to some $x\in\X$, denoted $x_n\rightharpoonup x$, iff $\ip{u,x_k}{}\to\ip{u,x}{}$ for every $u\in\X^*$. 
The set of weak sequential cluster points of a sequence $\seq{x_k}$ in $\X$ is defined as 
\nnewq{\label{eq:wkset}
\wkset{x}{k} \eqdef \brac{x\in\X: \exists \subseq{x_{k_j}}, x_{k_j}\rightharpoonup x}.
}

For a function $f \in \Gamma_0\left(\X\right)$, $\partial f: \X \to 2^\X$ is its subdifferential operator. When referring to the differentiability or the gradient of a function $f:\X\to\R$, it is meant in the sense of G\^ateaux. For a non-empty closed convex set $\C \subset \X$, $N_{\C}(x)$ is the normal cone of $\C$ at $x \in \C$. $\inte \C$ and $\conj{\C}$ denote the interior and the closure of a set $\C \subset\X$.

%Finally, when referring to the interior of a set $U$, denoted $\inte U$, or the closure, denoted $\conj{U}$, we mean with respect to the norm topology on $\X$.

\subsection{Bregman divergence notation}
We denote by $D$, without subscript, the Bregman divergence associated to $\phi(x,\mu)\eqdef \phi_p(x) + \phi_d(\mu)$; namely, given $w_i\eqdef\left(x_i,\mu_i\right)$ with $\left(x_i,\mu_i\right)\in\X_p \times \X_d$ for $i\in\brac{1,2}$,
\begin{equation*}
	\begin{split}
		D\left(w_1,w_2\right)& \eqdef D_p\left(x_1,x_2\right)+D_d\left(\mu_1,\mu_2\right).
	\end{split}
\end{equation*}
We proceed with some notions about regularity of functions.
\begin{definition}[Legendre function]
The function $\phi$ is called a Legendre function if $\partial \phi$ is both locally bounded and single-valued on its domain, $\para{\partial \phi}^{-1}$ is locally bounded on its domain, and $\phi$ is strictly convex on every convex subset of $\dom\para{\partial\phi}$.
\end{definition}
\begin{definition}[Relative smoothness]\label{gensmooth}
Given a function $\phi: \X \to \R\cup\brac{\pinfty}$ differentiable on $\inte\dom\para{\phi}$, we say that the function $f: \X \to \R\cup\brac{\pinfty}$ is $L$-smooth with respect to $\phi$ if it is differentiable on $\inte\dom\para{\phi}$ and $L\phi - f$ is convex on $\inte\dom\para{\phi}$; namely, if for every $x,y\in \inte\dom\para{\phi}$
	$$D_f(x,y)\leq L D_{\phi}(x,y).$$
\end{definition}
\begin{remark}
The relative smoothness property, used notably in \cite{BauschkeBolteTeboulle}, \cite{Nguyen17} and \cite{lu2018relatively}, implies the following fact which can be interpreted as a "generalized descent lemma": for every $x,y\in \inte\dom\para{\phi}$,
\begin{align}\label{gendl}
f(x)  \leq f(y) + \langle \nabla f (y),x-y\rangle + L D_{\phi}\left(x,y\right).
\end{align}
When $\phi$ is the euclidean square norm, or energy, relative smoothness is equivalent to Lipschitz-smoothness, i.e., Lipschitz-continuity of the gradient of $f$.
\end{remark}
\begin{definition}[Relative strong convexity]
\label{genstrong}
Given a function $\phi: \X \to \R\cup\brac{\pinfty}$ differentiable on $\inte\dom\para{\phi}$, and a non-empty closed convex set $\C \subset \dom(f) \cap \dom(\phi)$, we say that $f: \X \to \R\cup\brac{\pinfty}$ is $m$-strongly convex on $\C$ with respect to $\phi$ %if $f-m\phi$ is convex on $\inte\dom\para{\phi}$
if for every $x \in \C$ and $y \in \inte\dom\para{\phi} \cap \dom(\partial f) \cap \C$
\[
f(x) - f(y) - \iprod{u}{y-x} \geq m D_{\phi}(x,y) , \qforallq u \in \partial f(y) .
\]
\end{definition}
Note that the idea of relative strong convexity can be found in a footnote of \cite{chambolle2016ergodic} but it was not explored further. With these definitions, we use the following notation to improve readability
\begin{equation}\label{notation}
\begin{split}
\left(\frac{1}{\Lambda_k}-L\right) D\left(w_1,w_2\right)& \eqdef \left(\frac{1}{\lambda_k}-L_p\right) D_p\left(x_1,x_2\right)+\left(\frac{1}{\nu_k}-L_d\right)D_d\left(\mu_1,\mu_2\right)\\
\para{\frac{1}{\Lambda_\infty} - L}D\para{w_1,w_2} & \eqdef \para{\frac{1}{\lambda_\infty}-L_p}D_p\para{x_1,x_2}+\para{\frac{1}{\nu_{\infty}}-L_d}D_d\para{\mu_1,\mu_2}\\
M(w_1,w_2) & \eqdef \langle T (x_1-x_2),\mu_1-\mu_2 \rangle
\end{split}
\end{equation}
where $\lambda_k$ and $\nu_k$ are the step-sizes in Algorithm~\ref{alg:ibpds}, and $L_p$ and $L_d$ are the constants introduced in \ref{ass:legendre}.

\subsection{Probabilistic notation and preliminaries}
We denote by $\prspace$ a probability space with set of events $\events$, $\sigma$-algebra $\sigalg$, and probability measure $\prob$. Throughout, we assume that any real reflexive Banach space $\mc{X}$ is endowed with its Borel $\sigma$-algebra, $\borel{\mc{X}}$. Formally, we define the stochastic primal and dual errors at iteration $k$ as $\pnk$ and $\dnk$, i.e., $\pnk$ and $\dnk$ are measurable functions from $\events$ to $\X_p^*$ and $\X_d^*$ with their respective Borel $\sigma$-algebras. When it makes sense, we will also denote the combined error as $\wnk$ in the same way that we use $w_k$, e.g.,
\newq{
\ip{\wnk, w-w_k}{} \eqdef \langle\pnk, x-x_k\rangle + \langle\dnk, \mu-\mu_k\rangle.
}
We denote a filtration on $\prspace$ by $\Filt \eqdef \seq{\filt_k}$ where $\filt_k$ is a sub-$\sigma$-algebra satisfying\fekn $\filt_k\subset \filt_{k+1}\subset\sigalg$. Furthermore, given a set of random variables $\brac{a_0,\ldots,a_n}$ we denote by $\sigma\para{a_0,\ldots,a_n}$ the $\sigma$-algebra generated by $a_0,\ldots,a_n$. Finally, an expression $\para{P}$ is said to hold $\Pas$ if 
\newq{
\prob\para{\brac{\omega\in\events: \para{P}\mbox{ holds}}}=1.
}
Using the above notation, we denote the canonical filtration associated to the iterates of the algorithm as $\Filts \eqdef \seq{\filts_k}$ with, for all $k\in\N$,
\newq{
\filts_k \eqdef \sigma\brac{\para{x_0,\mu_0},\para{x_1,\mu_1},\ldots,\para{x_k,\mu_k}}
}
such that all iterates up to $\para{x_{k},\mu_{k}}$ are completely determined by $\filts_k$. 

For the remainder of the paper, all equalities and inequalities involving random quantities should be understood as holding $\Pas$ even if it is not explicitly written. 

\begin{definition}
Given a filtration $\Filt$, we denote by $\ell_+\para{\Filt}$ the set of sequences of $[0,+\infty[$-valued random variables $\seq{a_k}$ such that, for each $k\in\N$, $a_k$ is $\filt_k$ measurable. Then, we also define the following set of summable random variables,
\newq{
\ell^1_+\para{\Filt} \eqdef \brac{\seq{a_k}\in\ell_+\para{\Filt}:\sum\limits_{k\in\N} a_k<+\infty \Pas}.
}
The set of non-negative summable sequences is denoted $\ell^1_+$.
\end{definition}

%\paragraph{Preliminary background}

The following probabilistic results will be useful in the convergence analysis of \algref{alg:ibpds}.

\begin{lemma}[{Robbins-Siegmund, \cite[Theorem 1]{robbinssiegmund}}]\label{AuxPC}
Given a filtration $\Filt$ and the sequences of real-valued random variables $\seq{r_k}\in\ell_+\para{\Filt}$, $\seq{a_k}\in\ell_+\para{\Filt}$, and $\seq{z_k}\in\ell_+^1\para{\Filt}$ satisfying\fekn
\newq{
\EX{r_{k+1}}{\filt_k} - r_k \leq -a_k + z_k \Pas
}
it holds that $\seq{a_k}\in\ell^1_+\para{\Filt}$ and $\seq{r_k}$ converges $\Pas$ to a random variable with value in $[0,+\infty[$.
\end{lemma}
% \begin{proof}
% See \cite[Theorem 1]{robbinssiegmund}.
% \end{proof}

\begin{remark}\label{LemmaAux}
	In the deterministic case, \lemref{AuxPC} reduces to the following statement. Let $\seq{a_k}\in\ell_+, \seq{r_k}\in\ell_+$ and $\seq{z_k}\in\ell^1_+$ such that\fekn
	\newq{
	r_{k+1}-r_k\leq -a_k + z_k.
	}
	Then $\seq{a_k}\in\ell^1_+$ and $\seq{r_k}$ converges to $r\in[0,+\infty[$. This result is \cite[Lemma 2, page 44]{polyak1987introduction}.
\end{remark}

\begin{lemma}\label{AuxAS}
If $\seq{x_n}$ is a sequence of $\X$-valued random variables such that $\seq{\mathbb{E}\left(\norm{x_k}{\X}^q\right)}\in \ell^1_+$ for some $q\in]0,+\infty[$, then $x_k \to 0$ almost surely.
\end{lemma}

\begin{proof}
	For every $\varepsilon>0$, by Markov's inequality,
	\begin{equation}
	\begin{split}
	\sum_{n=0}^{N} \mathbb{P} \left(\norm{x_n}{\X}^q\geq \varepsilon\right) & \leq \frac{1}{\varepsilon}\sum_{n=0}^{N} \mathbb{E} \left(\norm{x_n}{\X}^q\right).
	\end{split}
	\end{equation}
	Taking the limit for $N\to+\infty$ and using the assumption $\seq{\mathbb{E}\left(\norm{x_k}{\X}^q\right)}\in \ell^1_+$, we get that, for every $\varepsilon>0$, it holds $\mathbb{P} \left(\norm{x_n}{\X}^q\geq \varepsilon\right)$ belongs to $\ell^1_+$. As a consequence of the Borel-Cantelli Lemma, $\norm{x_n}{\X}^q\to 0$ almost surely whence the claim follows.
\end{proof}

%-----------------------
% Footer to allow Gummi
% to maintain previews
% for individual pieces
%-----------------------
\ifdefined\COMPLETE
\else
\end{document}
\fi

%----------------------------------
% These commands allow Gummi to
% edit each piece individually
% and still maintain previews.
%----------------------------------
\ifdefined\COMPLETE
\else
\documentclass[12pt]{article}
\input{tex_package_header} %file containing all the used libraries
\begin{document}
\fi
%----------------------------------
% Remove at your own risk!
% There is also a footer that ends
% the document if the main wasn't
% loaded.
%----------------------------------
\section{Main assumptions and estimations}\label{sec:est}
\begin{subsection}{Main assumptions}\label{sec:assump}
We first state our assumptions and then remark on their motivations and common examples where they hold. Note that for several results, only a subset of these assumptions are needed; we will comment ont this hereafter. For brevity, throughout the remainder of the paper we employ the following notation
\begin{alignat*}{2}
\U_p &\eqdef  \inte\dom\para{\phi_p} \cap \dom \para{\partial g}\quad\quad&&\tilde{\U}_p\eqdef \dom\para{\phi_p}\cap\dom\para{\p g}\\ \U_d &\eqdef  \inte\dom \para{\phi_d} \cap \dom \para{\partial l^{*}}\quad &&\tilde{\U}_d\eqdef \dom\para{\phi_d}\cap\dom\para{\p l^*}.
\end{alignat*}
\begin{enumerate}[label=\bf{(\subscript{\mathbf{A}}{{\arabic*}})}]\label{ass:list}
	\item\label{ass:legendre}\label{relsmoothcond} The entropies $\phi_p$ and $\phi_d$ belong to  $\Gamma_0\left(\X_p\right)$ and $\Gamma_0\left(\X_d\right)$ with $\dom (\phi_p) \times \dom(\phi_d)=\C_p \times \C_d$ and with $f$ and $h^{*}$ being $L_p$ and $L_d$ - smooth wrt $\phi_p$ and $\phi_d$, respectively (see Definition \ref{gensmooth}). The $D$-$\prox$ mappings $\para{\nabla \phi_p+\lambda_k\partial g}^{-1}$ and $\para{\nabla\phi_d+\nu_k\partial l^*}^{-1}$ are well-defined (i.e., nonempty and single-valued) maps from $\inte\dom\para{\phi_p}$ and $\inte\dom\para{\phi_d}$ to $\inte\dom\para{\phi_p}$ and $\inte\dom\para{\phi_d}$, respectively.
	\item\label{ass:steps} The step size sequences $\seq{\lambda_k}$ and $\seq{\nu_k}$ are positive, nondecreasing, and bounded above with their limits denoted $\lim\limits_{k\to\infty}\lambda_k=\lambda_\infty$ and $\lim\limits_{k\to\infty}\nu_k=\nu_{\infty}$.
	\item\label{ass:littled} The step sizes satisfy \ref{ass:steps} and one of the following holds: \begin{enumerate}[label=\bf{(\Roman*)}]
				\item\label{littledweak} there is a function $d: \left(\X_p \times \X_d\right)^2 \to \R_+$ and $\eps \geq 0$ such that
	\begin{align}
		\inf_{\sst{w_1\in\tilde{\U}_p\times\tilde{\U}_d, \ w_2\in \U_p \times \U_d;\\ w_1 \neq w_2}} \frac{\left(\frac{1}{\Lambda_\infty}-L\right) D(w_1,w_2) - M(w_1,w_2)}{d(w_1,w_2)} \geq \eps;
		\label{hyp2}
	\end{align}
				\item\label{littledstrong} the above holds with $\eps >0$.
			\end{enumerate}
	\item\label{ass:unbiased} The error sequence $\seq{\Delta_k}$ is unbiased conditioned on the filtration $\Filts$, i.e.\fekn
\newq{
\EX{\pnk}{\filts_{k}} = \EX{\dnk}{\filts_{k}} = 0.
}
	\item\label{ass:error} One of the following holds:
		\begin{enumerate}[label=\bf{(\Roman*)}]
			\item\label{determinerror} for each $k\in\N$, the stochastic errors $\pnk$ and $\dnk$ are zero almost surely;
			\item\label{boundederror} the following sequences satisfy
			\newq{\seq{\EX{\norm{\pnk}{\X_p^*}}{\filts_k}}\in\ell^1_+\para{\Filts}&\qandq\seq{\TEX{\norm{\pnk}{\X_p^*}}}\in\ell^1_+\\
			\seq{\EX{\norm{\dnk}{\X_d^*}}{\filts_k}}\in\ell^1_+\para{\Filts}&\qandq\seq{\TEX{\norm{\dnk}{\X_d^*}}}\in\ell^1_+} and the sets $\U_p$ and $\U_d$ are bounded, i.e., $0<\diam_{\U_p}<+\infty$ and the same for $\diam_{\U_d}$;
			\item\label{strconvexerror} the entropies $\phi_p$ and $\phi_d$ are strongly convex with respect to $\norm{\cdot}{\X_p}^2$ and $\norm{\cdot}{\X_d}^2$ with moduli $\scphip$ and $\scphid$, respectively. Additionally, the step sizes $\seq{\lambda_k}$ and $\seq{\nu_k}$ satisfy \ref{ass:steps} with
			\newq{
			\nu_\infty\lambda_\infty < \frac{\scphip\scphid}{\norm{T}{p\to d^*}^2},
			}
			where $\norm{\cdot}{p\to d^*}$ is the standard operator norm between $\X_p$ and $\X_d^*$ and the following sequences satisfy
			\newq{
			\EX{\norm{\pnk}{\X_p^*}^2}{\filts_k}\in\ell^1_+\para{\Filts}\qandq\TEX{\norm{\pnk}{\X_p^*}^2}\in\ell^1_+.
			}
		\end{enumerate}
	\item\label{ass:dconsistent} For the function $d$ used in \eqref{hyp2} and all bounded sequences $\seq{v_k}$ and $\seq{z_k}$ in $\inte\dom\para{\phi} \eqdef \inte\dom\para{\phi_p}\times\inte\dom\para{\phi_d}$ 
	\begin{equation}\label{assumption_psi}
	d(v_k,z_k) \to 0 \  \ \ \  \Rightarrow \ \ \ \ v_k-z_k \to 0.
	\end{equation}
	\item\label{ass:coercive} For every $w\eqdef\para{x,\mu} \in \inte\dom\para{\phi}$, at least one of $D(w,\cdot)$ or $d(w,\cdot)$ is coercive.
%	\item\label{ass:maxmono} The set-valued operators $\partial g + \partial f + N_{\C_p}$ and $\partial l^* + \partial h^* + N_{\C_d}$ are maximal monotone.
	\item\label{ass:gradreg} For any bounded sequence $\seq{w_k}$ with $w_k\in\inte \dom \phi$ for each $k\in\N$, if $w_{k+1}-w_k \to 0$, then
	\begin{align*}
	\nabla\phi_p(x_{k+1})-\nabla\phi_p(x_k) \to 0 &\qandq \nabla f(x_{k+1})-\nabla f(x_{k}) \to 0; \\
	\nabla\phi_d(\mu_{k+1})-\nabla\phi_d(\mu_{k}) \to 0 &\qandq \nabla h^*(\mu_{k+1})-\nabla h^*(\mu_{k}) \to 0 .
	\end{align*}
	\item\label{ass:gradphiconsistent} For any sequence $\seq{w_k}$ with $w_k\in\inte \dom \phi$\fekn if $w_{k} \rightharpoonup w_{\infty}$, then
	\begin{align*}
	\nabla\phi(w_{k})\rightharpoonup \nabla\phi(w_{\infty}).
	\end{align*}
	\item\label{ass:Tconsistent} For an arbitrary sequence $\seq{w_k}\in \X_p \times\X_d$, if $w_{k} \rightharpoonup 0$, then
	\begin{align*}
	\langle T x_k,\mu_k \rangle \to 0.
	\end{align*}
	\item\label{ass:strconvex} 
		\begin{enumerate}[label=\bf{(\Roman*)}]
		\item At least one of the functions $f$ or $g$ is relatively strongly convex on $\C_p \cap \dom(g)$ wrt an entropy $\psi_p:\X_p\to\R\cup\brac{\pinfty}$ with constant $\scmf$ or $\scmg$, respectively (see Definition \ref{genstrong}). The entropy $\psi_p$ satisfies $\dom\para{\phi_p}\subseteq \dom\para{\psi_p}$. \label{ass:strconvexP}
		\item At least one of the functions $h^*$ or $l^*$ is relatively strongly convex on $\C_d \cap \dom(l^*)$ wrt an entropy $\psi_d:\X_d\to\R\cup\brac{\pinfty}$ with constant $\scmh$ or $\scml$, respectively (see Definition \ref{genstrong}). The entropy $\psi_d$ satisfies $\dom\para{\phi_d}\subseteq \dom\para{\psi_d}$. \label{ass:strconvexD}
		\end{enumerate}
\end{enumerate}
Analogously to \eqref{notation}, we will also use the shorhand notation using the relative strong convexity constants and entropies from \ref{ass:strconvex}:
\nnewq{\label{notationstr}
	\scs D'\left(w_1,w_2\right)  & \eqdef \scmf D_{\psi_p}\left(x_1,x_2\right)+\scmh D_{\psi_d}\left(\mu_1,\mu_2\right)\\
	\scns D'\left(w_1,w_2\right) & \eqdef \scmg D_{\psi_p}\left(x_1,x_2\right)+\scml D_{\psi_d}\left(\mu_1,\mu_2\right).
}
\begin{remark}[\ref{ass:legendre} and \ref{ass:steps}]
There are several, technical characterizations of sufficient conditions that ensure the latter half of \ref{ass:legendre} holds. Classical examples start by assuming that the spaces are reflexive and that $\phi_p$ and $\phi_d$ are Legendre functions and then add assumptions depending on the space being considered; see for instance, the comprehesinve treatment in \cite[Section~3]{BauschkeBorweinCombettes}. Notice that we do not require $\phi_p$ and $\phi_d$ to be Legendre in general, that is indeed incompatible with \ref{ass:gradphiconsistent} if the limit point is on the boundary. In practice, the latter half of \ref{ass:legendre} is required only for the existence and uniqueness of the sequence generated by the algorithm and is not used explicitly elsewhere in the convergence analysis. For \ref{ass:steps}, it is sufficient to take the step sizes $\seq{\lambda_k}$ and $\seq{\nu_k}$ to simply be constant.
\end{remark}
\begin{remark}[\ref{ass:littled}]
The infimum in \ref{ass:littled} is taken with $w_1\in\tilde{\U}_p\times\tilde{\U}_d$ and $w_2\in\U_p\times\U_d$ because, a priori, a solution $w^\star$ may lie in the boundary of $\tilde{\U}_p\times\tilde{\U}_d$ even if the iterates $\seq{w_k}$ themselves remain in $\U_p\times\U_d$ due to \ref{ass:legendre}. Since the Bregman divergence is still well defined when the first argument (but not the second) is in $\dom\para{\phi}\setminus\inte\dom\para{\phi}$, there is no issue with taking the infimum over this set. Observe that \ref{ass:littled} also entails that, for every $w_1\in\tilde{\U}_p\times\tilde{\U}_d$ and $w_2\in \U_p \times \U_d$\fekn
\begin{equation}\label{aux}
	\frac{1}{\Lambda_k} D\left(w_1,w_2\right)-M\left(w_1,w_2\right) \geq L D\left(w_1,w_2\right) + \eps  d\left(w_1,w_2\right) \geq 0 .
\end{equation}
\end{remark}
\begin{example}
	Suppose that $\varphi: \R_+ \to \R_+$ is a convex nondecreasing function with $\varphi^*$ its positive conjugate and $\gau$ a finite coercive gauge with domain $\R_+\para{\U_p -\U_p}\subset\X_d$ (in the Minkowski sense) and polar $\gau^{\circ}$. Assume that the quantities defined by
	\begin{equation*}
	\begin{split}
	\norm{T}{D_p}\eqdef \sup_{x_1,x_2\in \U_p; \ x_1\neq x_2} \ \frac{\varphi(\gau\left(T(x_2-x_1)\right))}{D_p\left(x_1,x_2\right)}\qquad\mbox{and}\qquad \norm{I}{D_d}\eqdef \sup_{\mu_1,\mu_2 \in \U_d; \ \mu_1\neq \mu_2} \ \frac{\varphi^*(\gau^{\circ}\left(\mu_2-\mu_1\right))}{D_d\left(\mu_1,\mu_2\right)}
	\end{split}
	\end{equation*}	
	are finite. We use the notation $\norm{\cdot}{D_p}$ and $\norm{\cdot}{D_p}$, but notice that they may not be norms.
	%, recalling that, for all $x_1, x_2\in \inte\ \dom\para{\phi_p}$,
	% \begin{equation}\label{bregprop}
	% \langle \nabla\phi_p(x_2)-\nabla\phi_p(x_1), x_2-x_1\rangle = D_p\left(x_1,x_2\right)+D_p\left(x_2,x_1\right).
	% \end{equation}	
	% Define also
	% \begin{equation*}
%	\alpha_p \eqdef  \inf_{x_1\neq x_2} \brac{ \frac{D_p\left(x_1,x_2\right)}{D_p\left(x_2,x_1\right)}: x_1,x_2\in\inte\ \dom\para{\phi_p}}\in[0,1]
%	\end{equation*}	
%	and the analogue for $\alpha_d$, the so called symmetry coefficients given in \cite[Definition~2]{BauschkeBolteTeboulle}. 
If, moreover, we suppose that the step sizes verify\fekn for some $\eps_k\geq 0$,
	\begin{equation}\label{eq:assparam}
	\left(\frac{1}{\lambda_k}-L_p\right) \geq \norm{T}{D_p} + \eps_k \qandq    
	\left(\frac{1}{\nu_k}-L_d\right) \geq \norm{I}{D_d} + \eps_k,
	\end{equation}	
	then \ref{ass:littled} is satisfied with $d\para{w_1,w_2}=D\para{w_1,w_2}$. Indeed, for any pair $w_1, w_2 \in \U_p \times \U_d$, we have\fekn
	\begin{equation}\label{eq:ineqpos}
	\begin{split}
	&\left(\frac{1}{\Lambda_k}-L\right) D\left(w_1,w_2\right)-M\left(w_1,w_2\right) \\
	\ \\
	= \ \ \ & \left(\frac{1}{\lambda_k}-L_p\right) D_p\left(x_1,x_2\right)+\left(\frac{1}{\nu_k}-L_d\right)D_d\left(\mu_1,\mu_2\right)-\langle T (x_1-x_2),\mu_1-\mu_2 \rangle\\
	 \ \\
	\geq \ \ \ & \norm{T}{D_p} D_p\left(x_1,x_2\right)+\norm{I}{D_d} D_d\left(\mu_1,\mu_2\right)-\gau\left(T (x_1-x_2)\right)\gau^{\circ}\left(\mu_1-\mu_2\right) + \eps_k D\left(w_1,w_2\right) \\
	\ \\
	\geq \ \ \ &  \varphi(\gau\left(T(x_1-x_2)\right)) +\varphi^{*}(\gau^{\circ}\left(\mu_1-\mu_2\right)) -\varphi(\gau\left(T (x_1-x_2)\right)) - \varphi^{*}(\gau^{\circ}\left(\mu_1-\mu_2\right)) + \eps_k D\left(w_1,w_2\right)\\
	\ \\
	% \geq \ \ \ & \left(1+\alpha_p^{-1}\right) \ \frac{\varphi(\gau\left(T(x_1-x_2)\right))}{(1+\alpha_p^{-1})D_p\left(x_1,x_2\right)} \ D_p\left(x_1,x_2\right)+(1+\alpha_d^{-1}) \ \frac{\varphi^{*}(\gau^{\circ}\left(\mu_1-\mu_2\right))}{(1+\alpha_d^{-1})D_d\left(\mu_1,\mu_2\right)} \ D_d\left(\mu_1,\mu_2\right)\\
	% &-\varphi(\gau\left(T (x_1-x_2)\right)) - \varphi^{*}(\gau^{\circ}\left(\mu_1-\mu_2\right)) + \eps_k D\left(w_1,w_2\right)\\
	% = \ \ \ & \varphi(\gau\left(T (x_1-x_2)\right)) + \varphi^{*}(\gau^{\circ}\left(\mu_1-\mu_2\right))-\varphi(\gau\left(T (x_1-x_2)\right)) - \varphi^{*}(\gau^{\circ}\left(\mu_1-\mu_2\right)) +  \eps_k D\left(w_1,w_2\right)\\
    = \ \ \ &  \eps_k D\left(w_1,w_2\right).
	\end{split}
	\end{equation}

Note that in this example we have taken the action of $T$ on the primal variables into the definition of $\norm{\cdot}{D_p}$. It is equally possible, and sometimes desirable, to define things such that the action of the adjoint $T^*$ on the dual variables is incorporated into $\norm{\cdot}{D_d}$ instead, which can change the values (and consequently step sizes) in a non-Hilbertian setting.
	\end{example}
\begin{remark}[\ref{ass:unbiased} and \ref{ass:error}]\label{remark_errorconv}
Notice that, using Lemma \ref{AuxAS}, \ref{ass:unbiased} and \ref{ass:error} (in any case) imply that $\seq{\pnk}$ and $\seq{\dnk}$ converge strongly (with respect to $\norm{\cdot}{\X_p^*}$ and $\norm{\cdot}{\X_d^*}$ respectively) to zero a.s. and that, furthermore, \faw $\seq{\TEX{\ip{\wnk, w-w_{k+1}}{}}}\in\ell^1_+$ and $\seq{\EX{\ip{\wnk,w-w_{k+1}}{}}{\filts_k}}\in\ell^1_+\para{\Filts}$ (see \lemref{lem:errorsum} for details). In \ref{ass:error}\ref{strconvexerror}, the norms $\norm{\cdot}{\X_p}$ and $\norm{\cdot}{\X_d}$ can be replaced with arbitrary norms as long as $\phi_p$ and $\phi_d$ are strongly convex with respect to their square. The different cases for \ref{ass:error} can be mixed for the primal and dual, e.g., one can take \ref{ass:error}\ref{strconvexerror} for the primal but have \ref{ass:error}\ref{boundederror} for the dual; the current presentation simply for convenience.
\end{remark}
\begin{remark}[\ref{ass:dconsistent}]
In the case where $d\para{x,y}$ is the Bregman divergence induced by the Shannon-Boltzman entropy, the Hellinger entropy, the fractional-power entropy, the Fermi-Dirac entropy, or the energy/euclidean entropy, \ref{ass:dconsistent} holds (see \cite[Remark 4]{BauschkeBolteTeboulle}. 

More generally, when $d=D_{\zeta}$ for some entropy $\zeta$ which is Legendre, we have from \cite[Example~4.10]{BauschkeBorweinCombettes} that \ref{ass:dconsistent} is satisfied whenever one of the following holds
\begin{itemize}
\item $\zeta$ is uniformly convex on bounded sets;
\item $\X_p\times\X_d$ is finite dimensional, $\dom\para{\zeta}$ is closed, and $\zeta\mid_{\dom\para{\zeta}}$ is strictly convex and continuous.
\end{itemize}
Thus, if $\zeta = \phi$, with $\phi$ Legendre, we require only $\dom\para{\phi}$ to be closed if $\X_p\times\X_d$ is finite dimensional.
\end{remark}
\begin{remark}[\ref{ass:coercive}]
Sufficient conditions for \ref{ass:coercive} to hold for Legendre functions in real reflexive Banach spaces are given in \cite[Lemma 7.3(viii) \& (ix)]{BauschkeBorweinCombettes01}. For example, $D_p\para{x,\cdot}$ is coercive if $\phi_p$ is supercoercive and $x\in \inte\dom\para{\phi_p}$; alternatively, if $\mc{X}_p$ is finite-dimensional, $\dom \para{\phi_p^*}$ is open, and $x\in \inte\dom \para{\phi_p}$. Similar conditions hold for $\phi_d$.
\end{remark}
%\begin{remark}[\ref{ass:maxmono}]
%First, observe that by assumption \ref{ass:hyp} and \cite[Theorem~3.1.11]{Zalinescu02}, all individual operators $\partial g$, $\partial f$, $\partial l^*$, $\partial h^*$, $N_{\C_p}$ and $N_{\C_d}$ are maximal monotone. Thus, there are numerous sufficient qualification conditions for \ref{ass:maxmono} to hold in real reflexive Banach spaces; see \cite{borwein07} or \cite{Simons08}. Relaxed conditions exist also for finite dimension (see e.g., \cite{chu1996}).
%%For instance, we can impose one of the following (see \cite{borwein07})\footnote{Similar conditions apply to $\partial l^*$ and $\nabla h^*$.}:
%%\begin{itemize}
%%\item $\inte\dom\para{\partial g}\cap\inte\dom\para{\partial f}\neq \emptyset$ and its analog for the dual.
%%\item $\dom\para{\partial f}\cap\inte\dom\para{\partial g} \neq \emptyset$ while $\dom\para{\nabla f}$ is closed and convex.
%%\item $\dom\para{\partial f}$ and $\dom\para{\partial g}$ are closed and convex and $0\in\mbox{core}\para{\mbox{conv}\para{\dom\para{\nabla f} - \dom\para{\partial g}}}$, where for the definition of $\mbox{core}$ and $\mbox{conv}$ see \cite[Definition~6.9 and 3.3]{BAUSCHCOMB}. 
%%\end{itemize}
%\end{remark}
\begin{remark}[\ref{ass:Tconsistent}]
Assumption \ref{ass:Tconsistent} is very mild and holds when the operator $T$ (or $T^*$) is for instance compact.
\end{remark}

\begin{remark}[\ref{ass:gradreg}, \ref{ass:gradphiconsistent}] \label{rem:finitedim}
Assumptions \ref{ass:gradreg}, \ref{ass:gradphiconsistent} and \ref{ass:Tconsistent} are required only for the pointwise weak convergence of the iterates, namely in Section~\ref{sec:pointwise}. \ref{ass:gradreg} and \ref{ass:gradphiconsistent} have been previously assumed by other authors to prove weak convergence of the iterates for the Bregman Forward-Backward algorithm on a real reflexive Banach space; see \cite{Nguyen17,Bui2021}. In particular, \ref{ass:gradphiconsistent} is a weak sequential continuity assumption on the gradients of the entropies, while \ref{ass:gradreg} can be obtained for instance from norm-to-norm uniform continuity on bounded sets of $\nabla \phi_p$, $\nabla \phi_d$, $\nabla f$, and $\nabla h^*$. A typical example where these assumptions hold is when $\X_p$ is the $\ell_q$ space\footnote{We focus on the primal space $\X_p$ but the same reasoning applies to $\X_d$.}, $q \in ]1,+\infty[$, and $\phi_p=\norm{\cdot}{\ell_q}^q/q$, in which case $\nabla \phi_p$ is the duality mapping on $\ell_q$. The latter is known in this case to be weakly continuous \cite{Browder1966} and norm-to-norm uniformly continuous on every bounded subset of $\ell^q$ \cite{Cioranescu1990}. However, if the duality mapping is replaced with the normalized duality mapping, i.e., $\phi_p =\norm{\cdot}{\ell_q}^2/2$, then \ref{ass:gradphiconsistent} fails unless $q=2$ (i.e., Hilbertian setting) while \ref{ass:gradreg} still holds for $\phi_p$; see \cite{Xu14}. 

On the other hand, \ref{ass:gradreg} is satisfied when $\X_p\times\X_d$ is finite dimensional. Indeed, in finite dimension not only do strong and weak convergence coincide but also $\nabla \phi_p$, $\nabla \phi_d$, $\nabla f$, and $\nabla h^*$ are all continuous on the interior of their domains by \cite[Corollary~9.20]{rockafellar1998variational} since $\phi_p, f\in\Gamma_0\para{\X_p}$ and $\phi_d, h^*\in\Gamma_0\para{\X_d}$. Again, \ref{ass:gradphiconsistent} is more subtle even in finite dimension since Legenderness of the entropy entails that if an interior sequence converges to a point on the boundary of the domain of the entropy, the sequence of gradients will diverge.
\end{remark}

We finish this section by providing an infinite-dimensional example where all asumptions hereabove are verified.
\begin{example}
We give an example of an infinite-dimensional Banach space $\X_p$ and an entropy $\phi_p$ for which assumptions \ref{ass:legendre}, \ref{ass:gradreg} and \ref{ass:gradphiconsistent} both hold. Consider $\mathcal{H}$ an infinite-dimensional Hilbert space and $\mathcal{V}$ a finite-dimensional Banach space, with respective norms $\norm{\cdot}{\mathcal{H}}$ and $\norm{\cdot}{\mathcal{V}}$, and define $\mathcal{X}_p = \mathcal{H}\times\mathcal{V}$ to be the Banach space with norm $\norm{(h,v)}{\mathcal{X}_p} = \sqrt{\norm{h}{\mathcal{H}}^2 + \norm{v}{\mathcal{V}}^2}$. Let $\C_p=\mathcal{X}_p$, we can pick the entropy $\phi_p\para{x} = \frac{1}{2}\norm{x}{\mathcal{X}_p}^2$ whose gradient is given, for $x=\para{h,v}$, by $\nabla \phi_p\para{x} = h + J_{\mathcal{V}}v$ where $J_{\mathcal{V}}$ is the normalized duality mapping for $\mathcal{V}$. Then, if we assume that $\mathcal{V}$ is a smooth and rotund space as in \cite[Lemma~6.2]{BauschkeBorweinCombettes01}, $\mathcal{X}_p$ will be a smooth and rotund space and we will have that $\phi_p$ is Legendre, i.e., \ref{ass:legendre} will be satisfied. Since $\mathcal{V}$ is finite-dimensional and $\X_p$ is open, \ref{ass:gradphiconsistent} is satisfied for $\nabla\phi_p$. Indeed, the limit point $x^\infty$ cannot lie on the boundary since the boundary is empty while $\mathcal{V}$ being finite-dimensional guarantees the continuity of $J_{\mathcal{V}}$.
\end{example}
\end{subsection}

\subsection{Main estimations}
The following results constitute the main estimations that will be used in the convergence analysis of \algref{alg:ibpds}.

\begin{lemma}\label{estimate}
Recall the notation of \eqref{notation}. Assume that \ref{ass:hyp} and \ref{ass:legendre}-\ref{ass:littled} hold, then we have the following energy estimation. For every $w\eqdef\para{x,\mu}\in (\C_p\times\C_d) \cap (\dom(g) \times \dom(l^*))$\fekn
\begin{multline}\label{eq:mainestimpos}
\LL{x_{k+1},\mu}-\LL{x,\mu_{k+1}}+\left[\frac{1}{\Lambda_{k+1}} D\left(w,w_{k+1}\right)-M\left(w,w_{k+1}\right)\right] + \ip{w_{k+1}-w, \wnk}{} + \eps d\left(w_{k+1},w_k\right) \\
\leq \left[\frac{1}{\Lambda_k} D\left(w,w_k\right)-M\left(w,w_k\right)\right].
\end{multline}
If, moreover, \ref{ass:strconvex}\ref{ass:strconvexP} and \ref{ass:strconvex}\ref{ass:strconvexD} hold, we have (using the notation of \eqref{notationstr}) for every $w\eqdef\para{x,\mu}\in (\C_p\times\C_d) \cap (\dom(g) \times \dom(l^*))$\fekn
\begin{multline}\label{eq:mainestimpos2}
\LL{x_{k+1},\mu}-\LL{x,\mu_{k+1}}+\left[\frac{1}{\Lambda_{k+1}} D\left(w,w_{k+1}\right)-M\left(w,w_{k+1}\right)\right] + \ip{w_{k+1}-w, \wnk}{}+ \eps d\left(w_{k+1},w_k\right)\\
\leq \left[\frac{1}{\Lambda_k} D\left(w,w_k\right)-M\left(w,w_k\right)\right] - \scs D'\para{w,w_k}- \scns D'\para{w,w_{k+1}}.
\end{multline}
\end{lemma}

\begin{proof}
	We will prove claim \eqref{eq:mainestimpos2} since \eqref{eq:mainestimpos} is a special case of it when $\scs=\scns=0$. For all $k\in\N$, the following holds by the definitions of $x_{k+1}$ and $\mu_{k+1}$ in \algref{alg:ibpds},
	\nnewq{\label{eq:subdiff}
	\frac{1}{\lambda_k}\para{\nabla\phi_p\para{x_k} - \nabla \phi_p\para{x_{k+1}}} -\nabla f\para{x_k} -\pnk -T^*\wtilde{\mu}_k\in \partial g\para{x_{k+1}}\\
	\frac{1}{\nu_k}\para{\nabla\phi_d\para{\mu_k} - \nabla \phi_d\para{\mu_{k+1}}} -\nabla h^*\para{\mu_k} -\dnk +T\wtilde{x}_k\in \partial l^*\para{\mu_{k+1}}.
	}
	Observe that by assumptions \ref{ass:hyp} and \ref{ass:strconvex}\ref{ass:strconvexP}, we have $\C_p \cap \dom(g) = \dom(\phi_p) \cap \dom(g) \subset \dom(\psi_p) \cap \dom(g)$. Morover, using also that $\dom(\partial g) \subset \dom(g)$, we have $\forall k \in \N$, $x_k \in \inte \dom(\phi_p) \cap \dom(\partial g) = \inte \dom(\phi_p) \cap \dom(\partial g) \cap \C_p \cap \dom(g) \subset \inte \dom(\psi_p) \cap \dom(\partial g) \cap \C_p \cap \dom(g)$. A similar reasoning is also valid replacing $(\C_p,g,\phi_p,\psi_p)$ with their dual counterparts $(\C_d,l^*,\phi_d,\psi_d)$ and invoking \ref{ass:strconvex}\ref{ass:strconvexD}. We are then in position to apply the relative strong convexity inequality of \defref{genstrong}, which holds at any $(x,\mu) \in (\C_p\times\C_d) \cap (\dom(g) \times \dom(l^*))$ and $(x_{k+1},\mu_{k+1})$, hence giving
	\nnewq{\label{eq:relstrconvexgl}
	g\para{x} &\geq g\para{x_{k+1}} + \ip{u, x-\xkp}{} + \scmg D_{\psi_p}\para{x,\xkp}\\
	l^*\para{\mu}&\geq l\para{\mu_{k+1}} + \ip{v, \mu-\mu_{k+1}}{} + \scml D_{\psi_d}\para{\mu,\mu_{k+1}}
	}
	for any $u\in \partial g\para{x_{k+1}}$ and $v\in\partial l^*\para{\mu_{k+1}}$. Combining \eqref{eq:subdiff} and \eqref{eq:relstrconvexgl} and applying the three-point identity for Bregman divergences \cite[Lemma 3.1]{chen1993convergence}, we have
	\begin{equation}\label{est}
	\begin{split}
	D_p\left(x,x_k\right) \geq \lambda_k\BPa{g(x_{k+1})-g(x)+\langle \nabla f (x_k) +\pnk, \ x_{k+1}-x\rangle +\langle T\left(x_{k+1}-x\right), \ \wtilde{\mu}_k\rangle}\\
	+\scmg\lambda_kD_{\psi_p}\para{x,x_{k+1}}+D_p\left(x,x_{k+1}\right) + D_p\left(x_{k+1},x_k\right);\\
	D_d\left(\mu,\mu_k\right) \geq \nu_k\BPa{l^{*}(\mu_{k+1})-l^{*}(\mu)+\langle \nabla h^{*} (\mu_k)+\dnk, \ \mu_{k+1}-\mu\rangle - \langle T\wtilde{x}_k, \  \mu_{k+1}-\mu\rangle}\\
	+\scml\nu_kD_{\psi_d}\para{\mu,\mu_{k+1}}+D_d\left(\mu,\mu_{k+1}\right)+D_d\left(\mu_{k+1},\mu_k\right).
	\end{split}
	\end{equation}	
	Moreover, from the relative smoothness assumed in \ref{ass:legendre} and the consequent generalized descent lemma \eqref{gendl}, we have\fekn
	\begin{equation}\label{ineq1}
	\begin{split}
	f(x_{k+1}) & \leq f(x_k) + \langle \nabla f (x_k), \ x_{k+1}-x_k\rangle + L_p D_p\left(x_{k+1},x_k\right)\\
	h^{*} (\mu_{k+1}) & \leq h^{*} (\mu_k) + \langle \nabla h^{*} (\mu_k), \ \mu_{k+1}-\mu_k\rangle+L_d D_d\left(\mu_{k+1},\mu_k\right).
	\end{split}
	\end{equation}
	To apply the relative strong convexity inequality to $f$ and $h^*$, we again check the required qualification conditions of \defref{genstrong}. First, from \ref{ass:hyp} and \ref{ass:strconvex}\ref{ass:strconvexP}, $\C_p \cap \dom(g) \subset \C_p = \dom(f) \cap \dom(\phi_p) \subset \dom(f) \cap \dom(\psi_p)$. In addition, $\forall k \in \N$, $x_k \in \inte \dom(\phi_p) \subset \inte \dom(\psi_p)$. Since $f$ is differentiable on $\inte \dom(\phi_p)$, we have $\inte \dom(\phi_p) \subset \dom(\partial f)$, i.e, $x_k \in \inte \dom(\psi_p) \cap \dom(\partial f)$. We have also argued above that $x_k \in \inte \dom(\phi_p) \cap \dom(\partial g) \subset \C_p \cap \dom(g)$, and thus $x_k \in \inte \dom(\psi_p) \cap \dom(\partial f) \cap \C_p \cap \dom(g)$ as required to apply the relative strong convexity inequality of $f$ at any $x \in \C_p \cap \dom(g)$ and $x_{k+1}$. The same reasoning remains valid replacing $(\C_p,f,g,\phi_p,\psi_p)$ with $(\C_d,h^*,l^*,\phi_d,\psi_d)$ and invoking \ref{ass:strconvex}\ref{ass:strconvexD}. We then have for any $(x,\mu) \in (\C_p\times\C_d) \cap (\dom(g) \times \dom(l^*))$\fekn
	\begin{equation}\label{ineq2}
	\begin{split}
	f(x) & \geq f(x_k) + \langle \nabla f (x_k), \ x-x_k\rangle + \scmf D_{\psi_p}\left(x,x_k\right)\\
	h^{*} (\mu) & \geq h^{*} (\mu_k) + \langle \nabla h^{*} (\mu_k), \ \mu-\mu_k\rangle + \scmh D_{\psi_d}\left(\mu,x_k\right).
	\end{split}
	\end{equation}
	Summing \eqref{ineq1} and \eqref{ineq2}, we obtain, for each $(x,\mu)\in\C_p\times\C_d$\fekn
	\begin{equation*}
	\begin{split}
	f(x_{k+1}) & \leq f(x) + \langle \nabla f (x_k), \ x_{k+1}-x\rangle + L_p D_p\left(x_{k+1},x_k\right) - \scmf D_{\psi_p}\para{x,x_k}\\
	h^{*} (\mu_{k+1}) & \leq h^{*} (\mu) + \langle \nabla h^{*} (\mu_k), \ \mu_{k+1}-\mu\rangle+L_d D_d\left(\mu_{k+1},\mu_k\right) - \scmh D_{\psi_d}\para{\mu,\mu_k}.
	\end{split}
	\end{equation*}
	Summing the latter with \eqref{est}, we obtain
	\newq{
		\lambda_k\sbrac{f(x_{k+1})+g(x_{k+1})-f(x)-g(x)+\langle T\left(x_{k+1}-x\right), \  \wtilde{\mu}_k\rangle + \ip{x_{k+1}-x, \pnk}{}} +D_p\left(x,x_{k+1}\right)\\
		+ (1 -L_p\lambda_k)D_p\left(x_{k+1},x_k\right)\leq D_p\left(x,x_k\right) -\scmf\lambda_kD_{\psi_p}\para{x,x_k} -\scmg\lambda_k D_{\psi_p}\para{x,x_{k+1}};
	}
	\newq{
		\nu_k\sbrac{h^{*} (\mu_{k+1})+l^{*}(\mu_{k+1})-h^{*} (\mu)-l^{*}(\mu) - \langle T\wtilde{x}_k, \  \mu_{k+1}-\mu\rangle+\ip{\mu_{k+1}-\mu,\dnk}{}}+D_d\left(\mu,\mu_{k+1}\right)\\
		+(1 -L_d\nu_k)D_d\left(\mu_{k+1},\mu_k\right)	\leq D_d\left(\mu,\mu_k\right)-\scmh \nu_kD_{\psi_d}\para{\mu,\mu_k}-\scml\nu_kD_{\psi_d}\para{\mu,\mu_{k+1}}.
	}
	\ \\
	Recall the notations of \eqref{notation}, \eqref{notationstr}, and that 
\newq{
\iprod{w_1 - w_2}{\wnk} \eqdef \iprod{x_{1}-x_2}{\pnk} + \iprod{\mu_{1}-\mu_2}{\dnk},
}
	then, for each $\para{x,\mu}\in\C_p\times\C_d$\fekn
	\begin{multline*}
	\LL{x_{k+1},\mu}-\LL{x,\mu_{k+1}}+\langle T\left(x_{k+1}-x\right), \wtilde{\mu}_k\rangle - \langle T\wtilde{x}_k, \mu_{k+1}-\mu\rangle + \ip{w_{k+1}-w, \wnk}{}\\
	+\frac{1}{\Lambda_k} D\left(w,w_{k+1}\right)  -\frac{1}{\Lambda_k} D\left(w,w_k\right)  + \left(\frac{1}{\Lambda_k}-L\right) D\left(w_{k+1},w_k\right)\\
	\leq \langle T x_{k+1}, \mu \rangle - \langle T x, \mu_{k+1} \rangle- \scns D'\para{w,w_{k+1}}- \scs D'\para{w,w_k}.
	\end{multline*}	
	Rearranging the terms, we arrive at
	\begin{equation*}
	\begin{split}
	&\LL{x_{k+1},\mu}-\LL{x,\mu_{k+1}}+\frac{1}{\Lambda_k} D\left(w,w_{k+1}\right) - \frac{1}{\Lambda_k} D\left(w,w_k\right)	+ \left(\frac{1}{\Lambda_k}-L\right) D\left(w_{k+1},w_k\right) + \ip{w_{k+1} - w,\wnk}{}\\
	&\leq \langle T x_{k+1}, \mu-\wtilde{\mu}_k \rangle +\langle T (\wtilde{x}_k-x),\mu_{k+1} \rangle + \langle T x, \wtilde{\mu}_k \rangle-\langle T \wtilde{x}_k, \mu \rangle - \scns D'\para{w,w_{k+1}} - \scs D'\para{w,w_k}\\
%	= \ \ \ & \langle T x_{k+1}, \mu-\wtilde{\mu}_k \rangle +\langle T (\wtilde{x}_k-x),\mu_{k+1} \rangle + \langle T x, \wtilde{\mu}_k \rangle-\langle T (\wtilde{x}_k-x), \mu \rangle-\langle T x, \mu \rangle\\
%	= \ \ \ & \langle T x_{k+1}, \mu-\wtilde{\mu}_k \rangle +\langle T (\wtilde{x}_k-x),\mu_{k+1}-\mu \rangle - \langle T x, \mu-\wtilde{\mu}_k \rangle\\
	&= \langle T (x_{k+1}-x), \mu-\wtilde{\mu}_k \rangle +\langle T (\wtilde{x}_k-x),\mu_{k+1}-\mu \rangle - \scns D'\para{w,w_{k+1}} - \scs D'\para{w,w_k}.
	\end{split}
	\end{equation*}	
%	Finally, for each $\para{x,\mu}\in\C_p\times\C_d$\fekn
%	\begin{equation*}
%	\begin{split}
%	\LL{x_{k+1},\mu}-\LL{x,\mu_{k+1}}+\frac{1}{\Lambda_k} D\left(w,w_{k+1}\right) - \frac{1}{\Lambda_k}D\left(w,w_k\right)	+ \left(\frac{1}{\Lambda_k}-L\right) D\left(w_{k+1},w_k\right)+\ip{w_{k+1}-w,\wnk}{}\\
%	\leq \langle T (x_{k+1}-x), \mu-\wtilde{\mu}_k \rangle +\langle T (\wtilde{x}_k-x),\mu_{k+1}-\mu \rangle - \scns D'\para{w,w_{k+1}} - \scs D'\para{w,w_k}.
%	\end{split}
%	\end{equation*}	
	Now we use that $\wtilde{x}_k=2x_{k+1}-x_k$ and $\wtilde{\mu}_k=\mu_k$, to obtain
	\begin{equation*}
	\begin{split}
	&\LL{x_{k+1},\mu}-\LL{x,\mu_{k+1}}+\frac{1}{\Lambda_k} D\left(w,w_{k+1}\right) - \frac{1}{\Lambda_k} D\left(w,w_k\right)+ \left(\frac{1}{\Lambda_k}-L\right) D\left(w_{k+1},w_k\right) + \ip{w_{k+1}-w, \wnk}{}\\
	&\leq  \langle T (x_{k+1}-x), \mu-\mu_k \rangle +\langle T (x_{k+1}-x),\mu_{k+1}-\mu \rangle+\langle T (x_{k+1}-x_k),\mu_{k+1} -\mu \rangle- \scns D'\para{w,w_{k+1}}\\
	&\quad  - \scs D'\para{w,w_k}\\
%	= \ \ \ & \langle T (x_{k+1}-x),\mu_{k+1}-\mu_k \rangle+\langle T (x_{k+1}-x_k),\mu_{k+1}-\mu \rangle\\
%	= \ \ \ & \left[\langle T (x_{k+1}-x_k),\mu_{k+1}-\mu_k \rangle+\langle T (x-x_{k+1}),\mu-\mu_{k+1}\rangle-\langle T (x-x_{k}),\mu-\mu_k \rangle\right]\\
%	& -\left[\langle T (x_{k+1}-x_k),\mu_{k+1}-\mu_k \rangle+\langle T (x-x_{k+1}),\mu-\mu_{k+1}\rangle-\langle T (x-x_{k}),\mu-\mu_k \rangle\right]\\
%	&+\langle T (x_{k+1}-x),\mu_{k+1}-\mu_k \rangle+\langle T (x_{k+1}-x_k),\mu_{k+1}-\mu \rangle\\
%	= \ \ \ & \left[\langle T (x_{k+1}-x_k),\mu_{k+1}-\mu_k \rangle+\langle T (x-x_{k+1}),\mu-\mu_{k+1}\rangle-\langle T (x-x_{k}),\mu-\mu_k \rangle\right]\\
%	&+\langle T (x_{k}-x),\mu_{k+1}-\mu_k \rangle+\langle T (x-x_k),\mu_{k+1}-\mu \rangle+\langle T (x-x_{k}),\mu-\mu_k \rangle\\
%	= \ \ \ & \left[\langle T (x_{k+1}-x_k),\mu_{k+1}-\mu_k \rangle+\langle T (x-x_{k+1}),\mu-\mu_{k+1}\rangle-\langle T (x-x_{k}),\mu-\mu_k \rangle\right]\\
%	&+\langle T (x-x_k),\mu_k-\mu_{k+1} \rangle+\langle T (x-x_k),\mu_{k+1}-\mu \rangle+\langle T (x-x_{k}),\mu-\mu_k \rangle\\
%	= \ \ \ & \left[\langle T (x_{k+1}-x_k),\mu_{k+1}-\mu_k \rangle+\langle T (x-x_{k+1}),\mu-\mu_{k+1}\rangle-\langle T (x-x_{k}),\mu-\mu_k \rangle\right]\\
%	&+\langle T (x-x_k),\mu_{k}-\mu \rangle+\langle T (x-x_{k}),\mu-\mu_k \rangle\\
	&= \BPa{\langle T (x_{k+1}-x_k),\mu_{k+1}-\mu_k \rangle+\langle T (x-x_{k+1}),\mu-\mu_{k+1}\rangle-\langle T (x-x_{k}),\mu-\mu_k \rangle} - \scns D'\para{w,w_{k+1}}\\
	&\quad  - \scs D'\para{w,w_k}.
	\end{split}
	\end{equation*}
	Equivalently, recalling that $M(w_1,w_2) \eqdef \langle T (x_1-x_2),\mu_1-\mu_2 \rangle$, we get
	\begin{multline}\label{eq:mainestim}
	\LL{x_{k+1},\mu}-\LL{x,\mu_{k+1}}+ \ip{w_{k+1}-w,\wnk}{}+\left[\frac{1}{\Lambda_k}D\left(w,w_{k+1}\right)-M\left(w,w_{k+1}\right)\right]\\
	-\left[\frac{1}{\Lambda_k} D\left(w,w_k\right)-M\left(w,w_k\right)\right]+\left[\left(\frac{1}{\Lambda_k}-L\right) D\left(w_{k+1},w_k\right)-M\left(w_{k+1},w_k\right)\right]\\
	\leq - \scns D'\para{w,w_{k+1}} - \scs D'\para{w,w_k}.
	\end{multline}	
	Recall that, by \ref{ass:steps}, $\seq{\lambda_k}$ and $\seq{\nu_k}$ are nondecreasing sequences, and thus
	\nnewq{\label{telescop}
	\frac{1}{\Lambda_{k+1}}D\left(w,w_{k+1}\right)\leq \frac{1}{\Lambda_{k}}D\left(w,w_{k+1}\right).
	}
	Finally, combining \eqref{eq:mainestim} with \eqref{telescop} and \ref{ass:littled} applied at the points $w_{k+1}$ and $w_k$, we get \eqref{eq:mainestimpos2}.
\end{proof}
%-----------------------------------------------------------
\begin{lemma}\label{lem:posproduct}
Assume \ref{ass:hyp} and \ref{ass:legendre}-\ref{ass:littled} hold and\fekn denote by $\wex_{k+1}$ the exact update of the algorithm, i.e.,
\nnewq{\label{wex}
\wex_{k+1} = \pmat{\xex_{k+1} \\ \muex_{k+1}}=\pmat{\para{\nabla\phi_p + \lambda_k\p g}^{-1}\para{\nabla \phi_p\para{x_k} - \lambda_k\para{\nabla f\para{x_k}} -\lambda_kT^*\mu_k}\\
\para{\nabla\phi_d + \nu_k\p l^*}^{-1}\para{\nabla\phi_d\para{\mu_k} -\nu_k\para{\nabla h^*\para{\mu_k}} + \nu_kT\para{2\xex_{k+1}-x_k}}}.
}
Then, the following holds\fekn
\nnewq{\label{posprodest}
\ip{\wnk,\wex_{k+1}-w_{k+1}}{} \geq \frac{1}{\Lambda_k}\para{D\para{\wex_{k+1},w_{k+1}} + D\para{w_{k+1},\wex_{k+1}}}-2M\para{\wex_{k+1},w_{k+1}}\geq 0.
}
\end{lemma}
\begin{proof}
By design of the algorithm, the following monotone inclusions hold\fekn
\nnewq{\label{pmono}
\nabla \phi_p\para{x_k} - \lambda_k\para{\nabla f\para{x_k} -T^*\mu_k}-\nabla\phi_p\para{\xex_{k+1}}	&\in\lambda_k\p g\para{\xex_{k+1}}\\
\nabla \phi_p\para{x_k} - \lambda_k\para{\nabla f\para{x_k} +\pnk -T^*\mu_k}-\nabla\phi_p\para{x_{k+1}}	&\in\lambda_k\p g\para{x_{k+1}}.
}
and similarly for the dual
\nnewq{
\label{dmono}
\nabla \phi_d\para{\mu_k} - \nu_k\para{\nabla h^*\para{\mu_k} +T\para{2\xex_{k+1}-x_k}}-\nabla\phi_d\para{\muex_{k+1}}	&\in\nu_k\p l^*\para{\muex_{k+1}}\\
\nabla \phi_d\para{\mu_k} - \nu_k\para{\nabla h^*\para{\mu_k} + \dnk + T\para{2x_{k+1}-x_k}}-\nabla\phi_d\para{\mu_{k+1}}	&\in\nu_k\p l^*\para{\mu_{k+1}}.
}
By monotonicity of the operators $\p l^*$ and $\p g$ combined with \eqref{dmono} and \eqref{pmono}, we then have\fekn
\nnewq{
\ip{\xex_{k+1} - x_{k+1}, \pnk\lambda_k - \nabla \phi_p\para{\xex_{k+1}} +\nabla\phi_p\para{x_{k+1}}}{} &\geq 0\\
\ip{\muex_{k+1} - \mu_{k+1}, \dnk \nu_k -\nabla\phi_d\para{\muex_{k+1}} +\nabla\phi_d\para{\mu_{k+1}} +2\nu_k T\para{\xex_{k+1}-x_{k+1}}}{}	&\geq 0.
}
We can rewrite the above using \defref{def:bregman} to have\fekn
\nnewq{\label{earlybound}
\ip{\xex_{k+1} - x_{k+1}, \pnk}{}	&\geq \frac{1}{\lambda_k}\para{D_p\para{\xex_{k+1},x_{k+1}} + D_p\para{x_{k+1},\xex_{k+1}}}\\
\ip{\muex_{k+1}-\mu_{k+1}, \dnk}{}		&\geq \frac{1}{\nu_k}\para{ D_d\para{\muex_{k+1},\mu_{k+1}} + D_d\para{\mu_{k+1},\muex_{k+1}}} - 2\ip{T\para{\xex_{k+1}-x_{k+1}}, \muex_{k+1}-\mu_{k+1}}{}.
}
Adding the above inequalities together gives\fekn
\nnewq{\label{noisebound}
\ip{\wnk, \wex_{k+1}-w_{k+1}}{}\geq \frac{1}{\Lambda_k}\para{D\para{\wex_{k+1},w_{k+1}} + D\para{w_{k+1},\wex_{k+1}}}-2M\para{\wex_{k+1},w_{k+1}}.
}
Using \ref{ass:littled} and \eqref{aux}, and the fact that $M$ is symmetric wrt its arguments\fekn
\newq{
\frac{1}{\Lambda_k}\para{D\para{\wex_{k+1},w_{k+1}} + D\para{w_{k+1},\wex_{k+1}}}-2M\para{\wex_{k+1},w_{k+1}}\geq 0.
}
\end{proof}

\begin{lemma}\label{lem:entropysc}
Assume \ref{ass:hyp}, \ref{ass:legendre}-\ref{ass:littled}, and \ref{ass:error}\ref{strconvexerror} all hold. \iffalse that the entropies $\phi_p$ and $\phi_d$ are strongly convex with respect to $\norm{\cdot}{\X_p}^2$ and $\norm{\cdot}{\X_d}^2$ with modulus $\scphip$ and $\scphid$, respectively, and that the step size limits $\lambda_\infty$ and $\nu_\infty$ satisfy
\newq{
\nu_\infty\lambda_\infty < \frac{\scphip\scphid}{\norm{T}{p\to d^*}^2}.
}\fi
One can choose $a>0$ so that\fekn
\newq{
\frac{\scphip}{\lambda_k}-\frac{\norm{T}{p\to d^*}^2}{a}>0\qandq\frac{\scphid}{\nu_k}-a>0
}
and the following holds\fekn
\newq{
\ip{\wnk, \wex_{k+1}-w_{k+1}}{}\leq \para{\frac{\scphip}{\lambda_k}-\frac{\norm{T}{p\to d^*}^2}{a}}^{-1}\norm{\pnk}{\X_p^*}^2 + \para{\frac{\scphid}{\nu_k}-a}^{-1}\norm{\dnk}{\X_d^*}^2.
}
\end{lemma}
\begin{proof}
It follows from the strong convexity of $\phi_p$ and $\phi_d$ given by \ref{ass:error}\ref{strconvexerror} that\fekn
\nnewq{
\frac{1}{\Lambda_k}\para{D\para{w_{k+1}, \wex_{k+1}} + D\para{\wex_{k+1}, w_{k+1}}} &= \frac{1}{\Lambda_k}\ip{\nabla \phi \para{w_{k+1}} - \nabla \phi \para{\wex_{k+1}}, w_{k+1}-\wex_{k+1}}{}\\
 &\geq \frac{\scphip}{\lambda_k}\norm{\xex_{k+1}-x_{k+1}}{\X_p}^2 + \frac{\scphid}{\nu_k}\norm{\muex_{k+1}-\mu_{k+1}}{\X_d}^2.
}
Substituting this result into \lemref{lem:posproduct} \eqref{posprodest} and applying Young's inequality with $a>0$ we get\fekn
\nnewq{\label{ckwk}
&\ip{\Delta_k,\wex_{k+1}-w_{k+1}}{} 	\\
&\geq \frac{\scphip}{\lambda_k}\norm{\xex_{k+1}-x_{k+1}}{\X_p}^2 + \frac{\scphid}{\nu_k}\norm{\muex_{k+1}-\mu_{k+1}}{\X_d}^2 - 2M\para{\wex_{k+1},w_{k+1}}\\
									&= \frac{\scphip}{\lambda_k}\norm{\xex_{k+1}-x_{k+1}}{\X_p}^2 + \frac{\scphid}{\nu_k}\norm{\muex_{k+1}-\mu_{k+1}}{\X_d}^2 - 2\para{\ip{T\para{\xex_{k+1}-x_{k+1}}, \muex_{k+1}-\mu_{k+1}}{}}\\
									&\geq \frac{\scphip}{\lambda_k}\norm{\xex_{k+1}-x_{k+1}}{\X_p}^2 + \frac{\scphid}{\nu_k}\norm{\muex_{k+1}-\mu_{k+1}}{\X_d}^2 - \frac{\norm{T}{p\to d^*}^2}{a}\norm{\xex_{k+1}-x_{k+1}}{\X_p}^2-a\norm{\muex_{k+1}-\mu_{k+1}}{\X_d}^2\\
									&= \para{\frac{\scphip}{\lambda_k}-\frac{\norm{T}{p\to d^*}^2}{a}}\norm{\xex_{k+1}-x_{k+1}}{\X_p}^2 + \para{\frac{\scphid}{\nu_k}-a}\norm{\muex_{k+1}-\mu_{k+1}}{\X_d}^2.
}
Then, since the step size sequences $\seq{\lambda_k}$ and $\seq{\nu_k}$ are bounded and nondecreasing by \ref{ass:steps}, and furthemore by \ref{ass:error}\ref{strconvexerror} are chosen small enough to satisfy 
\newq{
\nu_\infty\lambda_\infty < \frac{\scphip\scphid}{\norm{T}{p\to d^*}^2},
}
one can choose $a>0$ so that
\newq{
\frac{\scphip}{\lambda_\infty}-\frac{\norm{T}{p\to d^*}^2}{a}>0\qandq\frac{\scphid}{\nu_\infty}-a>0
}
and, by extension under \ref{ass:steps}\fekn
\newq{
\frac{\scphip}{\lambda_k}-\frac{\norm{T}{p\to d^*}^2}{a}>0\qandq\frac{\scphid}{\nu_k}-a>0.
}
Finally, we apply Young's inequality twice to the following to find\fekn
\newq{
\ip{\wnk, \wex_{k+1}-w_{k+1}}{} &= \ip{\pnk, \xex_{k+1}-x_{k+1}}{} + \ip{\dnk, \muex_{k+1}-\mu_{k+1}}{}\\
&\leq \frac{1}{2}\para{\frac{\scphip}{\lambda_k}-\frac{\norm{T}{p\to d^*}^2}{a}}^{-1}\norm{\pnk}{\X_p^*}^2 + \frac{1}{2}\para{\frac{\scphip}{\lambda_k}-\frac{\norm{T}{p\to d^*}^2}{a}}\norm{\xex_{k+1}-x_{k+1}}{\X_p}^2\\
&\qquad + \frac{1}{2}\para{\frac{\scphid}{\nu_k}-a}^{-1}\norm{\dnk}{\X_d^*}^2 + \frac{1}{2}\para{\frac{\scphid}{\nu_k}-a}\norm{\muex_{k+1}-\mu_{k+1}}{\X_d}^2\\
&\leq \frac{1}{2}\para{\frac{\scphip}{\lambda_k}-\frac{\norm{T}{p\to d^*}^2}{a}}^{-1}\norm{\pnk}{\X_p^*}^2 + \frac{1}{2}\para{\frac{\scphid}{\nu_k}-a}^{-1}\norm{\dnk}{\X_d^*}^2\\
&\qquad +\frac{1}{2}\ip{\wnk,\wex_{k+1}-w_{k+1}}{}
}
and the desired claim follows.
\end{proof}
\begin{remark}
In \lemref{lem:entropysc}, one can instead choose to use $\norm{T^*}{d\to p^*}^2$ to have\fekn
\newq{
\ip{\Delta_k,\wex_{k+1}-w_{k+1}}{}\leq \para{\frac{\scphip}{\lambda_k}-\frac{1}{a}} \norm{\pnk}{\X_p^*}^2 + \para{\frac{\scphid}{\nu_k}-a\norm{T^*}{d\to p^*}^2}\norm{\dnk}{\X_d^*}^2
}
if there is asymmtery in the size of $\scphip$ and $\scphid$.

In the event that only $\phi_p$ is strongly convex with respect to $\norm{\cdot}{\X_p}^2$ but the analog does not hold for $\phi_d$, we can make the following argument. Take \eqref{earlybound} from \lemref{lem:posproduct} and use strong convexity, to get
\newq{
\ip{\xex_{k+1}-x_{k+1},\pnk}{}\geq \frac{1}{\lambda_k}\para{D_p\para{\xex_{k+1},x_{k+1}}+D_p\para{x_{k+1},\xex_{k+1}}}\geq \frac{\scphip}{\lambda_k}\norm{\xex_{k+1}-x_{k+1}}{\X_p}^2
}
and so, by Cauchy-Schwarz,
\newq{
\norm{\xex_{k+1}-x_{k+1}}{\X_p} \leq \frac{\lambda_k}{\scphip}\norm{\pnk}{\X_p^*}.
}
Then, using again Cauchy-Schwarz and the previous inequality,
\newq{
\ip{\pnk,\xex_{k+1}-x_{k+1}}{}\leq \norm{\pnk}{\X_p^*}\norm{\xex_{k+1}-x_{k+1}}{\X_p} \leq \frac{\lambda_k}{\scphip}\norm{\pnk}{\X_p^*}^2
}
without the restriction on $\lambda_\infty$ and $\nu_\infty$ imposed in \lemref{lem:entropysc} because we no longer to need control the term $2M\para{\wex_{k+1},w_{k+1}}$. This term, $2M\para{\wex_{k+1},w_{k+1}}$, is a result of the way we have defined $\muex_{k+1}$ to depend on $\xex_{k+1}$, which is necessary to keep $\wex_{k+1}$ deterministic conditioned on the filtration $\filts_k$. Thus, if only one of the entropies can be chosen to be strongly convex, one is inclined to formulate the problem in such a way that the primal problem has the strongly convex entropy, and to deal with the dual problem using \ref{ass:error}\ref{determinerror} or \ref{ass:error}\ref{boundederror} for the dual.
%\newq{
%\frac{1}{\lambda_k}\para{D_p\para{x_{k+1},\xex_{k+1}} + D_p\para{\xex_{k+1},x_{k+1}}}\geq \frac{\scphip}{\lambda_k}\norm{\xex_{k+1}-x_{k+1}}{\X_p}^2
%}
%with
%\newq{
%\ip{\pnk,\xex_{k+1}-x_{k+1}}{} \geq \frac{1}{\lambda_k}\para{D_p\para{\xex_{k+1},x_{k+1}} + D_p\para{x_{k+1},\xex_{k+1}}}
%}
%giving
%\newq{
%\ip{\pnk,\xex_{k+1}-x_{k+1}}{} \geq \frac{\scphip}{\lambda_k}\norm{\xex_{k+1}-x_{k+1}}{\X_p}^2
%}
%and so
%\newq{
%\ip{\pnk, \xex_{k+1}-x_{k+1}}{} &\leq \frac{\lambda_k}{2\scphip}\norm{\pnk}{\X_p^*}^2 + \frac{\scphip}{2\lambda_k}\norm{\xex_{k+1}-x_{k+1}}{\X_p}^2\\
%&\leq \frac{\lambda_k}{2\scphip}\norm{\pnk}{\X_p^*}^2 + \frac{1}{2}\ip{\pnk,\xex_{k+1}-x_{k+1}}{}\\
%\implies \ip{\pnk,\xex_{k+1}-x_{k+1}}{}&\leq \frac{\lambda_k}{\scphip}\norm{\pnk}{\X_p^*}^2
%}
\end{remark}

\begin{lemma}\label{lem:errorsum}
Under \ref{ass:hyp} and \ref{ass:legendre}-\ref{ass:error}, the following sequences satisfy, \faw
\newq{
\seq{\EX{\ip{\wnk,w-w_{k+1}}{}}{\filts_k}}\in\ell^1_+\para{\Filts}\qandq \seq{\TEX{\ip{\wnk,w-w_{k+1}}{}}}\in\ell^1_+.
}
\end{lemma}
\begin{proof}
The assumption \ref{ass:error} has three cases with the first, \ref{ass:error}\ref{determinerror}, corresponding to the deterministic setting, i.e., the lemma holds trivially. For both of the following two cases we note that, by \lemref{lem:posproduct}\fekn\faw
\nnewq{\label{sameprod}
\EX{\ip{\wnk, w-w_{k+1}}{}}{\filts_k} &= \EX{\ip{\wnk, w-\wex_{k+1}}{} + \ip{\wnk,\wex_{k+1}-w_{k+1}}{}}{\filts_k}\\
&= \EX{\ip{\wnk,\wex_{k+1}-w_{k+1}}{}}{\filts_k} \geq 0
}
since, due to \ref{ass:unbiased}, $\wnk$ is unbiased conditioned on the filtration $\filts_k$. By the law of total expectation applied to the above, it follows that\fekn\faw
\newq{
\TEX{\ip{\wnk, w-w_{k+1}}{}} = \TEX{\ip{\wnk,\wex_{k+1}-w_{k+1}}{}} \geq 0
}
and thus the following sequences satisfy, \faw
\newq{
\seq{\EX{\ip{\wnk,w-w_{k+1}}{}}{\filts_k}}\in\ell_+\para{\Filts}\qandq \seq{\TEX{\ip{\wnk,w-w_{k+1}}{}}}\in\ell_+.
}

Now assume that \ref{ass:error}\ref{boundederror} holds, recall that\fekn
\newq{
\ip{\wnk,\wex_{k+1}-w_{k+_1}}{} \eqdef \ip{\pnk,\xex_{k+1}-x_{k+1}}{} + \ip{\dnk, \muex_{k+1}-\mu_{k+1}}{}.
}
By \ref{ass:error}\ref{boundederror}, the sets $\U_p$ and $\U_d$ are bounded and thus have finite diameters, $\diam_{\U_p}$ and $\diam_{\U_d}$ respectively. Furthermore, by \ref{ass:legendre} and the definition of the updates in the algorithm, the exact update $\wex_{k+1}$ will remain in $\U_p\times \U_d$ for all $k\in\N$. Then\fekn
\newq{
\EX{\ip{\pnk, \xex_{k+1}-x_{k+1}}{}}{\filts_k} \leq \EX{\norm{\pnk}{\X_p^*}\norm{\xex_{k+1}-x_{k+1}}{\X_p}}{\filts_k} \leq \diam_{\U_p}\EX{\norm{\pnk}{\X_p^*}}{\filts_k};\\
\EX{\ip{\dnk, \muex_{k+1}-\mu_{k+1}}{}}{\filts_k} \leq \EX{\norm{\dnk}{\X_d^*}\norm{\muex_{k+1}-\mu_{k+1}}{\X_d}}{\filts_k} \leq \diam_{\U_d}\EX{\norm{\dnk}{\X_d^*}}{\filts_k}.
}
Since $\seq{\EX{\norm{\pnk}{\X_p^*}}{\filts_k}}\in\ell^1_+\para{\Filts}$ and $\seq{\EX{\norm{\dnk}{\X_d^*}}{\filts_k}}\in\ell^1_+\para{\Filts}$ by \ref{ass:error}\ref{boundederror}, and noting \eqref{sameprod}, it holds that, \faw
\newq{
\seq{\EX{\ip{\wnk,w-w_{k+1}}{}}{\filts_k}}\in\ell^1_+\para{\Filts}.
}
Using the same argument with the law of total expectation together with the fact that $\seq{\TEX{\norm{\pnk}{\X_p^*}}}\in\ell^1_+$ and $\seq{\TEX{\norm{\dnk}{\X_d^*}}}\in\ell^1_+$ by \ref{ass:error}\ref{boundederror}, it then follows that, \faw
\newq{
\seq{\TEX{\ip{\wnk,w-w_{k+1}}{}}}\in\ell^1_+.}

Finally, in the case of \ref{ass:error}\ref{strconvexerror}, we assume that the entropies $\phi_p$ and $\phi_d$ are strongly convex with respect to $\norm{\cdot}{\X_p}^2$ and $\norm{\cdot}{\X_d}^2$ respectively. Using \lemref{lem:entropysc} and taking expectation conditioned on $\filts_k$, we have\fekn
\newq{
\EX{\ip{\wnk,\wex_{k+1}-w_{k+1}}{}}{\filts_k} &\leq \para{\frac{\scphip}{\lambda_k} - \frac{\norm{T}{p\to d^*}^2}{a}}^{-1}\EX{\norm{\pnk}{\X_p^*}^2}{\filts_k} + \para{\frac{\scphid}{\nu_k} - a}^{-1}\EX{\norm{\dnk}{\X_d^*}^2}{\filts_k}\\
&\leq \para{\frac{\scphip}{\lambda_\infty} - \frac{\norm{T}{p\to d^*}^2}{a}}^{-1}\EX{\norm{\pnk}{\X_p^*}^2}{\filts_k} + \para{\frac{\scphid}{\nu_\infty} - a}^{-1}\EX{\norm{\dnk}{\X_d^*}^2}{\filts_k}.
}
Thus by the summability assumption of \ref{ass:error}\ref{strconvexerror}, we have
\newq{
\seq{\EX{\norm{\pnk}{\X_p^*}^2}{\filts_k}}\in\ell^1_+\para{\Filts}\qandq \seq{\EX{\norm{\dnk}{\X_d^*}^2}{\filts_k}}\in\ell^1_+\para{\Filts}
}
and so, \faw 
\newq{
\seq{\EX{\ip{\wnk,w-w_{k+1}}{}}{\filts_k}}\in\ell^1_+\para{\Filts}.
}
Similarly, taking \lemref{lem:entropysc} with total expectation and the summability assumption of \ref{ass:error}\ref{strconvexerror} yields, \faw
\newq{
\seq{\TEX{\ip{\wnk,w-w_{k+1}}{}}}\in\ell^1_+.
}
\end{proof}
%-----------------------
% Footer to allow Gummi
% to maintain previews
% for individual pieces
%-----------------------
\ifdefined\COMPLETE
\else
\end{document}
\fi

%----------------------------------
% These commands allow Gummi to
% edit each piece individually
% and still maintain previews.
%----------------------------------
\ifdefined\COMPLETE
\else
\documentclass[12pt]{article}
\input{tex_package_header} %file containing all the used libraries
\begin{document}
\fi
%----------------------------------
% Remove at your own risk!
% There is also a footer that ends
% the document if the main wasn't
% loaded.
%----------------------------------
\section{Convergence Analysis}\label{sec:conv}
\subsection{Ergodic Convergence}
Define\fekn the \emph{ergodic} iterates $\bar{x}_{k}\eqdef\dfrac{1}{k} \sum\limits_{i=1}^k x_i$ and $\bar{\mu}_k\eqdef\dfrac{1}{k}\sum\limits_{i=1}^k\mu_i$. The following theorem characterizes the convergence of the algorithm for the Lagrangian optimality gap evaluated at the ergodic iterates; later results on pointwise convergence will also imply ergodic convergence.\footnote{By ``ergodic convergence'', we mean convergence of the Lagrangian optimality gap evaluated at the ergodic iterates; not any ergodic averaging of the Lagrangian values themselves.}
\begin{theorem}\label{thm:ergodic}
Let \ref{ass:hyp} and \ref{ass:legendre}-\ref{ass:unbiased} hold. Then we have the following convergence rate to a noise-dominated regime: for each $k\in\N$, for every $\left(x,\mu\right)\in(\C_p\times\C_d) \cap (\dom(g) \times \dom(l^*))$,
\begin{align}\label{eq:rateL}
		\mathbb{E}\left[\LL{\bar{x}_{k},\mu}-\LL{x,\bar{\mu}_{k}} \right]
		&\leq \frac{\frac{1}{\Lambda_0} D\left(w,w_0\right)-M\left(w,w_0\right)}{k}+\frac{\sum_{i=0}^{k-1}\TEX{\ip{\Delta_i,w-w_{i+1}}{}}}{k}.
		%&\leq \frac{\Ppa{\frac{1}{\lambda_0}+(1+\alpha_p^{-1})\norm{T}{p}} D_p\left(x,x_0\right)+\Ppa{\frac{1}{\nu_0}+(1+\alpha_d^{-1})\norm{I}{d}} D_d\left(\mu,\mu_0\right)}{k} .
\end{align}
In particular, if also \ref{ass:error} holds, every almost sure weak sequential cluster point of $\seq{\bar{w}_k}$ is optimal in mean; if $\bar{w}_{k_j} \rightharpoonup w_{\infty}$ almost surely, then $\mathbb{E}(w_{\infty})$ is a saddle point for the Lagrangian.
\end{theorem}

\begin{proof}
Let $w\eqdef\para{x,\mu}\in(\C_p\times\C_d) \cap (\dom(g) \times \dom(l^*))$. Beginning with \lemref{estimate}, 
%we have for every $w\eqdef\para{x,\mu}\in(\C_p\times\C_d) \cap (\dom(g) \times \dom(l^*))$\fekn
%	\begin{multline}\label{eq:mainestimposs}
%	\LL{x_{k+1},\mu}-\LL{x,\mu_{k+1}}+\left[\frac{1}{\Lambda_{k+1}} D\left(w,w_{k+1}\right)-M\left(w,w_{k+1}\right)\right]\\
%	+ \eps d\left(w_{k+1},w_k\right) \leq \left[\frac{1}{\Lambda_k} D\left(w,w_k\right)-M\left(w,w_k\right)\right] + \ip{\wnk,w-w_{k+1}}{}.
%	\end{multline}	
	taking the the total expectation of \eqref{eq:mainestimpos} and summing up from $0$ to $k-1$, discarding positive terms on the left hand side, we have
\begin{equation}
	\begin{split}
	\sum\limits_{i=0}^{k-1}\TEX{\LL{x_{i+1},\mu}-\LL{x,\mu_{i+1}}} &\leq \frac{1}{\Lambda_0} D\left(w,w_0\right)-M\left(w,w_0\right) + \sum\limits_{i=0}^{k-1}\TEX{\ip{\Delta_i,w-w_{i+1}}{}}.
	% &\leq \frac{1}{\Lambda_0} D\left(w,w_0\right)-M\left(w,w_0\right) + \sum\limits_{i=0}^{\infty}\TEX{\ip{\Delta_i,w-w_{i+1}}{}}.
	\end{split}
\end{equation}
	 Notice that $\sum\limits_{i=0}^{k-1}\TEX{\ip{\Delta_i,w-w_{i+1}}{}}$ is nonnegative by \ref{ass:unbiased} and \lemref{lem:posproduct}. Using Jensen's inequality with the convex-concave function $\mathcal{L}$, we have \eqref{eq:rateL}.\\
	 
	 Now, assuming also \ref{ass:error}, let $(\bar{x}_{k_j},\bar{\mu}_{k_j})\rightharpoonup(x_{\infty},\mu_{\infty})$ almost surely. First note that, 
	 by \lemref{lem:errorsum}, 
	 \newq{
	 \sum\limits_{i=0}^\infty \TEX{\ip{\Delta_i,w-w_{i+1}}{}}<+\infty.
	 } 
	 Then, for every $\para{x,\mu}\in(\C_p\times\C_d) \cap (\dom(g) \times \dom(l^*))$,
\begin{equation}\label{eq:ergodicsaddle}
\begin{split}
\LL{\mathbb{E}\left(x_{\infty}\right),\mu}-\LL{x,\mathbb{E}\left(\mu_{\infty}\right)}& \leq \mathbb{E}\left[\LL{x_{\infty},\mu}-\LL{x,\mu_{\infty}}\right]\\
& \leq \mathbb{E} \left[\liminf_{j \to \infty} \ \left[\LL{\bar{x}_{k_j},\mu}-\LL{x,\bar{\mu}_{k_j}}\right]\right]\\
& \leq \liminf_{j \to \infty} \ \mathbb{E} \left[\LL{\bar{x}_{k_j},\mu}-\LL{x,\bar{\mu}_{k_j}}\right] \\
& \leq 0,
\end{split}
\end{equation}
where we used Jensen's inequality, weak lower semicontinuity of $\mathcal{L}$, Fatou's Lemma and \eqref{eq:rateL} with \ref{ass:error} and \lemref{lem:errorsum}. Inequality \eqref{eq:ergodicsaddle} trivially holds outside $(\C_p\times\C_d) \cap (\dom(g) \times \dom(l^*))$, and so holds for any $\para{x,\mu} \in \X_p \times \X_d$, whence we get that $(\mathbb{E}\left(x_{\infty}\right),\mathbb{E}\left(\mu_{\infty}\right))$ is a saddle point for $\mathcal{L}$.
\end{proof}

\begin{remark}\label{rem:boundederror}
The term $k^{-1}\sum_{i=0}^{k-1}\TEX{\ip{\Delta_i,w-w_{i+1}}{}}$ in \thmref{thm:ergodic} is an averaging of the noise which dictates the radius of the noise-dominated region in some sense. For example, if we assume that there exists a constant $c\geq 0$ such that $\TEX{\ip{\Delta_i,w-w_{i+1}}{}} \leq c$ for all $i\in\N$ and for all $w\in\X_p\times\X_d$, then we have
\newq{
\frac{\sum_{i=0}^{k-1}\TEX{\ip{\Delta_i,w-w_{i+1}}{}}}{k}\leq c
}
for all $k\in\N$, i.e., the radius of the noise-dominated region in \thmref{thm:ergodic} is at most $c$.
\end{remark}

\begin{remark}
Consider the algorithm in the deterministic case, then choose $\para{x,\mu} = \para{\xs,\mus}$ for some saddle point $\para{\xs,\mus}$ in \eqref{eq:rateL}. In this case, the constant in the rate of convergence, $\frac{1}{\Lambda_0}D\para{\ws,w_0}-M\para{\ws,w_0}$, is given in terms of the Bregman divergence, in contrast to methods like \cite{chambolle2011first} which have constants in terms of the Euclidean norm. With this change in the geometry, the dependence of the constant on the dimension of the problem can be greatly reduced, even from linear to logarithmic dependence for some problems and appropriately chosen entropies.
\end{remark}

\subsection{Asymptotic Regularity}

\begin{theorem}\label{asreg}
Let \ref{ass:hyp}, \ref{ass:legendre}, \ref{ass:steps}, \ref{ass:littled}\ref{littledstrong}, \ref{ass:unbiased}, \ref{ass:error}, and \ref{ass:dconsistent} hold. Then the primal-dual sequence $\seq{x_k,\mu_k}$ is almost surely asymptotically regular, meaning that $x_{k+1}-x_k\to 0$ and $\mu_{k+1}-\mu_k\to 0$ almost surely.
\end{theorem}

\begin{proof}
	Use again \eqref{eq:mainestimpos} in \lemref{estimate} with $w$ equal to a saddle point $\ws\in\sadp$ and take the total expectation to get, for each $k\in\N$,
	\begin{multline}\label{ineq}
	\TEX{\LL{x_{k+1},\mus}-\LL{\xs,\mu_{k+1}}}+\TEX{\frac{1}{\Lambda_{k+1}} D\left(\ws,w_{k+1}\right)-M\left(\ws,w_{k+1}\right)} \\
	+ \eps \TEX{d\left(w_{k+1},w_k\right)} \leq \TEX{\frac{1}{\Lambda_k} D\left(\ws,w_k\right)-M\left(\ws,w_k\right)}+\TEX{\ip{\wnk, \ws-w_{k+1}}{}}.
	\end{multline}
	By the definition of saddle point in \eqref{sadp}, it holds\fekn
	\newq{
	\LL{x_{k},\mus}-\LL{\xs,\mu_{k}}\geq 0
	}
	and so, from \lemref{LemmaAux} with \ref{ass:unbiased}, \ref{ass:error}, \lemref{lem:errorsum}, and \ref{ass:littled}\ref{littledstrong}, 
	\newq{
	\begin{split}
	 \mathbb{E}\left[d\left(w_{k+1},w_k\right)\right]\in\ell^1_+.
	\end{split}
	}
	By \lemref{AuxAS}, $d\para{w_{k+1},w_k}\to 0$ almost surely. In view of \ref{ass:dconsistent}, we get that, almost surely,
	\begin{equation}\label{eq:asympreg}
	w_{k+1}-w_k \to 0,
	\end{equation}
	i.e., the primal-dual sequence $\seq{w_k}$ is almost surely asymptotically regular.
\end{proof}

\subsection{Pointwise Convergence}\label{sec:pointwise}
The main result of this section is related to the pointwise weak convergence of the primal-dual sequence $\seq{x_k,\mu_k}$ to a saddle point. These results require the stronger assumptions \ref{ass:gradreg}-\ref{ass:Tconsistent}, although they are verified in many situations (see the discussion in \remref{rem:finitedim} and example thereafter). We will also impose the following conditions, which are only necessary for this particular section in the stochastic case and can be dropped for the deterministic case or the other sections. 
\begin{enumerate}[label=\bf{(\subscript{\mathbf{PW}}{{\arabic*}})}]
\item \label{ass:separable}$\X_p$ and $\X_d$ are separable.% (this holds, in particular, if the space $\X_p \times \X_d$ itself is separable).
\item \label{ass:property}The Bregman divergence $D$ satisfies the following property. Let $\tilde{\Omega}$ be a full-measure subset of $\Omega$ ($\tilde{\Omega}\in\sigalg$ with $\mathbb{P}(\tilde{\Omega})=1$). Let $\ws\in\sadp$ and $\seqn{s_n}\subset\sadp$ such that $s_n\to \ws$. If, for every $n\in\N$ and for every $\omega\in\tilde{\Omega}$,
\newq{
\lim\limits_{k\to\infty}\Lambda_k^{-1}D\para{s_n,w_k\para{\omega}}-M\para{s_n,w_k\para{\omega}}=r_{s_n}\para{\omega}\in[0,\pinfty[,
}
then there exists a $[0,\pinfty[$-valued random variable $r_{\ws}$ such that, for any $\omega\in\tilde{\Omega}$,
\newq{
\lim\limits_{k\to\infty}\Lambda_k^{-1}D\para{\ws,w_k\para{\omega}}-M\para{\ws,w_k\para{\omega}}=r_{\ws}\para{\omega}.
}
\end{enumerate}
\begin{proposition}\label{wkclstropt}
Let \ref{ass:hyp}, \ref{ass:legendre}, \ref{ass:steps}, \ref{ass:littled}\ref{littledweak} or \ref{littledstrong}, and \ref{ass:unbiased}-\ref{ass:gradreg} hold. Then $\seq{(x_k,\mu_k)}$ is almost surely bounded and, recalling the notation of \eqref{eq:wkset} and \eqref{eq:solset}, $\wkset{w}{k}\subset\sadp$ $\Pasp$.
\end{proposition}
\begin{proof}
Evaluating \lemref{estimate} at a saddle point $w=\ws\in\sadp$ and taking expectation conditioned on the filtration $\filts_k$, we get\fekn
\begin{multline*}%\label{eq:mainestimpos}
\EX{\LL{x_{k+1},\mus}-\LL{\xs,\mu_{k+1}}}{\filts_k}+\EX{\frac{1}{\Lambda_{k+1}} D\left(\ws,w_{k+1}\right)-M\left(\ws,w_{k+1}\right)}{\filts_k} \\
+ \eps \EX{d\left(w_{k+1},w_k\right)}{\filts_k} \leq \left[\frac{1}{\Lambda_k} D\left(\ws,w_k\right)-M\left(\ws,w_k\right)\right]+\EX{\ip{\wnk,\ws-w_{k+1}}{}}{\filts_k}.
\end{multline*}
Then, by \ref{ass:unbiased}, \ref{ass:error}, \lemref{lem:errorsum}, and \lemref{AuxPC}, $\seq{\Lambda_k^{-1}D\left(\ws,w_k\right)-M\left(\ws,w_k\right)}$ is almost surely convergent to some $r\in[0,+\infty[$. In particular, from \ref{ass:littled} and \eqref{aux}, both $\seq{D\left(\ws,w_k\right)}$ and $\seq{d\left(\ws,w_k\right)}$ are almost surely bounded and the coercivity condition \ref{ass:coercive} entails that the sequence $\seq{w_k}$ is almost surely bounded in $\inte\dom\para{\phi}$. Since $\X_p$ and $\X_d$ are reflexive, $\wkset{w}{k} \neq \emptyset$ almost surely. Let $w_{\infty}=(x_{\infty},\mu_{\infty})$ be an almost sure weak sequential cluster point of $\seq{w_k}$, i.e., there is a subsequence $\subseqi{w_{k_i}}$ such that $w_{k_i}\rightharpoonup w_{\infty}$ almost surely. The updates of \algref{alg:ibpds} are equivalent to the following monotone inclusions,
\begin{small}
\begin{multline}
\label{eq:inclusionPD}
\begin{pmatrix}
\dfrac{\nabla\phi_p(x_{k_i})-\nabla\phi_p(x_{k_{i+1}})}{\lambda_k} + (\nabla f(x_{k_{i+1}})-\nabla f(x_{k_i})-\delta^p_{k_i})+T^*(\mu_{k_{i+1}}-\mu_{k_i}) \\
\dfrac{\nabla\phi_d(\mu_{k_i})-\nabla\phi_d(\mu_{k_{i+1}})}{\nu_k} + (\nabla h^*(\mu_{k_{i+1}})-\nabla h^*(\mu_{k_i})-\delta^d_{k_i})+T(x_{k_{i+1}}-x_{k_i})
\end{pmatrix}\\
\in 
\begin{pmatrix} \partial g + \nabla f & 0 \\ 0 & \partial l^* + \nabla h^* \end{pmatrix}
\begin{pmatrix} x_{k_{i+1}} \\ \mu_{k_{i+1}} \end{pmatrix}
+
\begin{pmatrix} 0& T^* \\ -T & 0 \end{pmatrix}
\begin{pmatrix} x_{k_{i+1}} \\ \mu_{k_{i+1}} \end{pmatrix} .
\end{multline} 
\end{small}
Since $\subseqi{w_{k_i}}$ lies in $\inte \C_p \times \inte \C_d$, we have $N_{\C_p}(x_{k_{i+1}})=0$ and $N_{\C_d}(\mu_{k_{i+1}})=0$. This together with \cite[Theorem~2.4.2(viii)]{Zalinescu02} implies
\begin{small}
\begin{multline}
\label{eq:operatorPD}
\begin{pmatrix} \partial g + \nabla f & 0 \\ 0 & \partial l^* + \nabla h^* \end{pmatrix}
\begin{pmatrix} x_{k_{i+1}} \\ \mu_{k_{i+1}} \end{pmatrix}
+
\begin{pmatrix} 0& T^* \\ -T & 0 \end{pmatrix}
\begin{pmatrix} x_{k_{i+1}} \\ \mu_{k_{i+1}} \end{pmatrix} \subset \\
\begin{pmatrix} \partial (g + f + \iota_{\C_p}) & 0 \\ 0 & \partial (l^* + h^* + \iota_{\C_d})\end{pmatrix}
\begin{pmatrix} x_{k_{i+1}} \\ \mu_{k_{i+1}} \end{pmatrix}
+
\begin{pmatrix} 0& T^* \\ -T & 0 \end{pmatrix}
\begin{pmatrix} x_{k_{i+1}} \\ \mu_{k_{i+1}} \end{pmatrix} .
\end{multline} 
\end{small}
The first operator on the right hand side of \eqref{eq:operatorPD} is maximal monotone thanks to \ref{ass:hyp} and \cite[Theorem~3.1.11]{Zalinescu02}. The second operator is a skew-symmetric linear operator which is then maximal monotone with full domain by \cite[Section~17]{Simons08}. By \cite[Theorem 24.1(a)]{Simons08}, we deduce that the operator in the right hand side of \eqref{eq:operatorPD} is maximal monotone. Hence its graph is sequentially closed in the weak-strong topology by \cite[Lemma~1.2]{Browder65}. Recall that, by \ref{ass:unbiased}, \ref{ass:error}, and \remref{remark_errorconv}, $\seq{\pnk}$ and $\seq{\dnk}$ converge strongly to zero almost surely. From \thmref{asreg} and the fact that $w_{k_i}\rightharpoonup w_\infty$, we have also that $\subseqi{(x_{k_i + 1},\mu_{k_i + 1})}$ converges weakly to $(x_{\infty},\mu_{\infty})$ almost surely. In addition, by \ref{ass:hyp}, $T$ is linear (and bounded) which, combined with \thmref{asreg}, yields 
\newq{
T(x_{k_{i+1}}-x_{k_i}) \to 0\qandq T^*(\mu_{k_{i+1}}-\mu_{k_i}) \to 0
}
almost surely. From \ref{ass:gradreg} combined with \thmref{asreg}, we deduce that, almost surely,
\begin{align*}
\nabla\phi_p(x_{k_{i+1}})-\nabla\phi_p(x_{k_i}) \to 0 &\qandq \nabla f(x_{k_{i+1}})-\nabla f(x_{k_i}) \to 0 \\
\nabla\phi_d(\mu_{k_{i+1}})-\nabla\phi_d(\mu_{k_i}) \to 0 &\qandq \nabla h^*(\mu_{k_{i+1}})-\nabla h^*(\mu_{k_i}) \to 0 .
\end{align*}
Now since both $\seq{\lambda_k}$ and $\seq{\nu_k}$ are bounded away fro zero by \ref{ass:steps}, we have shown that, almost surely, the left hand side of \eqref{eq:inclusionPD} converge strongly. Hence, by weak-strong sequential closedness of the graph of the operator in \eqref{eq:operatorPD} we have shown above, we get
	\begin{align*}
		\begin{pmatrix}
			0 \\
			0
		\end{pmatrix}
		&\in
		\begin{pmatrix} \partial(g + f + \iota_{\C_p}) & T^* \\ -T & \partial(l^* + h^* + \iota_{\C_d}) \end{pmatrix}
		\begin{pmatrix} x_{\infty} \\ \mu_{\infty} \end{pmatrix},
	\end{align*} 
holds almost surely, whence it follows that each weak sequential cluster point of $\seq{w_k}$ is a saddle point almost surely.
\end{proof}

The significance of the following proposition is in the order of the quantifiers; it guarantees that there exists a full-measure set $\tilde{\Omega}$ for which the conclusion holds for every solution $\ws$.
\begin{proposition}\label{limexists}
	Let \ref{ass:hyp}, \ref{ass:legendre}, \ref{ass:steps}, \ref{ass:littled}\ref{littledweak} or \ref{littledstrong}, and \ref{ass:unbiased}-\ref{ass:gradreg} hold as well as \ref{ass:separable} and \ref{ass:property}. Then, there exists $\tilde{\Omega}\in\sigalg$ such that $\prob\para{\tilde{\Omega}}=1$ and, for every $\ws\in\sadp$ and for every $\omega\in\tilde{\Omega}$, the sequence  $$\seq{\Lambda_k^{-1}D\para{\ws,w_k\para{\omega}}-M\para{\ws,w_k\para{\omega}}}$$ converges with limit in $[0,\pinfty[$.
\end{proposition}
\begin{proof}
By \ref{ass:separable}, there exists a countable set $S$ such that $\clos{S}=\sadp$.  Once again, as in the proof of \propref{wkclstropt}, for every $\ws\in\sadp$ there exist $\Omega_{\ws}\in\sigalg$ such that $\prob\para{\Omega_{\ws}}=1$ and, for every $\omega\in\Omega_{\ws}$, it holds
\newq{
\lim\limits_{k\to\infty}\Lambda_k^{-1}D\para{\ws,w_k\para{\omega}}-M\para{\ws,w_k\para{\omega}} = r^\star\para{\omega}\in[0,\pinfty[.
}
Let $\tilde{\Omega} = \bigcap\limits_{s\in S}\Omega_{s}$ and notice that $\prob\para{\tilde{\Omega}}=1$ since, by countability of $S$, we have
\newq{
\prob\para{\tilde{\Omega}}=1-\prob\para{\tilde{\Omega}^c}=1-\prob\para{\bigcup\limits_{s\in S}\Omega_{s}^c} \geq 1-\sum\limits_{s\in S}\prob\para{\Omega_{s}^c}=1.
}
Fix a particular $\ws\in \sadp$; since $\clos{S}=\sadp$, there exists a sequence $\seqn{s_n}$ in $S$ such that $s_n\to \ws$. At the same time, for each $n\in\N$, there exists $r_n$, a $[0,\pinfty[$-valued random variable such that, for each $\omega\in\tilde{\Omega}$,
\newq{
\lim\limits_{k\to\infty}\Lambda_k^{-1}D\para{s_n,w_k\para{\omega}}-M\para{s_n,w_k\para{\omega}}=r_n\para{\omega}\in[0,\pinfty[.
}
Applying now \ref{ass:property}, we find that, for any $\omega\in\tilde{\Omega}$, \newq{
\lim\limits_{k\to\infty} \Lambda_k^{-1}D\para{\ws,w_k\para{\omega}} - M\para{\ws,w_k\para{\omega}} = r_{\ws}\para{\omega} \in[0,\infty[.
}
\end{proof}

\begin{theorem}\label{wkconvergence}
Let \ref{ass:hyp}, \ref{ass:legendre}, \ref{ass:steps}, and \ref{ass:unbiased}-\ref{ass:Tconsistent} hold as well as \ref{ass:separable} and \ref{ass:property}. Suppose also that one of the following holds:
\begin{enumerate}[label=(\roman*)]
\item $\sadp$ is a singleton. \label{assum:wkconvergence1}
\item \ref{ass:littled}\ref{littledweak}, $\sadp \subset \inte \C_p \times \inte \C_d$ and $\phi_p$ and $\phi_d$ are Legendre. \label{assum:wkconvergence2}
\item \ref{ass:littled}\ref{littledstrong}, and $d(w^1,w^2)=0 \Rightarrow w^1=w^2$. \label{assum:wkconvergence3}
\end{enumerate}
Then, there exists $\bar{w}$, a $\sadp$-valued random variable, such that $\seq{w_k}\rightharpoonup \bar{w}$ $\Pasp$.
\end{theorem}
\begin{proof}
We use a standard reasoning inspired by Opial's lemma (see \cite{opial1967weak}). We recall the notation of \eqref{eq:wkset} for the set of weak cluster points of a sequence. By \propref{wkclstropt}, there exists $\Omega'\in\sigalg$ with $\prob\para{\Omega'}=1$ such that, for any $\omega\in\Omega'$, the following holds
	\newq{
	\mathfrak{W}\sbrac{\para{w_k\para{\omega}}}\subset \sadp}
	and the sequence $\seq{w_k\para{\omega}}$ is bounded, and thus $\mathfrak{W}\sbrac{\para{w_k\para{\omega}}} \neq \emptyset$ since the spaces are reflexive.
	Furthermore, by \propref{limexists}, there exists $\Omega''\in\sigalg$ with $\prob\para{\Omega''}=1$ such that, for any $\omega\in\Omega''$, for any $\ws\in\sadp$, it holds
	\newq{
	\lim\limits_{k\to\infty}\Lambda_k^{-1}D\para{\ws,w_k\para{\omega}}-M\para{\ws,w_k\para{\omega}} = r_{\ws}\para{\omega}\in[0,\pinfty[.
	}
	Let $\tilde{\Omega}=\Omega'\cap\Omega''$, for any $\omega\in\tilde{\Omega}$ we let $w^1\para{\omega}\in\mathfrak{W}\sbrac{\seq{w_k\para{\omega}}}$ and $w^2\para{\omega}\in\mathfrak{W}\sbrac{\seq{w_k\para{\omega}}}$ be two weak sequential cluster points of $\seq{w_k\para{\omega}}$, i.e., there exists two subsequences $\subseqi{w_{k_i}\para{\omega}}$ and $\subseq{w_{k_j}\para{\omega}}$ such that $w_{k_i}\para{\omega}\rightharpoonup w^1\para{\omega}$ and $w_{k_j}\para{\omega}\rightharpoonup w^2\para{\omega}$ almost surely. Since $\mathfrak{W}\sbrac{\seq{w_k\para{\omega}}}\subset\sadp$, $w^1\para{\omega}$ and $w^2\para{\omega}$ are saddle points.% We now prove the existence of a set $\tilde{\Omega}\subset\sigalg$ such that $\prob\para{\tilde{\Omega}}=1$ and, for all $\omega\in\tilde{\Omega}$, for any solution $\ws$,
%\newq{
%\lim\limits_{k\to\infty}\para{\Lambda_k^{-1}D\para{\ws,w_k}-M\para{\ws,w_k}}=r_{\ws}. 
%}
Thus, there exist $r_{w^1}\para{\omega},r_{w^2}\para{\omega}\in [0,+\infty[$ such that,
\newq{
\lim_{k \to \infty} \Ppa{\Lambda_k^{-1}D(w^1\para{\omega},w_k\para{\omega})-M(w^1\para{\omega},w_k\para{\omega})} = r_{w^1}\para{\omega}}
and
\newq{
\lim_{k \to \infty} \Ppa{\Lambda_k^{-1}D(w^2\para{\omega},w_k\para{\omega})-M(w^2\para{\omega},w_k\para{\omega})} = r_{w^2}\para{\omega}.
}
Using the three point identity, we have, for each $i\in\N$,
\nnewq{\label{threepointguy}
&\Lambda_{k_i}^{-1}D(w^1\para{\omega},w_{k_i}\para{\omega})-M(w^1\para{\omega},w_{k_i}\para{\omega}) - \Lambda_{k_i}^{-1}D(w^2\para{\omega},w_{k_i}\para{\omega})+M(w^2\para{\omega},w_{k_i}\para{\omega}) \\
&= \Lambda_{k_i}^{-1}\Ppa{D(w^1\para{\omega},w_{k_i}\para{\omega}) - D(w^2\para{\omega},w_{k_i}\para{\omega})} - \Ppa{M(w^1\para{\omega},w_{k_i}\para{\omega})-M(w^2\para{\omega},w_{k_i}\para{\omega})} \\
&= \Lambda_{k_i}^{-1}\Ppa{D(w^1\para{\omega},w^2\para{\omega})-\ip{\nabla \phi(w_{k_i}\para{\omega})-\nabla \phi(w^2\para{\omega}), w^1\para{\omega}-w^2\para{\omega}}{}} \\
&\quad\quad\quad\quad\quad\quad\quad\quad\quad\quad\quad\quad\quad\quad\quad\quad\quad\quad\quad\quad\quad- \Ppa{M(w^1\para{\omega},w_{k_i}\para{\omega})-M(w^2\para{\omega},w_{k_i}\para{\omega})} .
}
Recall that, by \ref{ass:steps}, both $\seq{\lambda_k}$ and $\seq{\nu_k}$ are nondecreasing and bounded above with limits $\lambda_{\infty}$ and $\nu_{\infty}$, respectively. We denote $\Lambda_{\infty}\eqdef\left(\lambda_{\infty},\nu_{\infty}\right)$. Then, recalling \ref{ass:gradphiconsistent} and \ref{ass:Tconsistent} and passing to the limit in \eqref{threepointguy} we get
\begin{align*}
r_{w^1}\para{\omega}-r_{w^2}\para{\omega}
&=\Lambda_{\infty}^{-1}\Ppa{D(w^1\para{\omega},w^2\para{\omega})-\ip{\nabla \phi(w^1\para{\omega})-\nabla \phi(w^2\para{\omega}), w^1\para{\omega}-w^2\para{\omega}}{}} \\
&\quad\quad\quad\quad\quad\quad\quad\quad\quad\quad\quad\quad\quad\quad\quad\quad\quad\quad\quad\quad\quad\quad\quad+ M(w^2\para{\omega},w^1\para{\omega}) \\
&=\Lambda_{\infty}^{-1}\Ppa{D(w^1\para{\omega},w^2\para{\omega})-D(w^1\para{\omega},w^2\para{\omega})-D(w^2\para{\omega},w^1\para{\omega})} \\
&\quad\quad\quad\quad\quad\quad\quad\quad\quad\quad\quad\quad\quad\quad\quad\quad\quad\quad\quad\quad\quad\quad\quad+ M(w^2\para{\omega},w^1\para{\omega}) \\
&= -\Lambda_{\infty}^{-1}D(w^2\para{\omega},w^1\para{\omega}) + M(w^2\para{\omega},w^1\para{\omega}) .
\end{align*}
Repeating this argument, replacing $w_{k_i}\para{\omega}$ by $w_{k_j}\para{\omega}$ above, we furthermore have
\begin{align*}
r_{w^1}\para{\omega}-r_{w^2}\para{\omega}
&=\Lambda_{\infty}^{-1}D(w^1\para{\omega},w^2\para{\omega}) - M(w^1\para{\omega},w^2\para{\omega}),
\end{align*}
which shows that
\[
\sbrac{\Lambda_{\infty}^{-1}D(w^1\para{\omega},w^2\para{\omega})- M(w^1\para{\omega},w^2\para{\omega})}+\sbrac{\Lambda_{\infty}^{-1}D(w^2\para{\omega},w^1\para{\omega}) -M(w^2\para{\omega},w^1\para{\omega})} = 0 .
\]
By \ref{ass:littled} and \eqref{aux}, we arrive at
\[
L \sbrac{D\left(w^1\para{\omega},w^2\para{\omega}\right) + D\left(w^2\para{\omega},w^1\para{\omega}\right)} + \eps  \sbrac{d\left(w^1\para{\omega},w^2\para{\omega}\right)+d\left(w^2\para{\omega},w^1\para{\omega}\right)} = 0 ,
\]
or equivalently, in view of \ref{ass:littled},
\begin{equation}\label{eq:identwlim}
\begin{aligned}
D\left(w^1\para{\omega},w^2\para{\omega}\right) + D\left(w^2\para{\omega},w^1\para{\omega}\right) &= 0 \qandq \\
\eps \sbrac{d\left(w^1\para{\omega},w^2\para{\omega}\right)+d\left(w^2\para{\omega},w^1\para{\omega}\right)} &= 0.
\end{aligned}
\end{equation}
To complete the proof, it remains to show that $w^1\para{\omega},w^2\para{\omega}$ for all $\omega\in\tilde{\Omega}$ since $\prob\para{\tilde{\Omega}}=1$.

\begin{itemize}
\item[\ref{assum:wkconvergence1}]
Thanks to \propref{wkclstropt}, we have $\mathfrak{W}\sbrac{\para{w_k\para{\omega}}} \subset \sadp = \{\bar{w}(\omega)\}$. 

\item[\ref{assum:wkconvergence2}]
In this case, we have $\sadp \subset \inte \dom(\phi_p) \times \inte \dom(\phi_p)$ thanks to \ref{ass:legendre}. Thus, in view of \propref{wkclstropt}, $w^i(\omega) \in \mathfrak{W}\sbrac{\seq{w_k\para{\omega}}} \subset \inte \dom(\phi_p) \times \inte \dom(\phi_p)$ for $i=1,2$. Moreover, \eqref{eq:identwlim} gives
\[
D\left(w^1\para{\omega},w^2\para{\omega}\right)+D\left(w^2\para{\omega},w^1\para{\omega}\right)=\iprod{\nabla \phi(w^1(\omega)) - \nabla \phi(w^2(\omega))}{w^1(\omega)-w^2(\omega)} = 0 .
\]
Unless $w^1(\omega)=w^2(\omega)$, this is in contradiction with strict monotonicity of $\nabla \phi$ on $\inte \dom(\phi_p) \times \inte \dom(\phi_p)$ since $\phi_p$ and $\phi_d$ are Legendre.

\item[\ref{assum:wkconvergence2}]
Now $\eps > 0$ and \eqref{eq:identwlim} entails
\[
d\left(w^1\para{\omega},w^2\para{\omega}\right) = 0 ,
\]
whence we conclude $w^1\para{\omega}=w^2\para{\omega}$ by assumption on $d$.
\end{itemize}
\end{proof}

\begin{remark}
The assumptions and results in \thmref{wkconvergence} can be restricted in a modular way, e.g., if only the set of primal solutions is a singleton (and not also the set of dual solutions) then we will retain weak convergence of the primal iterates to the solution to the primal problem.
\end{remark}

\subsection{Strong Convergence under Relative Strong Convexity}
In this part we assume that either $f$, $g$, or both are relatively strongly convex (see \defref{genstrong}) with respect to $\psi_p$ with constant $\scmf$, $\scmg$, or $\scmf+\scmg$, respectively, as in \ref{ass:strconvex}. For brevity, we analyze only the primal case but all of the analogous convergence results will hold for the dual case by making the corresponding assumptions on $h^*$, $l^*$, and $\psi_d$, as in \ref{ass:strconvex}. In addition, if the assumptions made here on the primal functions and entropies hold for the corresponding dual functions and entropies, we will have convergence results for the whole primal-dual sequence $\seq{w_k}$. 

Central to our arguments are the concepts of total convexity and sequential consistency which provide an elegant framework relating convergence in terms of the Bregman divergence and convergence in terms of the ambient norm of the space. We will assume that $\psi_p$ is sequentially consistent and totally convex, which we now go on to define. The following definitions come from \cite{butnariu} although earlier notions of total convexity and its modulus exist.
\begin{definition}
Define, for all $x\in\inte\dom\para{\psi_p}$ and $t\in[0,\infty[$,
\newq{
\Theta_{\psi_p}\para{x,t}\eqdef \inf\brac{D_{\psi_p}\para{x',x}:\norm{x-x'}{\X_p}=t}.
}
The function $\Theta$ is called the modulus of total convexity and it is clearly nondecreasing in $t$ (see \cite[Page 18]{butnariu}). We call a function $\psi_p$ totally convex at a point $x\in\inte\dom \para{\psi_p}$ iff $\Theta_{\psi_p}\para{x,t}>0$ for any $t>0$. We say the function $\psi_p$ is totally convex on a subset $X\subseteq \inte\dom\para{\psi_p}$ iff it is totally convex for each $x\in X$.
\end{definition}
Total convexity is a sort of generalization of strict convexity to functions defined on Banach spaces. Indeed, for finite-dimensional spaces, strict convexity and total convexity are equivalent for functions with full domain \cite[Proposition 1.2.6]{butnariu}. Examples of totally convex functions include the Shannon-Boltzmann entropy, the Hellinger entropy, the Fermi-Dirac entropy, the energy/euclidean entropy, and any strongly convex function as well. %We point out that \ref{ass:strconvex} assumes that there is a unique solution $\xs$ to \eqref{PProb} (and similarly for \eqref{DProb} if we have \ref{ass:strconvex} on the dual).
\begin{definition}
A function $\psi_p$ is called sequentially consistent on a subset $X\subseteq\inte\dom\para{\psi_p}$ iff for any bounded subset $V\subseteq X$, for any $t>0$, we have
\newq{
\inf\limits_{x\in V}\Theta_{\psi_p}\para{x,t} > 0.
}
\end{definition}

\begin{lemma}\label{lem:Dpto0}
Let $\para{\xs,\mus}\in\sadp$ be a saddle point. Assume \ref{ass:hyp}, \ref{ass:legendre}-\ref{ass:error}, and \ref{ass:strconvex}\ref{ass:strconvexP}. Then $D_{\psi_p}\para{\xs,x_k}\to 0$ almost surely. Similarly, if \ref{ass:strconvex}\ref{ass:strconvexD} holds, then $D_{\psi_d}\para{\mus,\mu_k}\to 0$ almost surely.
\end{lemma}
\begin{proof}
Under \ref{ass:strconvex}\ref{ass:strconvexP}, evaluating \eqref{eq:mainestimpos2} in \lemref{estimate} at a saddle point $w=\ws\in\sadp$ we have\fekn
%\begin{multline*}
%\frac{1}{\Lambda_{k}} D\para{\ws,w_{k}}-M\para{\ws,w_{k}} + \ip{\wnk, \ws-w_{k+1}}{} -\scs D'\para{w,w_k} - \scns D'\para{w,w_{k+1}} \\
%\geq \frac{1}{\Lambda_{k+1}} D\para{\ws,w_{k+1}} - M\para{\ws,w_{k+1}}
%\end{multline*}
%which we rewrite\fekn as
\begin{multline*}
\frac{1}{\Lambda_{k}}D\para{\ws, w_k}-\frac{1}{\Lambda_{k+1}}D\para{\ws, w_{k+1}} -M\para{\ws,w_{k}} + M\para{\ws, w_{k+1}} + \ip{\wnk, \ws-w_{k+1}}{} \\
\geq \scmg D_{\psi_p}\para{\xs, x_{k+1}} + \scmf D_{\psi_p}\para{\xs,x_k}.
\end{multline*}
We now break the proof into two cases based on whether $\scmg > 0$ or $\scmf > 0$, starting with $\scmf>0$. Taking the expectation conditioned on the filtration, we have\fekn
\begin{multline}\label{scfest}
\scmf D_{\psi_p}\para{\xs,x_k} \leq \frac{1}{\Lambda_{k}} D\para{\ws, w_k}-\frac{1}{\Lambda_{k+1}}\EX{D\para{\ws, w_{k+1}}}{\filts_k} -M\para{\ws,w_{k}}\\
+\EX{M\para{\ws, w_{k+1}}}{\filts_k}+ \EX{\ip{\wnk, \ws-w_{k+1}}{}}{\filts_k}.
\end{multline}
Applying \lemref{AuxPC} to \eqref{scfest} along with the assumption that $\scmf>0$, \ref{ass:unbiased}, and \ref{ass:error} with \lemref{lem:errorsum}, we find that $\seq{D_{\psi_p}\para{\xs, x_{k}}}\in\ell^1_+\para{\Filts}$ and, in particular, $D_{\psi_p}\para{\xs,x_k}\to 0$ almost surely.

Now, assuming $\scmg > 0$ gives\fekn
\newq{
\scmg D_{\psi_p}\para{\xs,x_{k+1}} \leq \frac{1}{\Lambda_k}D\para{\ws,w_k} - \frac{1}{\Lambda_{k+1}}D\para{\ws, w_{k+1}} - M\para{\ws, w_k} + M\para{\ws, w_{k+1}} + \ip{\wnk, \ws-w_{k+1}}{}.
}
Taking the expectation then leads to\fekn
\newq{
\scmg \TEX{D_{\psi_p}\para{\xs,x_{k+1}}} \leq \frac{1}{\Lambda_k}\TEX{D\para{\ws,w_k}} - \frac{1}{\Lambda_{k+1}}\TEX{D\para{\ws, w_{k+1}}} - \TEX{M\para{\ws, w_k}} + \TEX{M\para{\ws, w_{k+1}}}\\
 + \TEX{\ip{\wnk, \ws-w_{k+1}}{}}.
}
Then, by \remref{LemmaAux} with the assumption $\scmg>0$, \ref{ass:error} and \lemref{lem:errorsum}, we have that $\seq{\TEX{D_{\psi_p}\para{\xs,x_{k}}}}\in\ell^1_+$ and so, by \lemref{AuxAS}, we have that $D_{\psi_p}\para{\xs,x_{k}}\to 0$ almost surely.
\end{proof}

\begin{theorem}\label{thm:strconvergence}
Assume \ref{ass:hyp}, \ref{ass:legendre}-\ref{ass:error}, and \ref{ass:strconvex}\ref{ass:strconvexP} hold, that $\psi_p$ is sequentially consistent on $\U_p$, and assume $\xs$ is the unique solution to the primal problem (i.e., $\solsetp=\brac{\xs}$). Then, if the sublevel sets of $D_{\psi_p}\para{\xs,\cdot}$ are bounded, the sequence $\seq{x_k}$ converges strongly to the solution $\xs$ almost surely. Furthermore, if \ref{ass:strconvex}\ref{ass:strconvexD} holds, $\mus$ is the unique solution to the dual problem, $\psi_d$ is sequentially consistent on $\U_d$, and the sublevel sets of $D_{\psi_d}\para{\mus, \cdot}$ are bounded, then almost surely, the sequence $\seq{w_k}$ converges strongly to the saddle point $\ws$.
\end{theorem}
\begin{proof}
Under these assumptions, \lemref{lem:Dpto0} ensures $D_{\psi_p}\para{\xs,x_k}\to 0$ almost surely. The sublevel sets of $D_{\psi_p}\para{\xs,\cdot}$ are bounded and thus the sequence $\seq{x_k}$ is bounded. Since $\seq{x_k}$ also remains in $\U_p$ by \ref{ass:legendre}, there exists $U_p\subseteq \U_p$ a bounded set such that\fekn $x_k\in U_p$. Since $\psi_p$ is sequentially consistent on $\U_p$, we have, for any $t>0$,
\newq{
\inf\limits_{x\in U_p} \Theta_{\psi_p}\para{x, t}>0.
}
Assume now that $\seq{x_k}$ does not converge strongly to $\xs$. Then there exists a subsequence $\subseq{x_{k_j}}$ and $\epsilon > 0$ such that for all $j\in\N$ it holds,
\newq{
\norm{x_{k_j} - \xs}{\X_p} > \epsilon.
}
Since $\subseq{x_{k_j}}$ is a subsequence of $\seq{x_k}$, $\subseq{D_{\psi_p}\para{\xs, x_{k_j}}}$ is a subsequence of $\seq{D_{\psi_p}\para{\xs, x_k}}$ and so its limit is $0$. Since $\psi_p$ is sequentially consistent and $\norm{x_{k_j}-\xs}{}> \epsilon$, the following is true: for any $j\in\N$,
\nnewq{\label{strconvest}
D_{\psi_p}\para{\xs,x_{k_j}} \geq \Theta_{\psi_p}\para{x_{k_j}, \norm{x_{k_j}-\xs}{\X_p}} \geq \Theta_{\psi_p}\para{x_{k_j},\epsilon}\geq \inf\limits_{x\in \U_p}\Theta_{\psi_p}\para{x,\epsilon}>0,
}
which contradicts the fact that $\lim\limits_{j\to\infty}D_{\psi_p}\para{\xs,x_{k_j}}=0$ since the positive lower bound $\inf\limits_{x\in \U_p}\Theta_{\psi_p}\para{x,\epsilon}$ does not depend on $j$. Thus such a subsequence $\subseq{x_{k_j}}$ cannot exist and the desired claim follows.

Repeating this argument for the dual gives convergence of $\seq{\mu_k}$ to the solution of the dual problem $\mus$ and thus, if \ref{ass:strconvex} holds for the primal and the dual, we have that $\seq{w_k}$ converges to the saddle point $\ws$.
\end{proof}

\begin{remark}
The assumption that the sublevel sets of the the Bregman divergence be bounded, used in \thmref{thm:strconvergence}, holds for a wide class of entropies which includes the Shannon-Boltzmann entropy, the Hellinger entropy, the Fermi-Dirac entropy, the fractional power entropy, and energy/euclidean entropy (see \cite[Remark 4]{BauschkeBolteTeboulle}).
\end{remark}

\begin{remark}
In the statement of  \thmref{thm:strconvergence}, uniqueness of the solution $\xs$ is assumed only for clarity of presentation. Indeed, without the assumption the same argument used in the proof works for every solution $\xs$; and this implies that the solution to the primal problem under our considerations must be unique, as the sequence $x_k$ converges to any solution taken. We do not have a more direct proof for uniqueness of the solution in the general setting of \thmref{thm:strconvergence}, but we point at \propref{uniqueness} where we show a direct proof of uniqueness under the assumption that there exists a saddle point $\ws=\para{\xs,\mus}$ with $\xs \in \U_p$.
\end{remark}

%-----------------------
% Footer to allow Gummi
% to maintain previews
% for individual pieces
%-----------------------
\ifdefined\COMPLETE
\else
\end{document}
\fi

%----------------------------------
% These commands allow Gummi to
% edit each piece individually
% and still maintain previews.
%----------------------------------
\ifdefined\COMPLETE
\else
\documentclass[12pt]{article}
\input{tex_package_header} %file containing all the used libraries
\begin{document}
\fi
%----------------------------------
% Remove at your own risk!
% There is also a footer that ends
% the document if the main wasn't
% loaded.
%----------------------------------
\section{Applications and Numerical Experiments}\label{sec:app}
We examine two applications that satisfy our assumptions for \thmref{thm:ergodic}. The following results will be useful throughout the applications section, particularly when it comes to satisfying \ref{ass:littled}. In the rest of the section, $\norm{\cdot}{q}$, $q \in [1,+\infty]$, will stand for the $\ell^q$ norm on $\R^n$. $\calB^q_r$ is the $\ell^q$ ball of radius $r > 0$.

We begin with a famous result, Pinsker's inequality, which shows that the Kullback-Leibler divergence is strongly convex on the simplex wrt the $\ell^1$ norm.
\begin{lemma}[Pinsker's Inequality \cite{pinsker}]\label{lem:pinsker}
Let $x,y\in\Sigma^n\eqdef \brac{u\in\R^n:u\geq 0, u^T\mathds{1}=1}$ and let $K$ be the Shannon-Boltzmann entropy: $K\para{x} = \sum\limits_{i=1}^nx_i\log\para{x_i}$ on $\R_+^n$ with the convention that $0\log 0=0$. Then it holds
\newq{
\frac{1}{2}\norm{x-y}{1}^2 \leq D_K\para{x,y}.
}
\end{lemma}

\begin{lemma}\label{lipsimpstep}
Let $T: \R^n \to \R^m$, $\C_p=\R_+^n$, $\C_d=\R^n$, $g\para{x}=\iota_{\brac{1}}\para{x^T\mathds{1}}$ and $l \in \Gamma_0(\R^m)$. Choose $\phi_p\para{x} = \sum\limits_{i=1}^n x_i\log\para{x_i}$ on $\R_+^n$ with $0\log 0=0$, and $\phi_d\para{\mu} = \frac{1}{2}\norm{\mu}{2}^2$. Let $\gamma>0$. Then \ref{ass:littled} and \ref{ass:dconsistent} are satisfied with $\tilde{\U}_p = \Sigma^n$, $\U_p=\ri~\Sigma^n$, $\U_d=\tilde{\U}_d=\dom(\partial l^*)$,
\newq{
\lambda_\infty < \frac{1}{L_p + \gamma\norm{T}{2}^2}\qandq \nu_\infty < \frac{1}{L_d + \gamma^{-1}} ,
}
$\eps = \frac{1}{2}$ and
\newq{
d\para{w_1,w_2} = \para{\frac{1}{\lamsup}-L_p-\gamma\norm{T}{2}^2}\norm{x_1-x_2}{1}^2 + \para{\frac{1}{\nusup}-L_d-\frac{1}{\gamma}}\norm{\mu_1-\mu_2}{2}^2.
}
\end{lemma}
In the above, $\ri$ denotes the relative interior.
\begin{proof}
The expressions of $\tilde{\U}_p$, $\U_p$, $\tilde{\U}_d$ and $\U_d$ are immediate. By definition (see \eqref{notation}), for any $w\in\tilde{\U}_p\times\tilde{\U}_d$ and $w'\in\U_p\times\U_d$, we have
\begin{multline*}
\para{\frac{1}{\Lambda_\infty} - L} D\para{w,w'} - M\para{w,w'} = \para{\frac{1}{\lambda_\infty}-L_p}D_p\para{x,x'} + \para{\frac{1}{\nu_\infty} - L_d}\frac{1}{2}\norm{\mu-\mu'}{2}^2 \\
- \ip{T \para{x-x'},\mu-\mu'}{}.
\end{multline*}
Using \lemref{lem:pinsker}, it holds for any $x\in \tilde{\U}_p$ and $x'\in\U_p$,
\newq{
D_p\para{x,x'} \geq \frac{1}{2}\norm{x-x'}{1}^2.
}
By Young's inequality, for any $\gamma>0$, we also have
\begin{align*}
-\ip{T\para{x-x'},\mu-\mu'}{}\geq -\frac{\gamma}{2}\norm{T\para{x-x'}}{2}^2 - \frac{1}{2\gamma}\norm{\mu-\mu'}{2}^2 &\geq -\frac{\gamma}{2}\norm{T}{2}^2\norm{\para{x-x'}}{2}^2 - \frac{1}{2\gamma}\norm{\mu-\mu'}{2}^2 \\
&\geq -\frac{\gamma}{2}\norm{T}{2}^2\norm{\para{x-x'}}{1}^2 - \frac{1}{2\gamma}\norm{\mu-\mu'}{2}^2 .
\end{align*}
Combining the two, we find, for any $w\in\tilde{\U}_p\times\tilde{\U}_d$ and $w'\in\U_p\times\U_d$
\begin{multline*}
\para{\frac{1}{\Lambda_\infty} - L} D\para{w,w'} - M\para{w,w'} 	
%&\geq \frac{1}{2}\sbrac{\para{\frac{1}{\lambda_\infty}-L_p-\gamma\norm{T}{2}^2}\norm{x-x'}{2}^2 + %\para{\frac{1}{\nu_\infty}-L_d-\frac{1}{\gamma}} \norm{\mu-\mu'}{2}^2}\\
\geq \frac{1}{2}\Bigg[\para{\frac{1}{\lambda_\infty}-L_p-\gamma\norm{T}{2}^2}\norm{x-x'}{1}^2 \\
+ \para{\frac{1}{\nu_\infty}-L_d-\frac{1}{\gamma}} \norm{\mu-\mu'}{2}^2\Bigg]
\end{multline*}
which gives \eqref{hyp2}. Checking \ref{ass:dconsistent} is immediate.
%where $\norm{T}{2}^2$ is the square of the classical operator norm and the desired claim follows.
\end{proof}

\subsection{Linear Inverse Problems on the Simplex}

In \cite{chambolle2016ergodic}, the problem of least squares regression on the simplex was considered as an application of the Chambolle-Pock algorithm. A natural extension for \algref{alg:ibpds} is to replace the euclidean norm with the Kullback-Leibler divergence. The Kullback-Leibler divergence is not Lipschitz-smooth and so the Chambolle-Pock algorithm of \cite{chambolle2011first} and \cite{chambolle2016ergodic} cannot be applied, although \cite{chambolle2016ergodic} does allow one to use an entropy in computing the $D$-proximal mapping associated to $g$. 

Consider the problem,
\nnewq{\label{lipsimp}
\min\limits_{\sst{x\in \R^n_+\\ x^T\mathds{1} = 1}} D_K\para{Ax,b} + \beta \norm{B x}{1}
}
where $A \in \R^{m\times n}_{+}$ is a matrix which does not contain any rows which are identically $0$, $b\in\R^m_{++}$, $K$ is the Shannon-Boltzmann entropy with the convention that $0\log 0=0$,
\newq{
K\para{x} = \sum\limits_{i=1}^n x_i\log\para{x_i}, \qwithq \dom(K)=\R_+^n,
}
and $B:\R^n\to \R^{n-1}$ is the linear operator given by
\newq{
B x = \begin{pmatrix}x_2-x_1\\ \vdots \\ x_n-x_{n-1} \end{pmatrix}.
}
It is known that the term $\norm{B x}{1}$ in \eqref{lipsimp} is intended to promote piecewise-constant solutions \cite{rudin1992nonlinear}. Rewriting \eqref{lipsimp}, the associated saddle point problem is given by,
\newq{
\min\limits_{x\in\R^n_+}\max\limits_{\mu\in\R^{n-1}} D_K\para{Ax,b} + \iota_{\brac{1}}\para{x^T\mathds{1}}+\ip{B x,\mu}{} - \iota_{\calB^\infty_\beta}\para{\mu}.
}
Problem~\ref{lipsimp} can be put in a form solvable using the primal-dual algorithm of \cite{chambolle2011first}. But, in addition to working over higher dimensional spaces, this algorithm does not exploit the geometry underlying the problem hence requiring computing (euclidean) prox mappings which are computationally more demanding. 

We can apply \algref{alg:ibpds} with the following choices,
\newq{
f\para{x} = D_K\para{Ax,b},\quad &g\para{x} = \iota_{\brac{1}}\para{x^T \mathds{1}},\quad T=B,\quad h^*\equiv 0,\quad l^*\para{\mu} = \iota_{\calB^\infty_\beta}\para{\mu},\\
&\C_p = \R^n_+,\qandq \C_d = \R^{n-1}.
}
We choose $\phi_p$ and $\phi_d$ (with the same convention $0\log 0=0$) to be
\newq{
\phi_p\para{x} = \sum\limits_{i=1}^n x_i\log\para{x_i}\qandq
\phi_d\para{\mu} = \frac{1}{2}\norm{\mu}{2}^2
}
which induces the divergences $D_p$ and $D_d$
\newq{
D_p\para{x,x'} = \sum\limits_{i=1}^n x_i\log\para{\frac{x_i}{x'_i}}-x_i+x'_i\qandq
D_d\para{\mu,\mu'} = \frac{1}{2}\norm{\mu-\mu'}{2}^2.
}
This gives us the following $D$-prox operator for our problem,
\newq{
\prox_{\lambda_k g}^{D_p}\para{x}\eqdef \argmin\limits_{u\in\C_p}\brac{\lambda_k g\para{u} + D_p\para{u,x}} = \argmin\limits_{\sst{u\in\R^n_+\\ u^T\mathds{1}=1}}\brac{D_p\para{u,x}}= \para{\frac{\exp\para{x_i}}{\sum\limits_{j=1}^n\exp\para{x_j}}}_{i=1}^n.
}
The main hypothesis \ref{ass:hyp} is clearly satisfied in this problem. In order to satisfy \ref{ass:legendre}, we must find a constant $L_p>0$ such that $L_p\phi_p\para{x}-f\para{x} $ is convex for all $x\in \inte\para{\dom\phi_p} = \R^n_{++}$. This is precisely what is shown in \cite[Lemma~8]{BauschkeBolteTeboulle}, which we include here for clarity.
\begin{lemma}
Let $\phi_p\para{x} = \sum\limits_{i=1}^n x_i\log\para{x_i}$, $f\para{x} = D_K\para{Ax,b}$, and $A\in\R^{m\times n}_{+}$ such that none of the rows of $A$ are completely $0$. Then, for any $L_p$ such that
\newq{
L_p\geq\max\limits_{1\leq j\leq m}\para{\sum\limits_{i=1}^nA_{i,j}},
}
$L_p\phi_p - f$ is convex on $\R^n_{++}$.
\end{lemma}
\begin{proof}
See \cite[Lemma~8]{BauschkeBolteTeboulle}
\end{proof}
It remains to choose step sizes $\seq{\lambda_k}$ and $\seq{\nu_k}$ such that \ref{ass:steps} and \ref{ass:littled} are satisfied, for which we refer to \lemref{lipsimpstep}.
\begin{remark}
Notice that the constant $\gamma>0$ in \lemref{lipsimpstep} is arbitrary. For the experiments, we took $\gamma = \norm{B}{2}^{-1}$ to have symmetric step sizes,
\newq{
\lambda_k = \frac{1}{L_p + \norm{B}{2}}\qandq \nu_k = \frac{1}{L_d + \norm{B}{2}}
}
since $L_d=0$ in this problem.
\end{remark}
\begin{figure}[h]
\centering
\includegraphics[width=0.49\linewidth]{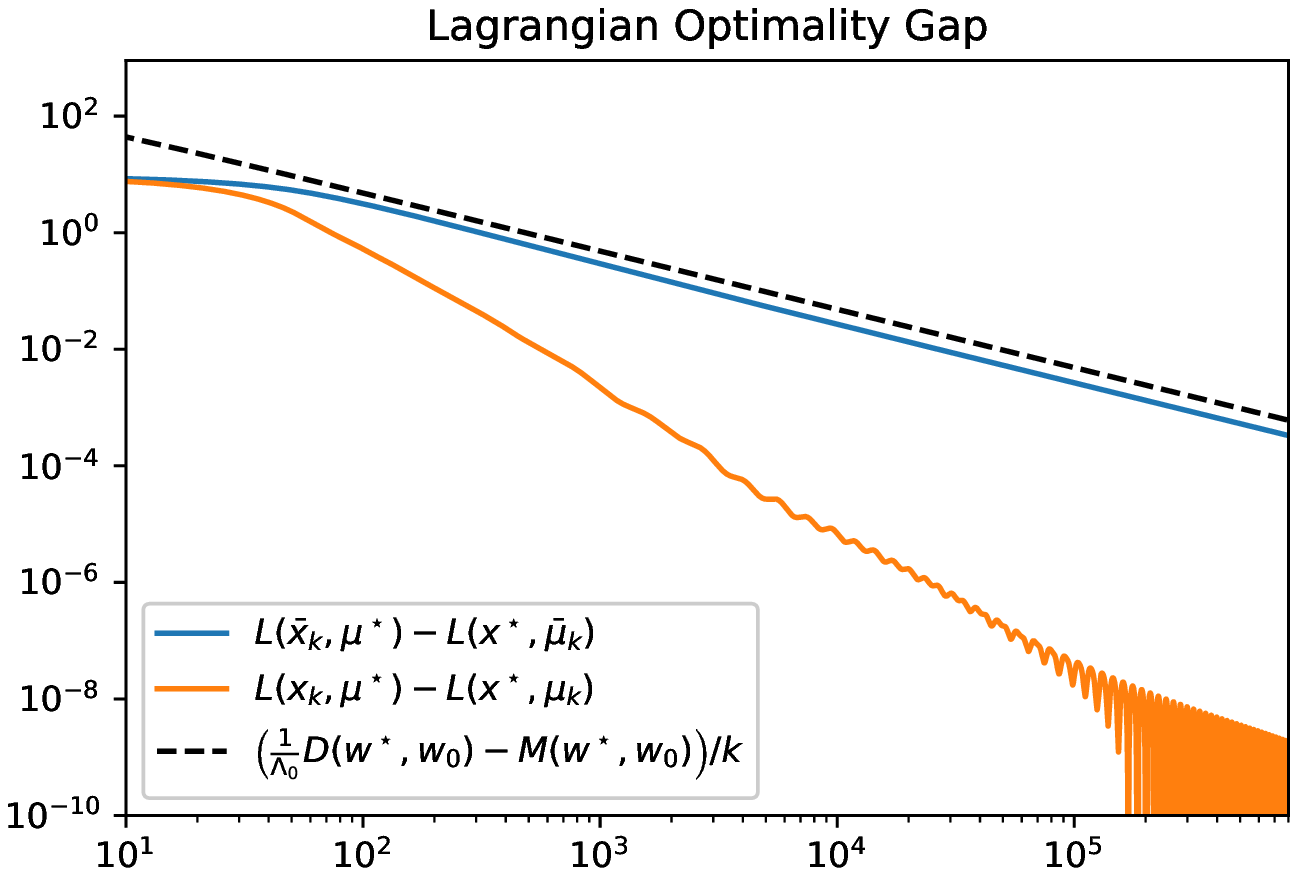}\includegraphics[width=0.485\linewidth]{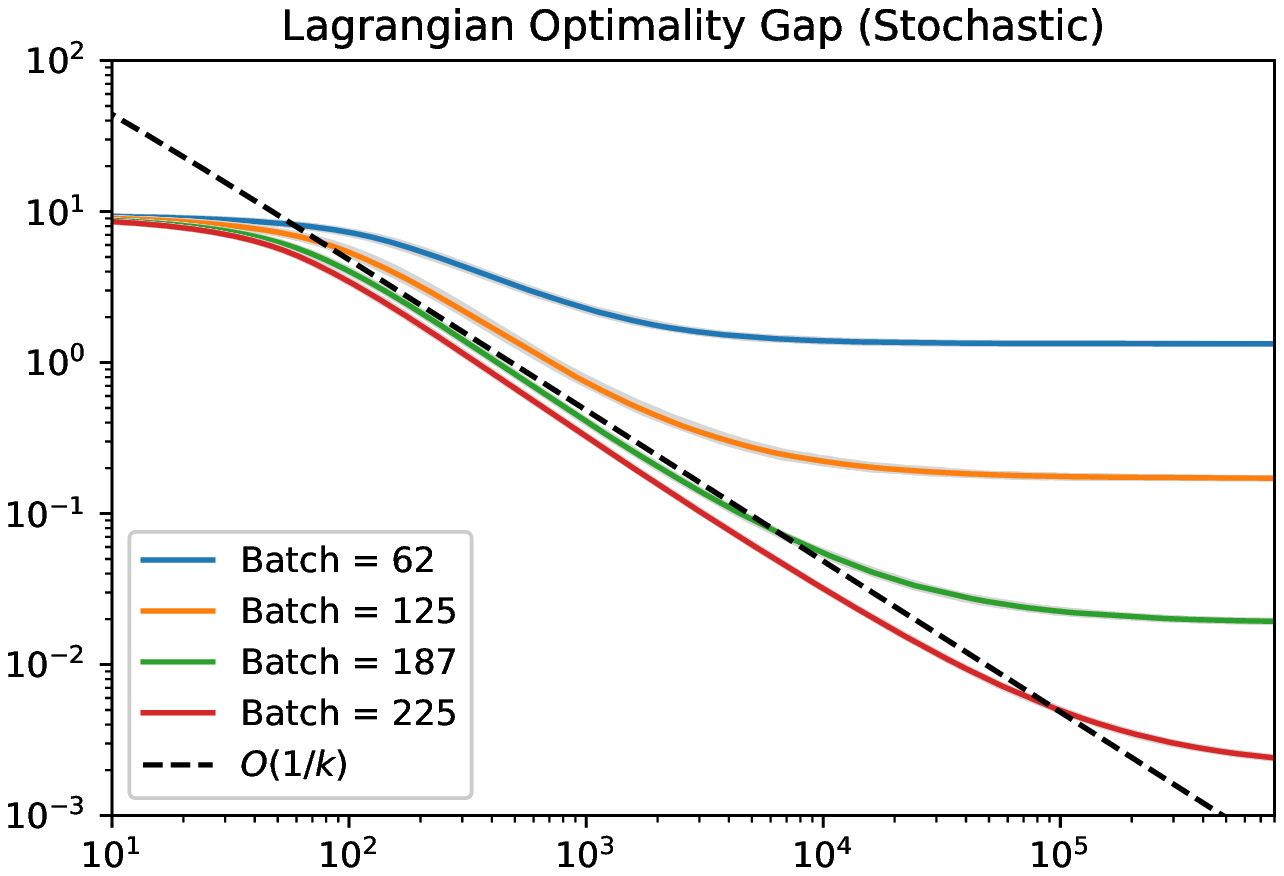}\\
\caption{(Left) Lagrangian optimality gap, computed for both the pointwise iterates in orange and the ergodic iterates in blue, for \algref{alg:ibpds} applied deterministically to the linear inverse problem on the simplex in dimension $n=250$. (Right) The average Lagrangian optimality gap for the ergodic iterates of \algref{alg:ibpds} applied stochastically for various batch sizes, showing convergence to a noise dominated region as predicted by \thmref{thm:ergodic}. The colored lines are the average Lagrangian values for 20 runs of the algorithm, with individual runs displayed in light gray. The $O(1/k)$ theoretical rate is the same in both plots, as given in \thmref{thm:ergodic}.}
\label{fig:lipsimpfig}
\end{figure}

We now apply \algref{alg:ibpds} to solve \eqref{lipsimp} using the step size and entropy choices discussed above. We take $n=250$ and $m=250$, generate $A$ with $a_{i,j}\in[0.01,1.01]$ uniformly i.i.d., and generate $b$ with entries uniformly i.i.d. in $]0,1]$. We initialize with $x_0 = \para{\frac{1}{n},\ldots,\frac{1}{n}}$ and $\mu_0=0$ with the constant step sizes $\lambda_k = \frac{1}{L_p + \norm{B}{2}}$ and $\nu_k = \frac{1}{\norm{B}{2}}$. We also consider $D_K\para{Ax,b} = \sum\limits_{i=1}^m \para{Ax}_i\log\para{\frac{\para{Ax}_i}{b_i}}$ as a finite-sum for which we can sample $f_i\para{x} = \para{Ax}_i\log\para{\frac{\para{Ax}_i}{b_i}}$ in batches, uniformly, when computing the gradient. \thmref{thm:ergodic} ensures convergence of the Lagrangian optimality gap in the deterministic setting for the ergodic iterates, and convergence in expectation to a noise-dominated region for stochastic sampling if the error is bounded in expectation as discussed in \remref{rem:boundederror}. In Lemma~\ref{lem:bounderrbatch} in the appendix, we prove that this is indeed the case.

On the left of \figref{fig:lipsimpfig}, the Lagrangian optimality gap is presented for both the ergodic and pointwise iterates in the deterministic case. We show the same gaps for the stochastic version of the algorithm with batch sampling in \figref{fig:lipsimpfig} on the right. To plot these gaps, we first run the deterministic version of the algorithm for a high number of iterations to find an (approximate) saddle point $\para{\xs,\mus}\in\sadp$ and then rerun the algorithm for 80\% of the number of inital iterations, computing the gap at each iteration. For the stochastic version, we run the algorithm 20 times for each batch size and then plot the average of the gap for the ergodic iterates over these 20 runs in color, with individual runs represented in light gray. Clearly as the batch size increases, the radius of the noise-dominated region shrinks.

\subsection{Variational problems with the entropic Wasserstein distance}
Consider the optimal transport problem between two discrete measures, $\rho$ and $\theta$, defined on two metric spaces $\calX$ and $\calY$. Let $C \in \R^{n \times m}$ be the ground cost on $\calX \times \calY$. The cost $C$ is typically application-dependent, and reflects some prior knowledge on the data to be processed. We regularize the optimal transport problem by subtracting in the objective the entropy of the transport plan $\pi$,
\newq{
E\para{\pi} = -\sum\limits_{i=1}^n\sum\limits_{j=1}^m\pi_{i,j}\log\para{\pi_{i,j}}.
}
The idea of regularizing the optimal transport problem by including the entropy of the transport plan $\pi$ is not new. It was popularized by \cite{cuturi2013sinkhorn} and then explored, for example, in \cite{cuturi2014fast} for computing entropic Wasserstein barycenters,  in \cite{peyre2015entropic} for approximating entropic Wasserstein gradient flows, in \cite{cuturi2016smoothed} for variational Wasserstein problems, in \cite{Cuturi_2018}, etc. For $\gamma > 0$, the entropic regularization of the Kantorovich formulation of optimal transport can be  written as the convex optimization problem
\begin{equation}\label{eq:OTW}
W_{\gamma}\para{\rho,\theta} \eqdef \inf_{\pi \in \Pi\para{\rho,\theta}} \brac{\ip{C,\pi}{} + \gamma\sum\limits_{i=1}^n\sum_{j=1}^m\pi_{i,j}\log\para{\pi_{i,j}} = \gamma\sum\limits_{i=1}^n\sum\limits_{j=1}^m\pi_{i,j}\log\para{\frac{\pi_{i,j}}{\xi_{i,j}}}} ,
\end{equation}
where $\Pi\para{\rho,\theta}\eqdef\brac{\pi\in\R^{n\times m}_+: \pi \mathds{1} = \rho, \pi^T\mathds{1} = \theta}$ is the so-called transportation polytope and $\xi_{i,j} \eqdef \exp\pa{\frac{-C_{i,j}}{\gamma}}$ is the Gibbs Kernel. When $\calX = \calY$, $\gamma = 0$ and $C = d^p$, where $d$ is a distance on $\calX$, then $W_0^{1/p}$ is the well-known $p$-Wasserstein distance. 

We consider solving the following variational problem over discrete measures, i.e., vectors in the simplex $\Sigma^n\eqdef\brac{x: x\geq 0, x^T\mathds{1}=1}$,
\begin{equation}\label{eq:OTvarP}
\min_{\rho \in \Sigma^n} W_{\gamma}(F \rho,\theta) + J \circ B(\rho),
\end{equation}
where $J \in \Gamma_0(\R^p)$, $F: \Sigma^n \to \Sigma^m$ and $B: \R^n \to \R^p$ are both linear operators. Seen as a matrix, $F$ is typically column-stochastic while $\rho\in\Sigma^n$ is a discrete measure over the metric space $\calX$ and $\theta\in\Sigma^m$ is the fixed observed discrete measure over the metric space $\calY$.

Problem~\eqref{eq:OTvarP} is a natural way to solve inverse problems on discrete measures where one assumes that
\newq{
\theta \approx F\rho_0,
}
where $\rho_0$ is an unknown discrete measure over $\calY$ to recover from the observed $\theta$. When $F = \Id$ and $\gamma=0$, \eqref{eq:OTvarP} is closely related to computing the Wasserstein gradient flow (aka JKO flow \cite{JKOpaper}) of $J \circ B$. The JKO flow was first studied in \cite{JKOpaper} as it relates to the Fokker-Planck equation before being generalized (cf. \cite{ambrosiobook}, \cite{santambrogio2016}). Entropic regularization, i.e., with $\gamma >0$, was studied in \cite{peyre2015entropic} to compute Wasserstein gradient flows over spaces of probability distributions with the topology induced by the Wasserstein metric.

Applying Fenchel-Rockafellar duality to \eqref{eq:OTW} (see~\cite[Proposition~2.4]{PC19} for the unregularized case and \cite[Section~5.1]{cuturi2014fast} for the entropic case), it is straightforward to see that problem~\eqref{eq:OTvarP} reads also
\nnewq{\label{eq:OTvarPdual}
\min_{\rho \in \Sigma^n}\sup_{\tau \in \R^m, \eta \in \R^m} \ip{\tau,F\rho}{} + \ip{\eta,\theta}{}  - \gamma \sum_{j=1}^m\sum\limits_{i=1}^m \exp\Ppa{\frac{\tau_i + \eta_j - C_{i,j}}{\gamma}} + J \circ B(\rho) .
}
Taking the supremum over $\eta$, one can easily show that (see also \cite[Proposition~2.1]{Genevay16}),
\nnewq{\label{eq:OTvarPsemidual}
\min_{\rho \in \Sigma^n}\sup_{\tau \in \R^m} \ip{\tau,F \rho}{} - \gamma \sum_{j=1}^m \theta_j \log\Ppa{\sum\limits_{i=1}^m\exp\Ppa{\frac{\tau_i - C_{i,j}}{\gamma}}} + J \circ B(\rho) .
}

\begin{remark}
Observe in \eqref{eq:OTvarPsemidual} that the smooth term in $\tau$ (excluding the inner product $\ip{\tau, F\rho}{}$) is actually a log-sum-exp smooth approximation of the $\max$ function, which would appear naturally when marginalizing with respect to $\eta$ in the case $\gamma=0$.
\end{remark}

Now, dualizing on $J$, we finally get that \eqref{eq:OTvarP} is equivalent to
\begin{equation}\label{eq:OTvarPprimaldual}
\min_{\rho \in \R_+^n}\sup_{\tau \in \R^m, \zeta \in \R^p} \iota_{\brac{1}}(\rho^T\mathds{1}) + \ip{\para{\tau,\zeta},\para{F \rho,B \rho}}{} - \gamma \sum_{j=1}^m \theta_j \log\Ppa{\sum_{i=1}^m\exp\Ppa{\frac{\tau_i - C_{i,j}}{\gamma}}} - J^*(\zeta) .
\end{equation}
The problem in \eqref{eq:OTvarPprimaldual} is a saddle point problem which can be solved with \algref{alg:ibpds} by taking
\newq{
\calC_p = \R_+^n,\quad\calC_d = \R^{m+p},\quad T\rho=\para{F\rho,B\rho},\quad f\para{\rho} = 0,\quad g\para{\rho} = \iota_{\brac{1}}\para{\rho^T\mathds{1}},\\
l^*\para{\mu}=l^*\para{\zeta} = J^*\para{\zeta}, \qandq h^*\para{\mu}=h^*\para{\tau} = \gamma \sum_{j=1}^m \theta_j \log\para{\sum_{i=1}^m\exp\Ppa{\frac{\tau_i - C_{i,j}}{\gamma}}} .
}
The natural choice for the entropies is, again,
\newq{
\phi_p\para{x} = \sum\limits_{i=1}^nx_i\log\para{x_i}\qandq \phi_d\para{\mu} = \frac{1}{2}\norm{\mu}{2}^2.
}
\begin{lemma}
The function $h^*\para{\mu}$ is $L_d$ Lipschitz-smooth for $L_d \geq \gamma^{-1}\sum\limits_{j=1}^m\theta_j=\gamma^{-1}$.
\end{lemma}
\begin{proof}
The log-sum-exp function (with temperature constant $\gamma$),
\newq{
\lse_\gamma\para{x} \eqdef \gamma\log\para{\sum\limits_{i=1}^n\exp\para{\frac{x_i}{\gamma}}},
}
is $C^2$ and convex on $\R^n$ (see \cite[Lemma~4]{gao2017}, \cite[Example~2.16, page~48]{rockafellar1998variational}) and thus so is $h^*\para{\tau,\zeta}$. The gradient, $\nabla \lse_\gamma\para{x}$, is given, component-wise, for each $k\in\brac{1,\ldots,n}$ by
\newq{
\comp{\para{\sigma_\gamma\para{x}}}{k}= \frac{\exp\para{x_k/\gamma}}{\sum\limits_{i=1}^n\exp\para{x_i/\gamma}}.
}
The function $\sigma_\gamma\para{x}$ is called the softmax function with temperature constant $\gamma$ and is Lipschitz-continuous in the euclidean norm with Lipschitz constant $\gamma^{-1}$ (see \cite[Proposition~4]{gao2017}). Thus, to see that the function $h^*$ is Lipschitz-smooth, denote the $j$th column of $C$ as $C_{\cdot,j}$ and notice
\newq{
h^*\para{\mu} = h^*\para{\tau} = \sum\limits_{j=1}^m\theta_j\lse_\gamma\para{\tau - C_{\cdot,j}} \implies \nabla h^*\para{\mu} = \nabla h^*\para{\tau}= \sum\limits_{j=1}^m \theta_j\sigma_\gamma\para{\tau-C_{\cdot,j}}.
}
With this we write,
\newq{
\norm{\nabla h^*\para{\mu} - \nabla h^*\para{\mu'}}{2} &= \norm{\sum\limits_{j=1}^m\theta_j\para{\sigma_\gamma\para{\tau - C_{\cdot,j}} - \sigma_\gamma\para{\tau' - C_{\cdot,j}}}}{2}\\
						&\leq \para{\sum\limits_{j=1}^m\theta_j}\norm{\sigma_\gamma\para{\tau - C_{\cdot,j}}-\sigma_\gamma\para{\tau' - C_{\cdot,j}}}{2}\\
						&\leq \gamma^{-1}\para{\sum\limits_{j=1}^m\theta_j}\norm{\tau - \tau'}{2}
}
and the desired claim follows.
\end{proof}

\begin{figure}[h]
\centering
\includegraphics[width=0.49\linewidth]{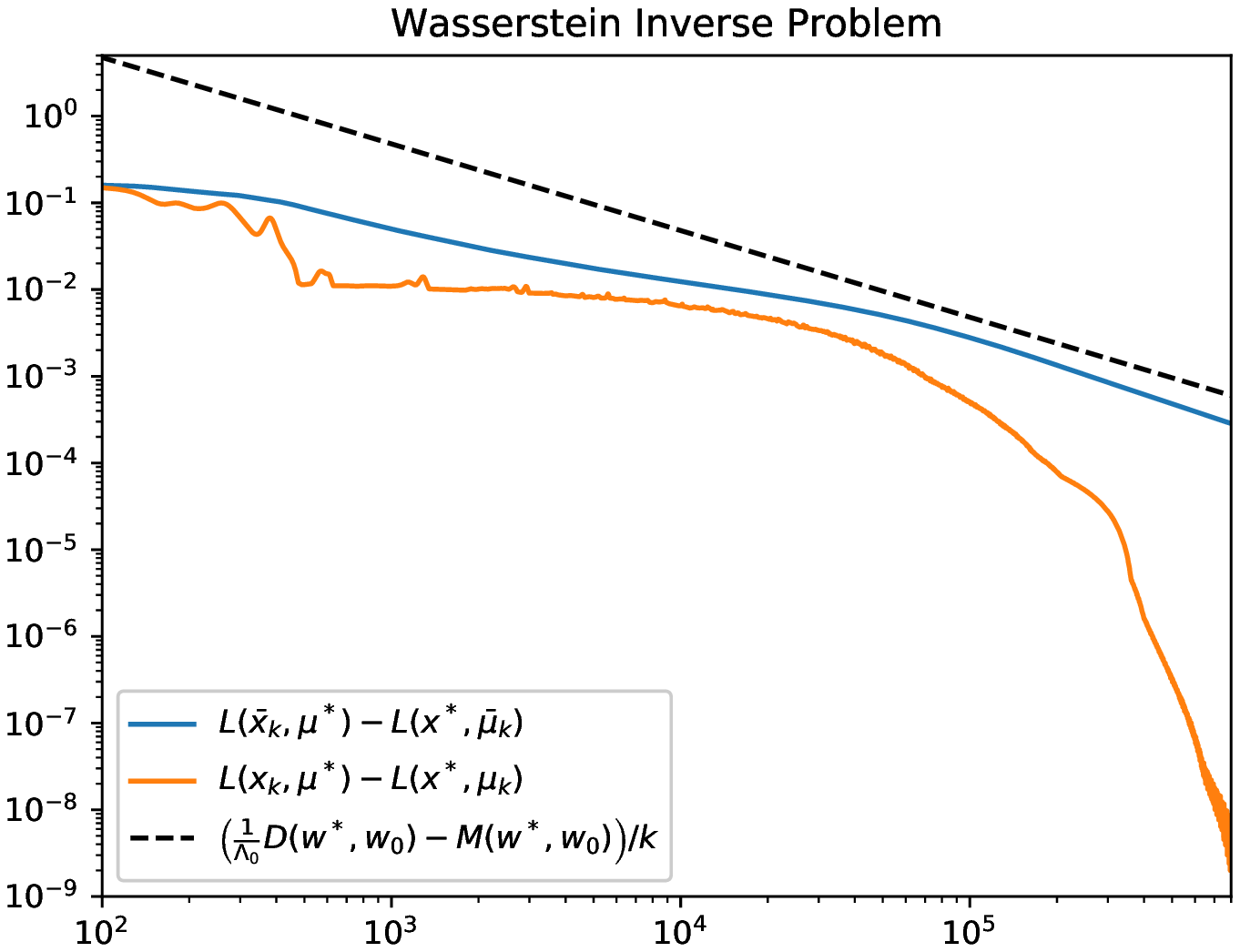}
\includegraphics[width=0.49\linewidth]{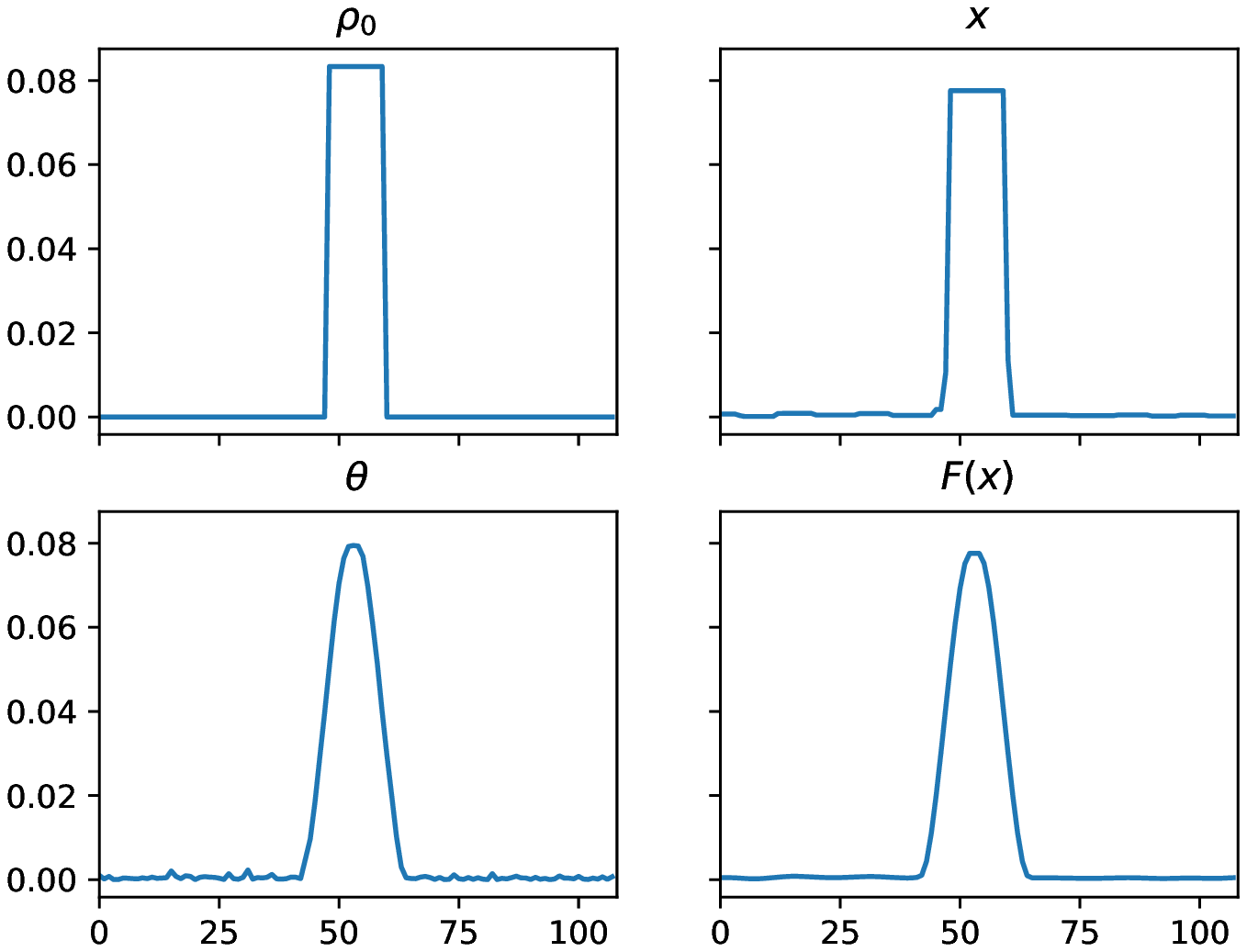}
\caption{(Left) Ergodic and pointwise convergence profiles for \algref{alg:ibpds} applied to the Wasserstein inverse problem with entropic regularization parameter $\gamma=1$, total variation regularization parameter $\beta=1$, and $n=108$. (Right) The ground truth measure $\rho_0$, the recovered measure $x$, the corrupted observation $\theta$, and the image $Fx$ of the recovered measure $x$.}
\label{fig:OT}
\end{figure}

It is clear that \ref{ass:hyp} holds in this setting. It remains to find suitable step sizes $\seq{\lambda_k}$ and $\seq{\nu_k}$ to satisfy \ref{ass:steps} and \ref{ass:littled}. Since the entropies here are exactly the same as in the linear inverse problem on the simplex, we refer again to \lemref{lipsimpstep}. With these step sizes, we consider a one-dimensional instance of the problem with $n=108$, $C_{i,j} = \frac{1}{2}\norm{i-j}{2}^2$, $F$ a convolution operator with kernel $\exp\para{-\frac{1}{1-t^2}}$ for $t\in]-,1,1[$ and $0$ otherwise, $J\circ B$ the total variation \cite{rudin1992nonlinear}, and $\theta \approx F\rho_0$ our observation of $F\rho_0$ is corrupted by Dirichlet distributed noise. We take $x_0=\para{\frac{1}{n},\ldots,\frac{1}{n}}$ and $\mu_0=0$. The results are displayed in \figref{fig:OT}.

%----------------------------
%barycenter stuff
%----------------------------
\begin{remark}
A chief advantage of \eqref{eq:OTvarPprimaldual}, in contrast to optimizing with respect to the transport plan $\pi$, is the significant difference in computational complexity, since the former is operating over $n+m+p$ variables only rather than $nm$. Indeed, one can rewrite problem~\eqref{eq:OTvarP}
as
\newq{
\min\limits_{\sst{\rho\in\Sigma^n\\\pi^T\mathds{1}=\theta\\\pi\mathds{1}=F\rho\\\pi\in\R^{n\times m}_+}}D_K\para{\pi,\xi} + J\circ B\para{\rho} = \min\limits_{\sst{\rho\in\R^n\\ \pi\in\R^{n\times m}}} g\para{\rho,\pi} + f\circ L \para{\rho,\pi}
}
where $g\para{\rho,\pi} \eqdef \iota_{\R^n_+}\para{\rho} + \iota_{\R^{n\times m}_+}\para{\pi} + D_K\para{\pi,\xi}$, $L$ is a linear operator defined as
\newq{
%L\eqdef \pmat{B & 0\\-F & \mathds{1}\\ 0 & ^T\mathds{1}},\quad 
L\para{\rho,\pi} \eqdef \pmat{B\rho\\ \rho^T\mathds{1} \\ -F\rho + \pi\mathds{1}\\\pi^T\mathds{1}}
}
and $f\para{s,t,r,u} \eqdef J\para{s} + \iota_{\brac{1} \times \brac{0} \times \brac{\theta}}\para{t,r,u}$. This formulation is solvable using the Chambolle-Pock algorithm of \cite{chambolle2011first} but, in addition to working with much more variables, over higher dimensional spaces, does not exploit the geometryof the simplex, and requires computing $\prox$ mappings which are computationally more demanding. The $\prox$ operator associated to %$\iota_{\Sigma^n}\para{\rho}$ requires sorting while the prox of 
$D_K\para{\pi,\xi} + \iota_{\R^{n\times m}_+}\para{\pi}$ will require the Lambert $W$ function which is a special function (see \cite{el2013proximal} for more). Even starting from the semidualized form \eqref{eq:OTvarPprimaldual} will require either sorting or incrteasing the number of dual variables if euclidean splitting methods like in \cite{chambolle2011first} are applied.
\end{remark}

\begin{remark}
Although we considered here only a simple Wasserstein inverse problem involving a single observed measure, \algref{alg:ibpds} and our problem framework readily extend to more complicated settings such as computing the Wasserstein barycenter of indirectly observed measures. Wasserstein barycenter problems were first introduced in \cite{agueh2011barycenters} without entropic regularization of the Wasserstein distance. Later, the use of entropic regularization of the Wasserstein distance to speed up computation of barycenters was put forth in \cite{cuturi2014fast}, however the barycenter itself was not regularized; such developments would come later, e.g., \cite{cazelles2017regularized}, \cite{bigot2019data}, etc, and even then the problems considered did not include the possibility of observing the image of the measure $\theta$ under a linear operator $F$ rather than observing the measure $\theta$ itself.

Let $q\in\N$ and consider $q$ reference measures $\theta^i\in \R^{n_i}$ with $n_i\in\N$ for each $1\leq i\leq q$, each having been observed through some linear operator $F^i:\R^{n}\to\R^{n_i}$ applied to an unknown discrete measure $\rho^i\in\Sigma^n$, i.e., $\theta^i\approx F^i\rho^i$. Let $\alpha \in \Sigma^q$. Then we can write the regularized Wasserstein barycenter problem as
\newq{
\min\limits_{\rho\in\Sigma^n} \sum\limits_{k=1}^q \alpha_qW_{\gamma_k}\para{F^k\rho, \theta^k} + \sum\limits_{r=1}^{q'}J^r\circ B^r\para{\rho}
}
which is equivalent to the following,
\begin{multline*}
\min_{\rho\in\Sigma^n}\sup\limits_{\sst{\tau^1\in\R^{n_1},\ldots, \tau^q\in\R^{n_q}\\ \zeta^1\in\R^{m_1},\ldots,\zeta^{q'}\in\R^{m_{q'}}}} \sum\limits_{k=1}^q\sbrac{ \ip{\alpha_k\tau_k, F_k\rho}{} - \alpha_k\gamma_k\sum\limits_{j=1}^{n_k}\theta^k_j\log\para{\sum\limits_{i=1}^{n_k}\exp\para{\frac{\tau^k_{i} - C^k_{i,j}}{\gamma_k}}}} \\
+ \sum\limits_{r=1}^{q'}\sbrac{\ip{\zeta^r, B^r \rho}{} - \para{J^r}^*\para{\zeta^r}}.
\end{multline*}
This formulation of the problem can be solved with  with \algref{alg:ibpds} by taking
\newq{
\C_p = \R_+^n,\quad \C_d &= \R^{n_1}\times\cdots\times\R^{n_q}\times\R^{m_1}\times\cdots\times\R^{m_{q'}},\quad f\para{\rho} = 0,\quad g\para{\rho} = \iota_{\ens{\sum\limits_{i=1}^n\rho_i=1}}\para{\rho},\\
l^*\para{\mu} &= l^*\para{\zeta_1,\ldots,\zeta_{q'}} = \sum\limits_{l=1}^{q'} J_l^*\para{\zeta_l},\qandq\\
h^*\para{\mu} &= h^*\para{\tau^1,\ldots,\tau^q} =  \sum\limits_{k=1}^q\alpha_k\gamma_k \sum\limits_{j=1}^{n_k} \theta^k_j \log\para{\sum\limits_{i=1}^{n_k}\exp\para{\frac{\tau^k_i - C^k_{i,j}}{\gamma_k}}},
}
with the same entropy choices as we took for \eqref{eq:OTvarPprimaldual}.
\end{remark}

\begin{remark}
Consider the same setup as in the previous remark with $\para{\theta^{1},\ldots,\theta^{q}}$ and let $\beta\in\R_{+}$. Another interesting formulation of the regularized Wasserstein barycenter problem that can be solved using \algref{alg:ibpds} is the following
\newq{
\min\limits_{\rho \in \Sigma^n}\min\limits_{\rho_1,\ldots,\rho_q \in \Sigma^n} \sum\limits_{i=1}^q \BPa{W_{\gamma_i}\para{\theta^i, F^i\rho_i} + J\circ A\para{\rho^i}} + \beta\sum\limits_{i=1}^q \alpha_iW_{\gamma_i}\para{\rho_i,\rho}.
}
This problem is simultaneously solving the Wasserstein inverse problem for each observed measure $\theta^i$ while also finding a barycenter $\rho$ among the proposed solutions $\rho^i$ of the Wasserstein inverse problems.
\end{remark}
%-----------------------
% Footer to allow Gummi
% to maintain previews
% for individual pieces
%-----------------------
\ifdefined\COMPLETE
\else
\end{document}
\fi

\section*{Acknowledgements}
ASF was supported by the ERC Consolidated grant NORIA and by the Air Force Office of Scientific Research, Air Force Material Command, USAF, under grant number FA9550-19-1-7026. JF was partly supported by Institut Universitaire de France. CM was supported by Project MONOMADS funded by Conseil R\'egional de Normandie.

\appendix
%----------------------------------
% These commands allow Gummi to
% edit each piece individually
% and still maintain previews.
%----------------------------------
\ifdefined\COMPLETE
\else
\documentclass[12pt]{article}
\input{tex_package_header} %file containing all the used libraries
\begin{document}
\fi
%----------------------------------
% Remove at your own risk!
% There is also a footer that ends
% the document if the main wasn't
% loaded.
%----------------------------------
\section{Appendix}\label{sec:appendix}
\begin{lemma}\label{lem:bounderrbatch}
For each $k\in\N$, for a fixed batch size $0<q< m$, let $B_k\subset\brac{1,\ldots,m}$ be the batch of $q$ indices sampled at iteration $k$ and define $B_k^c \eqdef \brac{1,\ldots,m}\setminus B_k$. Consider the error term induced by the batch sampling:
\newq{
\pnk = -\nabla \para{\sum\limits_{i\in B_k^c} \para{Ax}_i\log\para{\frac{\para{Ax}_i}{b_i}}}.
}
If the entries of $A$ are positive then $\seq{\TEX{\ip{\pnk,x-x_{k+1}}{}}}$ is bounded for all $x\in \Sigma^n$.
\end{lemma}
\begin{proof}
By Cauchy-Schwarz we have
\newq{
\TEX{\ip{\pnk, x-x_{k+1}}{}} \leq \TEX{\norm{\pnk}{}\norm{x-x_{k+1}}{}}\leq \diam_{\Sigma^n}\TEX{\norm{\pnk}{}}
}
and so it suffices to bound $\TEX{\norm{\pnk}{}}$ for all $k\in\N$. Rather than bound the expectation itself, we will provide a coarse bound which holds deterministically for an arbitrary batch. For any batch $B\subset \brac{1,\ldots,m}$ of size $0<q< m$, define $\tilde{A}\in\R^{(m-q)\times n}_{++}$ to be the matrix composed of rows which were not sampled in the batch $B$ and similarly for $\tilde{b}\in\R^{m-q}_{++}$,
\newq{
\tilde{A}\eqdef \begin{bmatrix}A_{i,\cdot}\end{bmatrix}_{i\in B^c}\quad\quad\mbox{and}\quad\quad \tilde{b}\eqdef \para{b_i}_{i\in B^c}.
}
Using component wise $\log$ and division we have, for all $k\in\N$,
\newq{
\norm{\nabla \para{\sum\limits_{i\in B^c} \para{Ax_k}_i\log\para{\frac{\para{Ax_k}_i}{b_i}}}}{} = \norm{\tilde{A}^T \log\para{\frac{\tilde{A}x_k}{\tilde{b}}} }{} &\leq \norm{\tilde{A}}{} \norm{\log\para{\frac{\tilde{A}x_k}{\tilde{b}}}}{}\\
&\leq \norm{\tilde{A}}{}\para{\norm{\log\para{\tilde{A}x_k}}{} + \norm{\log\para{\tilde{b}}}{}}
}
As $A$ and $b$ are fixed from the problem data, $\norm{\tilde{A}}{}$ and $\norm{\log\para{\tilde{b}}}{}$ are bounded. All that remains is to bound $\norm{\log\para{\tilde{A}x_k}}{}$, for which we first recall that $x_k\in\Sigma^n\cap \R^n_{++}$ for all $k\in\N$ by design of the algorithm and the choice of $\phi_p$. Let $\underline{a} = \min\limits_{i,j}A_{i,j} > 0$ and $\overline{a} = \max\limits_{i,j}A_{i,j}$, then for each $i\in\brac{1,\ldots,m-q}$
\newq{
\log\para{\ip{\tilde{A}_{i,\cdot}, x_k}{}} \leq \log\para{\norm{\tilde{A}_{i,\cdot}}{\infty}\norm{x_k}{1}}=\log\para{\norm{\tilde{A}_{i,\cdot}}{\infty}}\leq \log\para{\overline{a}}
}
as well as
\newq{
\log\para{\ip{\tilde{A}_{i,\cdot}, x_k}{}} \geq \log\para{\min\limits_{j}\tilde{A}_{i,j}} \geq \log\para{\underline{a}}
}
so that the components of $\log\para{\tilde{A}x}$ are contained in a ball of radius $\max\brac{\absv{\log\para{\underline{a}}}, \absv{\log\para{\overline{a}}}}$, which is finite since the entries of $A$ are positive and $A$ is fixed. Thus $\norm{\log\para{\tilde{A}x}}{}$ is bounded and the proof is complete.
\end{proof}

\begin{proposition}\label{uniqueness}
Assume \ref{ass:hyp}, \ref{ass:legendre}, and \ref{ass:strconvex} hold and that $\psi_p$ is totally convex and sequentially consistent on $\U_p$. Moreover, suppose that $\solsetp \cap \U_p\neq\emptyset$; meaning that there exists at least one solution $\xs$ of \eqref{PProb} in $\U_p$. Then, there exists a unique solution to the primal problem (i.e., $\solsetp=\brac{\xs}$).
\end{proposition}

\begin{proof}
First notice that, from \ref{ass:strconvex}, $\dom\para{\phi_p} \subseteq \dom \para{\psi_p}$ and so $\U_p \subseteq \inte\dom\para{\phi_p}\subseteq \inte \dom \para{\psi_p}$. As $\psi_p$ is sequentially consistent on $\U_p$, we have that, for any bounded subset $V\subseteq \U_p$ and for any $t>0$,
\newq{
\inf\limits_{x\in V} \Theta_{\psi_p}\para{x,t} > 0.
}
Suppose for instance that \ref{ass:strconvex} holds specifically with $f$ relatively strongly convex with respect to $\psi_p$. Then, by \defref{genstrong}, for any $x,y\in \inte\dom\para{\psi_p}$,
\newq{
\scmf D_{\psi_p}\para{x,y} \leq  D_{f} \para{x,y}.
}
Then, for any bounded subset $V\subseteq \U_p$ and for any $t>0$,
\newq{
  % \leq \Theta_{f}\para{x,t}\implies \inf\limits_{x\in V}\scmf \Theta_{\psi_p}\para{x,t} \leq 
  0< \scmf \inf\limits_{x\in V} \Theta_{\psi_p}\para{x,t}\leq \inf\limits_{x\in V}\Theta_{f}\para{x,t},
}
and so $f$ is sequentially consistent on $\U_p$. Recall from \cite[Proposition page 50-51]{butnariu2003uniform} that sequential consistency of $f$ on the set $\U_p$ implies uniform convexity of $f$ on any bounded subset $V\subseteq \U_p$. Denote by $\xs$ a point in $\solsetp\cap\U_p$. As $\xs\in\U_p$, that is an open set, there is a ball $\calB_{\tau}(\xs)$, for $\tau > 0$, which is bounded and contained in $\U_p$. In particular, $f$ is uniformly convex on $\calB_{\tau}(\xs)$ and so is the objective function in \eqref{PProb}. Suppose by contradiction that there exists another solution $\bar{x}\in\solsetp$ with $\bar{x}\neq\xs$. By convexity, the segment connecting $\xs$ and $\bar{x}$ is contained in $\solsetp$. Then, the intersection of this segment with $\calB_{\tau}(\xs)$ has more than one element and is contained both in $\solsetp \cap \calB_{\tau}(\xs)$. This is a contradiction with the uniform convexity of the objective function in \eqref{PProb} on $\calB_{\tau}(\xs)$. % Indeed, as $f$ is uniformly convex on $B$,  which defines $\solsetp$.

% and, as it is contained in $B$, $f$ is uniformly convex on $Q$. But $Q$ is also contained in $\solsetp$. So we found a set $Q$ where $f$ is uniformly convex and has more than one minimizer, that is a contradiction.
\end{proof}
%-----------------------
% Footer to allow Gummi
% to maintain previews
% for individual pieces
%-----------------------
\ifdefined\COMPLETE
\else
\end{document}
\fi

\small
%%%%%%%%%%%%%%%%%%%%%%%%%%%%%%%%%%%
\begin{small}
\bibliographystyle{plain}
\bibliography{ibpd}

\begin{thebibliography}{10}

\bibitem{agueh2011barycenters}
Martial Agueh and Guillaume Carlier.
\newblock Barycenters in the wasserstein space.
\newblock {\em SIAM Journal on Mathematical Analysis}, 43(2):904--924, 2011.

\bibitem{ambrosiobook}
L.~Ambrosio, N.~Gigli, and G.~Savare.
\newblock {\em Gradient Flows}.
\newblock Lectures in Mathematics. ETH Z{\"u}rich. Birkh{\"a}user Basel, 2008.

\bibitem{BauschkeBorweinCombettes01}
Heinz~H. Bauschke, Jonathan~M. Borwein, and Patrick~L. Combettes.
\newblock Essential smoothness, essential strict convexity, and legendre
  functions in banach spaces.
\newblock {\em Communications in Contemporary Mathematics}, 03(04):615--647,
  2001.

\bibitem{BauschkeBorweinCombettes}
Heinz~H. Bauschke, Jonathan~M. Borwein, and Patrick~L. Combettes.
\newblock Bregman monotone optimization algorithms.
\newblock {\em SIAM Journal on Control and Optimization}, 42(2):596--636, 2003.

\bibitem{bigot2019data}
J{\'e}r{\'e}mie Bigot, Elsa Cazelles, and Nicolas Papadakis.
\newblock Data-driven regularization of wasserstein barycenters with an
  application to multivariate density registration.
\newblock {\em Information and Inference: A Journal of the IMA}, 8(4):719--755,
  2019.

\bibitem{birnbaum2011distributed}
Benjamin Birnbaum, Nikhil~R Devanur, and Lin Xiao.
\newblock Distributed algorithms via gradient descent for fisher markets.
\newblock In {\em Proceedings of the 12th ACM conference on Electronic
  commerce}, pages 127--136, 2011.

\bibitem{Browder65}
F.~E. Browder.
\newblock Multi-valued monotone nonlinear mappings and duality mappings in
  banach spaces.
\newblock {\em Trans. Amer. Math. Soc.}, 118:338--351, 1965.

\bibitem{Browder1966}
F.~E. Browder.
\newblock Fixed point theorems for nonlinear semicontractive mappings in banach
  spaces.
\newblock {\em Arch. Rational Mech. Anal.}, 21:259--269, 1966.

\bibitem{Bui2021}
M.~N. B\`ui and P.~L. Combettes.
\newblock Bregman forward-backward operator splitting.
\newblock {\em Set-Valued and Variational Analysis}, 29(3):583--603, 2021.

\bibitem{butnariu2003uniform}
Dan Butnariu, Alfredo Iusem, and Constantin Zalinescu.
\newblock On uniform convexity, total convexity and convergence of the proximal
  point and outer bregman projection algorithms in banach spaces.
\newblock {\em Journal of Convex Analysis}, 10, 01 2003.

\bibitem{butnariu}
Dan Butnariu and Alfredo~N Iusem.
\newblock {\em Totally convex functions for fixed points computation and
  infinite dimensional optimization}, volume~40.
\newblock Springer Science \& Business Media, 2000.

\bibitem{cazelles2017regularized}
Elsa Cazelles, J{\'e}r{\'e}mie Bigot, and Nicolas Papadakis.
\newblock Regularized barycenters in the wasserstein space.
\newblock In {\em International Conference on Geometric Science of
  Information}, pages 83--90. Springer, 2017.

\bibitem{chambolle2011first}
A.~Chambolle and T.~Pock.
\newblock A first-order primal-dual algorithm for convex problems with
  applications to imaging.
\newblock {\em Journal of Mathematical Imaging and Vision}, 40(1):120--145,
  2011.

\bibitem{chambolle2018stochastic}
Antonin Chambolle, Matthias~J. Ehrhardt, Peter Richt{\'a}rik, and
  Carola-Bibiane Sch{\"o}nlieb.
\newblock Stochastic primal-dual hybrid gradient algorithm with arbitrary
  sampling and imaging applications.
\newblock {\em SIAM Journal on Optimization}, 28(4):2783--2808, 2018.

\bibitem{chambolle2016ergodic}
Antonin Chambolle and Thomas Pock.
\newblock On the ergodic convergence rates of a first-order primal--dual
  algorithm.
\newblock {\em Mathematical Programming}, 159(1-2):253--287, 2016.

\bibitem{chen1993convergence}
Gong Chen and Marc Teboulle.
\newblock Convergence analysis of a proximal-like minimization algorithm using
  bregman functions.
\newblock {\em SIAM Journal on Optimization}, 3(3):538--543, 1993.

\bibitem{Cioranescu1990}
I.~Cioranescu.
\newblock {\em Geometry of Banach spaces, DualityMappings and Nonlinear
  Problems}.
\newblock Kluwer, Dordrecht, 1990.

\bibitem{combettes2012}
Patrick~L. Combettes and Jean-Christophe Pesquet.
\newblock Primal-dual splitting algorithm for solving inclusions with mixtures
  of composite, lipschitzian, and parallel-sum type monotone operators.
\newblock {\em Set-Valued and Variational Analysis}, 20(2):307--330, Jun 2012.

\bibitem{Nguyen15}
P.L. Combettes and Q.~Nguyen.
\newblock Solving composite monotone inclusions in reflexive banach spaces by
  constructing best bregman approximations from their kuhn-tucker set.
\newblock {\em Journal of Convex Analysis}, 23:481--510, 05 2016.

\bibitem{condat2012primal}
L.~Condat.
\newblock A primal--dual splitting method for convex optimization involving
  lipschitzian, proximable and linear composite terms.
\newblock {\em Journal of Optimization Theory and Applications}, pages 1--20,
  2012.

\bibitem{cuturi2013sinkhorn}
Marco Cuturi.
\newblock Sinkhorn distances: Lightspeed computation of optimal transport.
\newblock In {\em Advances in neural information processing systems}, pages
  2292--2300, 2013.

\bibitem{cuturi2014fast}
Marco Cuturi and Arnaud Doucet.
\newblock Fast computation of wasserstein barycenters.
\newblock {\em Journal of Machine Learning Research}, 2014.

\bibitem{cuturi2016smoothed}
Marco Cuturi and Gabriel Peyr{\'e}.
\newblock A smoothed dual approach for variational wasserstein problems.
\newblock {\em SIAM Journal on Imaging Sciences}, 9(1):320--343, 2016.

\bibitem{Cuturi_2018}
Marco Cuturi and Gabriel Peyr\'{e}.
\newblock Semidual regularized optimal transport.
\newblock {\em SIAM Review}, 60(4), Jan 2018.

\bibitem{el2013proximal}
Mireille El~Gheche, Jean-Christophe Pesquet, and Joumana Farah.
\newblock A proximal approach for optimization problems involving kullback
  divergences.
\newblock In {\em 2013 IEEE International Conference on Acoustics, Speech and
  Signal Processing}, pages 5984--5988. IEEE, 2013.

\bibitem{gao2017}
Bolin Gao and Lacra Pavel.
\newblock On the properties of the softmax function with application in game
  theory and reinforcement learning, 2017.

\bibitem{Genevay16}
Aude Genevay, Marco Cuturi, Gabriel Peyr\'{e}, and Francis Bach.
\newblock Stochastic optimization for large-scale optimal transport.
\newblock In {\em Advances in Neural Information Processing Systems 29}, pages
  3440--3448. Curran Associates, Inc., 2016.

\bibitem{hamedani2018primal}
Erfan~Yazdandoost Hamedani and Necdet~Serhat Aybat.
\newblock A primal-dual algorithm for general convex-concave saddle point
  problems.
\newblock {\em arXiv preprint arXiv:1803.01401}, 2018.

\bibitem{BauschkeBolteTeboulle}
J.~Bolte H.H.~Bauschke and M.~Teboulle.
\newblock A descent lemma beyond lipschitz gradient continuity: first-order
  methods revisited and applications.
\newblock {\em Math. Oper. Res.}, 42(2):330--348, 2017.

\bibitem{jiang2021bregman}
Xin Jiang and Lieven Vandenberghe.
\newblock Bregman primal--dual first-order method and application to sparse
  semidefinite programming.
\newblock {\em arXiv preprint}, 2021.

\bibitem{JKOpaper}
Richard Jordan, David Kinderlehrer, and Felix Otto.
\newblock The variational formulation of the fokker--planck equation.
\newblock {\em SIAM Journal on Mathematical Analysis}, 29(1):1--17, 1998.

\bibitem{pesquetplaying}
N.~{Komodakis} and J.~{Pesquet}.
\newblock Playing with duality: An overview of recent primal-dual approaches
  for solving large-scale optimization problems.
\newblock {\em IEEE Signal Processing Magazine}, 32(6):31--54, 2015.

\bibitem{lebedev67}
VN~Lebedev and NT~Tynjanski{\i}.
\newblock Duality theory of concave-convex games.
\newblock In {\em Soviet Math. Dokl}, volume~8, pages 752--756, 1967.

\bibitem{lu2019relative}
Haihao Lu.
\newblock "relative continuity" for non-lipschitz nonsmooth convex optimization
  using stochastic (or deterministic) mirror descent.
\newblock {\em INFORMS Journal on Optimization}, 1(4):288--303, 2019.

\bibitem{lu2020generalized}
Haihao Lu and Robert~M. Freund.
\newblock Generalized stochastic frank--wolfe algorithm with stochastic
  ``substitute'' gradient for structured convex optimization.
\newblock {\em Mathematical Programming}, Mar 2020.

\bibitem{lu2018relatively}
Haihao Lu, Robert~M. Freund, and Yurii Nesterov.
\newblock Relatively smooth convex optimization by first-order methods, and
  applications.
\newblock {\em SIAM Journal on Optimization}, 28(1):333--354, 2018.

\bibitem{mclinden74}
L.~McLinden.
\newblock An extension of fenchel's duality theorem to saddle functions and
  dual minimax problems.
\newblock {\em Pacific J. Math.}, 50(1):135--158, 1974.

\bibitem{moreau64}
J.-J. Moreau.
\newblock Th\'eor\`emes ``inf-sup,''.
\newblock {\em C. R. Acad. Sci. Paris S\'er. A Math.}, 258:2720--2722, 1964.

\bibitem{Nguyen17}
Q.~V. Nguyen.
\newblock Forward-backward splitting with bregman distances.
\newblock {\em Vietnam J. Math.}, 45(519--539), 2017.

\bibitem{opial1967weak}
Z.~Opial.
\newblock Weak convergence of the sequence of successive approximations for
  nonexpansive mappings.
\newblock {\em Bulletin of the American Mathematical Society}, 73(4):591--597,
  1967.

\bibitem{peyre2015entropic}
Gabriel Peyr{\'e}.
\newblock Entropic approximation of wasserstein gradient flows.
\newblock {\em SIAM Journal on Imaging Sciences}, 8(4):2323--2351, 2015.

\bibitem{PC19}
Gabriel Peyr{\'e}, Marco Cuturi, et~al.
\newblock Computational optimal transport: With applications to data science.
\newblock {\em Foundations and Trends{\textregistered} in Machine Learning},
  11(5-6):355--607, 2019.

\bibitem{pinsker}
Mark~Semenovich Pinsker.
\newblock The information stability of gaussian random variables and processes
  (in russian)).
\newblock In {\em Doklady Akademii Nauk}, volume 133, pages 28--30. Russian
  Academy of Sciences, 1960.

\bibitem{polyak1987introduction}
B.~T. Polyak.
\newblock {\em Introduction to optimization}.
\newblock Optimization Software, 1987.

\bibitem{rasch2020inexact}
Julian Rasch and Chambolle Antonin.
\newblock Inexact first-order primal--dual algorithms.
\newblock {\em Computational Optimization and Applications}, 76(2):381--430,
  2020.

\bibitem{robbinssiegmund}
H.~Robbins and D.~Siegmund.
\newblock A convergence theorem for non negative almost supermartingales and
  some applications.
\newblock In Jagdish~S. Rustagi, editor, {\em Optimizing Methods in
  Statistics}, pages 233 -- 257. Academic Press, 1971.

\bibitem{rockafellar64}
R.~T. Rockafellar.
\newblock Minimax theorems and conjugate saddle-functions.
\newblock {\em Mathematica Scandinavica}, 14(2):151--173, 1964.

\bibitem{rockafellar1998variational}
R.~T. Rockafellar and R.~Wets.
\newblock {\em Variational analysis}, volume 317.
\newblock Springer Verlag, 1998.

\bibitem{rudin1992nonlinear}
L.~I. Rudin, S.~Osher, and E.~Fatemi.
\newblock Nonlinear total variation based noise removal algorithms.
\newblock {\em Physica D: Nonlinear Phenomena}, 60(1):259--268, 1992.

\bibitem{santambrogio2016}
Filippo Santambrogio.
\newblock $\{$Euclidean, metric, and Wasserstein$\}$ gradient flows: an
  overview.
\newblock {\em Bulletin of Mathematical Sciences}, 7(1):87--154, 2017.

\bibitem{silvetisiopt}
Antonio Silveti-Falls, Cesare Molinari, and Jalal Fadili.
\newblock Generalized conditional gradient with augmented lagrangian for
  composite minimization.
\newblock {\em SIAM Journal on Optimization}, 30(4):2687--2725, 2020.

\bibitem{Simons08}
S.~Simons.
\newblock {\em From Hahn-Banach to Monotonicity}, volume 1693 of {\em Lecture
  Notes in Math.}
\newblock Springer-Verlag, New York, 2008.

\bibitem{vu2011splitting}
B.~C. V{\~u}.
\newblock A splitting algorithm for dual monotone inclusions involving
  cocoercive operators.
\newblock {\em Advances in Computational Mathematics}, pages 1--15, 2011.

\bibitem{Xu14}
Hong Xu, Tae-Hwa Kim, and Ximing Yin.
\newblock Weak continuity of the normalized duality map.
\newblock {\em Journal of Nonlinear and Convex Analysis}, 15(3):595--604, 2014.

\bibitem{Zalinescu02}
C.~Z\u{a}linescu.
\newblock {\em Convex Analysis in General Vector Spaces}.
\newblock World Scientific Publishing, River Edge, NJ, 2002.

\end{thebibliography}
\end{small}
%%%%%%%%%%%%%%%%%%%%%%%%%%%%%%%%%%%
\end{document}